\documentclass[10pt]{article}

\usepackage{graphicx,hyperref}
\usepackage{amsmath}
\usepackage[version=4]{mhchem}
\usepackage{siunitx}
\usepackage{longtable,tabularx}
\usepackage{pstricks}
\usepackage{amsfonts}
\usepackage{latexsym}
\usepackage{epsfig}
\usepackage{subfig}
\usepackage{framed,fancybox}
\usepackage{bbm}
\usepackage{multirow}
\usepackage[T1]{fontenc}
\usepackage{palatino}
\usepackage{threeparttable}
\usepackage{rotating}
\makeatletter
\def\myVCENTER#1{\vcenter{\hbox{$\m@th#1$}}}
\makeatother

\long\def\symbolfootnote[#1]#2{\begingroup\def\thefootnote{\fnsymbol{footnote}}\footnote[#1]{#2}\endgroup}

\definecolor{shadecolor}{gray}{0.99}
{\endMakeFramed}

\definecolor{shadecolor}{gray}{0.99}
\long\def\symbolfootnote[#1]#2{\begingroup\def\thefootnote{\fnsymbol{footnote}}\footnote[#1]{#2}\endgroup}
\def\qed{\hfill{$\vcenter{\hrule height1pt \hbox{\vrule width1pt height5pt
    \kern5pt \vrule width1pt} \hrule height1pt}$} \medskip}
\newcommand{\CGPOPS}{{\mathbb{CGPOPS}}}
\newcommand{\m}[1]{{\bf{#1}}}
\newcommand{\g}[1]{\boldsymbol #1}
\newcommand{\bb}[1]{\mathbb #1}
\newcommand{\tr}{^{\sf T}}
\newcommand{\C}[1]{{\cal {#1}}}
\newcommand{\pd}[2]{\frac{\partial #1}{\partial #2}}

\title{\bf $\CGPOPS$: A C++ Software for Solving Multiple-Phase Optimal Control Problems Using Adaptive Gaussian Quadrature Collocation and Sparse Nonlinear Programming}

\author{Yunus M.~Agamawi\thanks{Ph.D.~Student, Department of Mechanical and Aerospace Engineering.  E-mail:  yagamawi@ufl.edu.} \\ Anil V.~Rao\thanks{Associate Professor, Department of Mechanical and Aerospace Engineering, Erich Farber Faculty Fellow and University Term Professor.  E-mail: anilvrao@ufl.edu.  Associate Fellow, AIAA.  Corresponding Author.} \vspace{12pt} \\ {\em University of  Florida} \\ {\em Gainesville, FL, 32611-6250}}

\date{}
\begin{document}

\renewcommand{\baselinestretch}{1}
\normalsize\normalfont 
\maketitle
\begin{abstract}
A general-purpose C++ software program called $\mathbb{CGPOPS}$ is described for solving multiple-phase optimal control problems using adaptive Gaussian quadrature collocation.  The software employs a Legendre-Gauss-Radau direct orthogonal collocation method to transcribe the continuous-time optimal control problem into a large sparse nonlinear programming problem.  A class of $hp$ mesh refinement methods are implemented which determine the number of mesh intervals and the degree of the approximating polynomial within each mesh interval to achieve a specified accuracy tolerance.  The software is interfaced with the open source Newton NLP solver IPOPT.  All derivatives required by the NLP solver are computed using either central finite differencing, bicomplex-step derivative approximation, hyper-dual derivative approximation, or automatic differentiation.  The key components of the software are described in detail and the utility of the software is demonstrated on five optimal control problems of varying complexity.  The software described in this article provides a computationally efficient and accurate approach for solving a wide variety of complex constrained optimal control problems.  
\end{abstract}

\section{Introduction\label{sect:introduction}}


Optimal control problems arise in a wide variety of subjects including virtually all branches of engineering, economics, and medicine.  Over the past few decades, the subject of optimal control has transitioned from theory to computations as a result of the increasing complexity of optimal control applications and the inability to solve them analytically.  In particular, computational optimal control has become a science in and of itself, resulting in a variety of numerical methods and corresponding software implementations of those methods.  To date, the vast majority of software implementations of optimal control involve direct transcription of a continuous optimal control problem to a nonlinear programming problem (NLP).  The resulting NLP from discretizing the continuous optimal control problem may then be solved using well established techniques.  Examples of well-known software for solving optimal control problems include {\em SOCS} \cite{Betts2}, {\em DIRCOL} \cite{vonStryk1}, {\em GESOP} \cite{Gesop1}, {\em OTIS} \cite{Vlases1}, {\em MISER} \cite{Goh1}, {\em PSOPT} \cite{Becerra1}, {\em GPOPS} \cite{Rao8}, {\em ICLOCS} \cite{ICLOCS2010}, {\em ACADO} \cite{Houska2011}, and $\mathbb{GPOPS-II}$ \cite{Patterson2014}.

Over the past few decades, direct collocation methods have become popular in the numerical solution of nonlinear optimal control problems.  In a direct collocation method, the state and control are parameterized at an appropriately chosen set of discrete points along the time interval of interest.  The continuous optimal control problem is then transcribed to a finite-dimensional NLP.  The resulting NLP may then be solved using well known software such as {\em SNOPT} \cite{Gill1}, {\em IPOPT} \cite{Biegler2}, and {\em KNITRO} \cite{Byrd1}.  Direct collocation methods were originally developed as $h$ methods (such as Euler or Runge-Kutta methods) where the time interval of interest is divided into a mesh and the state is approximated using a fixed-degree polynomial in each mesh interval.  Convergence in an $h$ method is then achieved by increasing the number and placement of the mesh points \cite{Betts3,Jain1,Zhao2}.  More recently, a great deal of research has been done in the class of direct {\em Gaussian quadrature orthogonal collocation} methods \cite{Elnagar1,Elnagar2,Benson2,Huntington1,Huntington2,Gong2,Rao8,Garg1,Garg2,Garg3,Kameswaran1,DarbyGargRao2011,Patterson2012}.  In a Gaussian quadrature collocation method, the state is typically approximated using a Lagrange polynomial where the support points of the Lagrange polynomial are chosen to be points associated with a Gaussian quadrature.  Originally, Gaussian quadrature collocation methods were implemented as $p$ methods using a single interval.  Convergence of the $p$ method was then achieved by increasing the degree of the polynomial approximation.  For problems whose solutions are smooth and well-behaved, a Gaussian quadrature orthogonal collocation method converges at an exponential rate \cite{HagerHouRao15a,HagerHouRao16a,HagerLiuMohapatraWangRao18a,DuChenHager2019,HagerHouMohapatraRaoWang2019}. The most well developed $p$ Gaussian quadrature methods are those that employ either Legendre-Gauss (LG) points \cite{Benson2,Rao8}, Legendre-Gauss-Radau (LGR) points \cite{Kameswaran1,Garg1,Garg2}, or Legendre-Gauss-Lobatto (LGL) points \cite{Elnagar1}.

In this paper a new optimal control software called $\CGPOPS$ is described that employs $hp$ direct orthogonal collocation methods.  An $hp$ method is a hybrid between an $h$ and a $p$ method in that both the number of mesh intervals and the degree of the approximating polynomial within each mesh interval can be varied in order to achieve a specified accuracy.  As a result, in an $hp$ method it is possible to take advantage of the exponential convergence of a Gaussian quadrature method in regions where the solution is smooth and introduce mesh points only when necessary to deal with potential nonsmoothness or rapidly changing behavior in the solution.  Originally, $hp$ methods were developed as finite-element methods for solving partial differential equations \cite{Babuska2,Babuska3,Babuska4,Babuska5}.  In the past few years the problem of developing $hp$ methods for solving optimal control problems has been of interest \cite{Darby3,Darby2,Patterson2015,Liu2015,Liu2018}.  The work of \cite{Darby3,Darby2,Patterson2015,Liu2015,Liu2018} provides examples of the benefits of using an $hp$-adaptive method over either a $p$ method or an $h$ method.  This recent research has shown that convergence using $hp$ methods can be achieved with a significantly smaller finite-dimensional approximation than would be required when using either an $h$ or a $p$ method.

It is noted that previously the software $\mathbb{GPOPS-II}$ was developed as described in Ref.~\cite{Patterson2014}.  Although both the $\mathbb{GPOPS-II}$ and $\CGPOPS$ software programs implement Gaussian quadrature collocation with $hp$ mesh refinement, $\CGPOPS$ is a fundamentally different software from $\mathbb{GPOPS-II}$.  First, $\mathbb{GPOPS-II}$ is a MATLAB software program, while $\CGPOPS$ is a C++ software program.  Furthermore, because $\CGPOPS$ is implemented in C++, it has the potential for improved computational efficiency and portability over a MATLAB software such as $\mathbb{GPOPS-II}$.  Second, while $\mathbb{GPOPS-II}$ employs both sparse finite-differencing and automatic differentiation using the software {\em ADiGator}~\cite{Weinstein2017}, $\CGPOPS$ employs the following four derivative estimation methods: central finite differencing, bicomplex-step~\cite{Lantoine1}, hyper-dual~\cite{Fike2011}, and automatic differentiation~\cite{Griewank0}.  Both the bicomplex-step and hyper-dual derivative approximations are referred to as semi-automatic differentiation methods and are implemented via source code transformation and operator overloading.  Third, while $\mathbb{GPOPS-II}$ is only capable of identifying the first-order derivative dependencies and over-estimates the dependencies of the second derivatives, $\CGPOPS$ is able to exactly identify both the first- and second-order derivative dependencies of the continuous optimal control problem functions when the derivatives are approximated using the hyper-dual method.  The improvement in determining the dependencies at the level of second derivatives further improves computational efficiency over $\mathbb{GPOPS-II}$.  

The objective of this paper is to describe a computationally efficient general-purpose C++ optimal control software that accurately solves a wide variety of constrained continuous optimal control problems.  In particular, the software described in this paper employs a differential form of the multiple-interval version of the {\em Legendre-Gauss-Radau (LGR) collocation method} \cite{Garg1,Garg2,Garg3,Patterson2012}.  The LGR collocation method is chosen for use in the software because it provides highly accurate state, control, and costate approximations while maintaining a relatively low-dimensional approximation of the continuous problem.  The key components of the software are then described and the software is demonstrated on five examples from the open literature.  Each example demonstrates different capabilities of the software.  The first example is the hyper-sensitive optimal control problem from Ref.~\cite{Rao4} and demonstrates the ability of the software to accurately solve a problem whose optimal solution changes rapidly in particular regions of the solution.  The second example is the reusable launch vehicle entry problem taken from Ref.~\cite{Betts3} and demonstrates the ability of $\CGPOPS$ to compute an accurate solution using a relatively coarse mesh.  The third example is the space station attitude control problem taken from Refs.~\cite{Pietz2003,Betts3} and demonstrates the ability of the software to generate accurate solutions to a problem whose solution is not intuitive.  The fourth example is a free-flying robot problem taken from Refs.~\cite{Betts3,Sakawa} and shows the ability of the software to handle bang-bang optimal control problems using the novel bang-bang control mesh refinement method included in the software.  The fifth example is a launch vehicle ascent problem taken from Refs.~\cite{Benson1,Rao8,Betts3} that demonstrates the ability of the software to solve a multiple-phase optimal control problem.  In order to validate the results, the solutions obtained using $\CGPOPS$ are compared against the solutions obtained using the software $\mathbb{CGPOPS}$ \cite{Patterson2014}.

This paper is organized as follows.  In Section \ref{sect:multi-phase} the general multiple-phase optimal control problem is present.  In Section \ref{sect:RPM} the Legendre-Gauss-Radau collocation method that is used as the basis of $\CGPOPS$ is described.  In Section \ref{sect:software} the key components of $\CGPOPS$ are described.  In Section \ref{sect:examples} the results obtained using the software on the five aforementioned examples are shown.  In Section \ref{sect:capabilities} a discussion of the capabilities of the software that are demonstrated by the results obtained using the software is provided.  In Section \ref{sect:limitations} possible limitations of the software are discussed.  Finally, in Section \ref{sect:conclusions} conclusions on the work described in this paper are provided. 


\section{General Multiple-Phase Optimal Control Problems\label{sect:multi-phase}}

\begin{figure}[ht!]
  \centering
  \epsfig{file=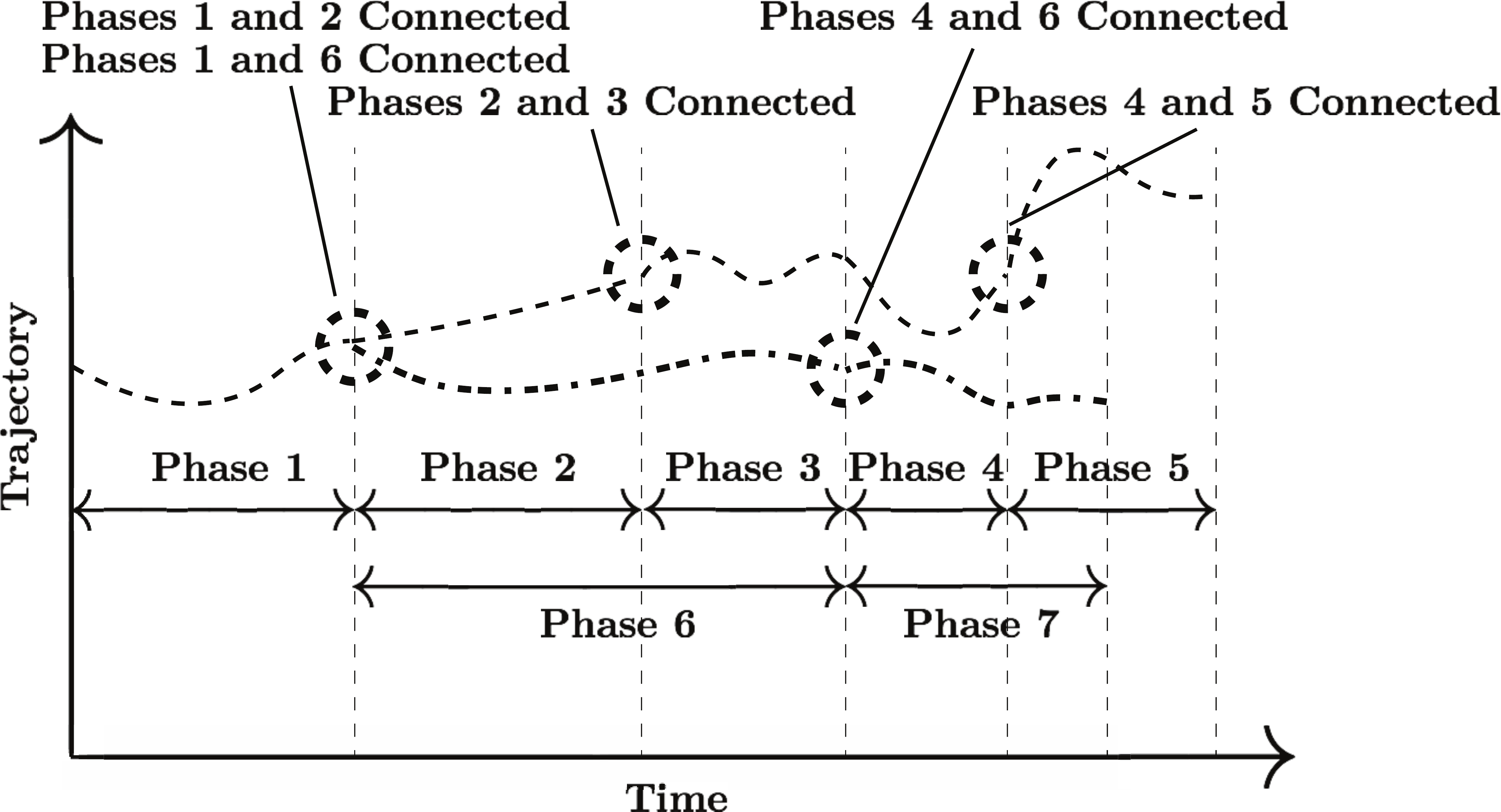,height=2.6in}
  \caption{Schematic of linkages for multiple-phase optimal control problem.  The example shown in the picture consists of seven phases where the ends of phases 1, 2, and 4 are linked to the starts of phases 2, 3, and 5, respectively, while the ends of phase 1 and 6 are linked to the starts of phase 6 and 4, respectively. \label{fig: linkages}}
\end{figure}

The general multiple-phase optimal control problem that can be solved by $\CGPOPS$ is given as follows.  Without loss of generality, consider the following general multiple-phase optimal control problem where each phase is defined on the interval $t\in[t_0^{(p)},t_f^{(p)}]$.  First let $p \in \{1,\dots,P\}$ be the phase number where $P$ is the total number of phases.  Determine the state $\m{y}^{(p)}(t)\in\bb{R}^{1~\times~n_y^{(p)}}$, the control $\m{u}^{(p)}(t)\in\bb{R}^{1~\times~n_u^{(p)}}$, the integrals $\m{q}^{(p)}\in\bb{R}^{1~\times~n_q^{(p)}}$, the start times $t_0^{(p)} \in\bb{R}$, and the terminus times $t_f^{(p)} \in\bb{R}$ in all phases $p \in \{1,\dots,P\}$, along with the static parameters $\m{s} \in\bb{R}^{1~\times~n_s}$ that minimize the objective functional
\begin{equation}\label{eq:multi-cost}
  \cal{J}=\cal{\phi}\left(\m{e}^{(1)}, \dots, \m{e}^{(P)}, \m{s}\right)~,
\end{equation}
subject to the dynamic constraints
\begin{equation}\label{eq:multi-dyn}
\frac{d\m{y}^{(p)}}{dt} \equiv \dot{\m{y}}^{(p)} = \m{a}^{(p)}\left(\m{y}^{(p)}(t),\m{u}^{(p)}(t), t, \m{s}\right)
~,\quad p \in \{1,\dots,P\}~,
\end{equation}
the event constraints
\begin{equation}\label{eq:multi-event}
\m{b}_{\min} \leq \m{b}\left(\m{e}^{(1)}, \dots, \m{e}^{(P)}, \m{s}\right) \leq \m{b}_{\max}~,
\end{equation}
the inequality path constraints
\begin{equation}\label{eq:multi-path}
\m{c}_{\min}^{(p)} \leq \m{c}^{(p)}\left(\m{y}^{(p)}(t),\m{u}^{(p)}(t), t, \m{s}\right)\leq \m{c}_{\max}^{(p)}
~,\quad p \in \{1,\dots,P\}~,
\end{equation}
the integral constraints
\begin{equation}\label{eq:multi-integral}
\m{q}_{\min}^{(p)} \leq \m{q}^{(p)} \leq \m{q}_{\max}^{(p)}
~,\quad  p \in \{1,\dots,P\}~,
\end{equation}
and the static parameter constraints
\begin{equation}\label{eq:multi-static}
\m{s}_{\min}  \leq \m{s} \leq \m{s}_{\max}~,
\end{equation}
where 
\begin{equation}\label{eq:multi-endpoint-vec}
\m{e}^{(p)} = \left[ \m{y}^{(p)}(t_0^{(p)}), t_0^{(p)}, \m{y}^{(p)}(t_f^{(p)}), t_f^{(p)}, \m{q}^{(p)} \right]
~,\quad p \in \{1,\dots,P\}~,
\end{equation}
and the integral vector components in each phase are defined as
\begin{equation}\label{eq:multi-approx}
\begin{array}{rcl}
q_j^{(p)} & = & \displaystyle \int_{t_0^{(p)}}^{t_f^{(p)}} g_j^{(p)} \left(\m{y}^{(p)}(t),\m{u}^{(p)}(t), t, \m{s}\right) dt~,
\\
 & &  j \in \{1,\dots,n_q^{(p)}\}~,~ p \in \{1,\dots,P\}~.
\end{array}
\end{equation}
It is important to note that the event constraints of Eq.~\eqref{eq:multi-event} contain functions which can relate information at the start and/or terminus of any phase (including any relationships involving any integral or static parameters), with phases not needing to be in sequential order to be linked.  Moreover, it is noted that the approach to linking phases is based on well-known formulations in the literature such as those given in Ref.~\cite{Betts1} and \cite{Betts3}.  A schematic of how phases can potentially be linked is given in Fig.~\ref{fig: linkages}.

\section{Legendre-Gauss-Radau Collocation Method\label{sect:RPM}}


As stated at the outset, the objective of this research is to provide researchers a computationally efficient general-purpose optimal control software for solving of a wide variety of complex constrained continuous optimal control problems using direct collocation.  While in principle any collocation method can be used to approximate the optimal control problem given in Section \ref{sect:multi-phase}, in this research the Legendre-Gauss-Radau (LGR) collocation method \cite{Garg1,Garg2,Garg3,Kameswaran1,Patterson2015,Liu2015,Liu2018} is employed.  It is noted that the NLP arising from the LGR collocation method has an elegant sparse structure which can be exploited as described in Refs.~\cite{Patterson2012,Patterson2014,Agamawi2018}.  In addition, the LGR collocation method has a well established convergence theory as described in Refs.~\cite{HagerHouRao15a,HagerHouRao16a,HagerLiuMohapatraWangRao18a,DuChenHager2019,HagerHouMohapatraRaoWang2019}.

In the context of this research, a multiple-interval form of the LGR collocation method is chosen.  In the multiple-interval LGR collocation method, for each phase $p$ of the optimal control problem (where the phase number $p\in\{1,\ldots,P\}$ has been omitted in order to improve clarity of the description of the method), the time interval $t\in[t_0,t_f]$ is converted into the domain $\tau\in[-1,+1]$ using the affine transformation,
\begin{equation}\label{eq:affine-transformation}
  \begin{array}{lcl}
    t & = & \displaystyle \frac{t_f-t_0}{2}\tau + \frac{t_f+t_0}{2}~, \vspace{3pt} \\
    \tau & = &  \displaystyle 2\frac{t-t_0}{t_f-t_0}-1~.
  \end{array}
\end{equation}
The domain $\tau\in[-1,+1]$ is then divided into $K$ mesh intervals, $\cal{S}_k=[T_{k-1},T_k]\subseteq [-1,+1],\;k\in\{1,\ldots,K\}$ such that
\begin{equation}
\bigcup_{k=1}^K\cal{S}_k=[-1,+1]
~,\quad\bigcap_{k=1}^K\cal{S}_k=\{T_1,\ldots,T_{K-1}\}~,
\end{equation}
and $-1=T_0<T_1<\ldots<T_{K-1}<T_K=+1$.
For each mesh interval, the LGR points used for collocation are defined in the domain of $[T_{k-1},T_k]$ for $k\in\{1,\ldots,K\}$.  The control is parameterized at the collocation points within each mesh interval.  The state of the continuous optimal control problem is then approximated in mesh interval $\cal{S}_k,\;k\in\{1,\ldots,K\}$, as 
\begin{equation}\label{eq:LGR-state-approximation}
\m{y}^{(k)}(\tau)  \approx \m{Y}^{(k)}(\tau) = \sum_{j=1}^{N_k+1} \m{Y}_{j}^{(k)}
\ell_{j}^{(k)}(\tau)~, \quad  \ell_{j}^{(k)}(\tau) = \prod_{\stackrel{l=1}{l\neq j}}^{N_k+1}\frac{\tau-\tau_{l}^{(k)}}{\tau_{j}^{(k)}-\tau_{l}^{(k)}}~, 
\end{equation}  
where
$\ell_{j}^{(k)}(\tau)$ for $ j\in\{1,\ldots,N_k+1\}$ is a basis of Lagrange polynomials,
$\left(\tau_1^{(k)},\ldots,\tau_{N_k}^{(k)}\right)$ are the set of $N_k$ Legendre-Gauss-Radau (LGR) \cite{Abramowitz1} collocation points in the interval $[T_{k-1},T_k)$ in $\cal{S}_k$, $\tau_{N_k+1}^{(k)}=T_k$ is a noncollocated support point, and $\m{Y}_{j}^{(k)} \equiv  \m{Y}^{(k)}(\tau_j^{(k)})$.  Differentiating $\m{Y}^{(k)}(\tau)$ in Eq.~(\ref{eq:LGR-state-approximation}) with respect to $\tau$ gives
\begin{equation}\label{eq:LGR-diff-state-approximation}
  \frac{d\m{Y}^{(k)}(\tau)}{d\tau} = \sum_{j=1}^{N_k+1}\m{Y}_{j}^{(k)}\frac{d\ell_j^{(k)}(\tau)}{d\tau}~.
\end{equation}
The dynamics are then approximated at the $N_k$ LGR points in mesh interval $k\in\{1,\ldots,K\}$ as
\begin{equation}\label{eq:LGR-defect}
 \sum_{j=1}^{N_k+1}D_{ij}^{(k)} \m{Y}_j^{(k)} = \frac{t_f-t_0}{2}\m{a}\left(\m{Y}_i^{(k)},\m{U}_i^{(k)},t (\tau_i^{(k)},t_0,t_f),\m{s}\right)
~,\quad i \in \{1,\dots,N_k\}~,
\end{equation}
where 
\begin{equation*}
  D_{ij}^{(k)} = \frac{d\ell_j^{(k)}(\tau_i^{(k)})}{d\tau}
  ~,\quad i \in \{1,\dots,N_k\}~,~j \in \{1,\dots,N_k+1\}~,
\end{equation*}
are the elements of the $N_k\times (N_k+1)$ {\em Legendre-Gauss-Radau differentiation matrix} \cite{Garg1} in mesh interval $\cal{S}_k$, $\;k\in\{1,\ldots,K\}$, and $\m{U}_i^{(k)}$ is the parameterized control at the $i^{th}$ collocation point in mesh interval $k$.  Finally, reintroducing the phase notation $p\in\{1,\ldots,P\}$, the phases of the problem are linked together by the event constraints
\begin{equation}\label{eq:LGR-event}
  \m{b}_{\min} \leq \m{b}\left(\m{E}^{(1)}, \dots, \m{E}^{(P)}, \m{s}\right) \leq \m{b}_{\max}~,
\end{equation}
where $\m{E}^{(p)}$ is the endpoint approximation vector for phase $p$ defined as
\begin{equation}\label{eq:NLP-endpoint-vec}
\m{E}^{(p)} = \left[ \m{Y}_{1}^{(p)}, t_0^{(p)}, \m{Y}_{N^{(p)}+1}^{(p)}, t_f^{(p)}, \m{Q}^{(p)} \right]~,
\end{equation}
such that $N^{(p)}$ is the total number of collocation points used in phase $p$ given by, 
\begin{equation}
N^{(p)} = \sum_{k=1}^{K^{(p)}} N_k^{(p)}~,
\end{equation}
and $\m{Q}^{(p)}\in\mathbb{R}^{1~\times~n_q^{(p)}}$ is the integral approximation vector in phase $p$.

The aforementioned LGR discretization then leads to
the following NLP.  Minimize the objective function
\begin{equation}\label{eq:NLP-cost}
  \cal{J}=\cal{\phi}\left(\m{E}^{(1)}, \dots, \m{E}^{(P)}, \m{s}\right)~,
\end{equation}
subject to the defect constraints
\begin{equation}\label{eq:NLP-defect}
\m{\Delta}^{(p)} =
\m{D}^{(p)}\m{Y}^{(p)} - \frac{t_f^{(p)}-t_0^{(p)}}{2}\m{A}^{(p)}
=\m{0}
~,\quad p \in \{1,\dots,P\}~,
\end{equation}
the path constraints
\begin{equation}\label{eq:NLP-path}
\m{c}_{\min}^{(p)} \leq
\m{C}_{i}^{(p)}
\leq \m{c}_{\max}^{(p)}
~,\quad i \in \{1,\dots,N^{(p)}\}, p \in \{1,\dots,P\}~,
\end{equation}
the event constraints
\begin{equation}\label{eq:NLP-event}
  \m{b}_{\min} \leq \m{b}\left(\m{E}^{(1)}, \dots, \m{E}^{(P)}, \m{s}\right) \leq \m{b}_{\max}~,
\end{equation}
the integral constraints
\begin{equation}\label{eq:NLP-integral}
\m{q}_{\min}^{(p)} \leq \m{Q}^{(p)} \leq \m{q}_{\max}^{(p)}
~,\quad p \in \{1,\dots,P\}~,
\end{equation}
the static parameter constraints
\begin{equation}\label{eq:NLP-static}
\m{s}_{\min}  \leq \m{s} \leq \m{s}_{\max}~,
\end{equation}
and integral approximation constraints
\begin{equation}\label{eq:NLP-approx}
\g{\rho}^{(p)} =
\m{Q}^{(p)} - \frac{t_f^{(p)}-t_0^{(p)}}{2}\begin{bmatrix}
\m{w}^{(p)}
\end{bmatrix}^{\tr}
\m{G}^{(p)}
= \m{0}
~,\quad p\in\{1,\dots,P\}~,
\end{equation}
where
\begin{equation}\label{eq:NLP-A-def}
\m{A}^{(p)} = \begin{bmatrix}
\m{a}^{(p)}\left(\m{Y}_{1}^{(p)},\m{U}_{1}^{(p)},t_1^{(p)},\m{s}\right)\\
\vdots\\
\m{a}^{(p)}\left(\m{Y}_{N^{(p)}}^{(p)},\m{U}_{N^{(p)}}^{(p)},t_{N^{(p)}}^{(p)},\m{s}\right)
\end{bmatrix}
\in \mathbb{R}^{N^{(p)}~\times~n_y^{(p)}}~,
\end{equation}
\begin{equation}\label{eq:NLP-C-def}
\m{C}^{(p)} = \begin{bmatrix}
\m{c}^{(p)}\left(\m{Y}_{1}^{(p)},\m{U}_{1}^{(p)},t_1^{(p)},\m{s}\right)\\
\vdots\\
\m{c}^{(p)}\left(\m{Y}_{N^{(p)}}^{(p)},\m{U}_{N^{(p)}}^{(p)},t_{N^{(p)}}^{(p)},\m{s}\right)
\end{bmatrix}
\in \mathbb{R}^{N^{(p)}~\times~n_c^{(p)}}~,
\end{equation}
\begin{equation}\label{eq:NLP-G-def}
\m{G}^{(p)} = \begin{bmatrix}
\m{g}^{(p)}\left(\m{Y}_{1}^{(p)},\m{U}_{1}^{(p)},t_1^{(p)},\m{s}\right)\\
\vdots\\
\m{g}^{(p)}\left(\m{Y}_{N^{(p)}}^{(p)},\m{U}_{N^{(p)}}^{(p)},t_{N^{(p)}}^{(p)},\m{s}\right)
\end{bmatrix}
\in \mathbb{R}^{N^{(p)}~\times~n_q^{(p)}}~,
\end{equation}
$\m{D}^{(p)} \in \mathbb{R}^{N^{(p)}~\times~[N^{(p)}+1]}$ is the LGR differentiation matrix in phase $p \in \{1,\ldots,P\}$, and $\m{w}^{(p)} \in \mathbb{R}^{N^{(p)}~\times~1}$ are the LGR weights at each node in phase $p$.  It is noted that $\m{a}^{(p)} \in \mathbb{R}^{1~\times~n_y^{(p)}}$, $\m{c}^{(p)} \in \mathbb{R}^{1~\times~n_c^{(p)}}$, and $\m{g}^{(p)} \in \mathbb{R}^{1~\times~n_q^{(p)}}$ correspond, respectively, to the functions that define the right-hand side of the dynamics, the path constraints, and the integrands in phase $p\in\{1,\ldots,P\}$, where $n_y^{(p)}$, $n_c^{(p)}$, and $n_q^{(p)}$ are, respectively, the number of state components, path constraints, and integral components in phase $p$.  Finally, the state matrix, $\m{Y}^{(p)}\in \mathbb{R}^{[N^{(p)}+1]~\times~n_y^{(p)}}$, and the control matrix, $\m{U}^{(p)}\in \mathbb{R}^{N^{(p)}~\times~n_u^{(p)}}$, in phase $p\in\{1,\ldots,P\}$ are formed as 
\begin{equation}\label{eq:NLP-Var-Mat}
\m{Y}^{(p)} = \begin{bmatrix}
\m{Y}^{(p)}_{1}\\
\vdots \\
\m{Y}^{(p)}_{N^{(p)}+1}
\end{bmatrix}
\text{ and }
\m{U}^{(p)} = 
\begin{bmatrix}
\m{U}^{(p)}_{1}  \\
\vdots	\\
\m{U}^{(p)}_{N^{(p)}}
\end{bmatrix}~,
\end{equation}
respectively, where $n_u^{(p)}$ is the number of control components in phase $p$.


\section{Major Components of $\mathbb{CGPOPS}$\label{sect:software}}


In this section we describe the major components of the C++ software $\CGPOPS$ that implements the aforementioned LGR collocation method.  In Section \ref{subsect:nlp}, the large sparse nonlinear programming problem (NLP) associated with the LGR collocation method is described.  In Section \ref{subsect:sparse-structure}, the structure of the NLP described in Section \ref{subsect:nlp} is shown.  In Section \ref{subsect:scaling} the method for scaling the NLP via scaling of the optimal control problem is over-viewed.  In Section \ref{subsect:derivative-computation}, the approach for estimating the derivatives required by the NLP solver is explained.  In Section \ref{subsect:dependencies}, the method for determining the dependencies of each optimal control function in order to provide the most sparse NLP to the NLP solver is presented.  In Section \ref{subsect:mesh-refinement} the $hp$ mesh refinement methods that are included in the software in order to iteratively determine a mesh that satisfies a user-specified accuracy tolerance are described.  Finally, in Section \ref{subsect:flow} we provide a high level description of the algorithmic flow of $\CGPOPS$.  


\subsection{Sparse NLP Arising from Radau Collocation Method\label{subsect:nlp}}


The resulting nonlinear programming problem (NLP) that arises when using LGR collocation to discretize the continuous optimal control problem is given as follows.  Determine the NLP decision vector, $\m{z}$, that minimizes the NLP objective function,
\begin{equation}\label{eq:NLP-objective}
  f(\m{z})~,
\end{equation}
subject to the constraints
\begin{equation}\label{eq:NLP-constraint}
  \m{H}_{\min} \leq \m{H}(\m{z}) \leq \m{H}_{\max}~,
\end{equation}
and the variable bounds
\begin{equation}\label{eq:NLP-variable-bound}
  \m{z}_{\min} \leq \m{z} \leq \m{z}_{\max}~.
\end{equation}
It is noted that, while the size of the NLP arising from the LGR collocation method changes depending upon the number of mesh intervals and LGR points used in each phase, the structure of the NLP remains the same regardless of the size of the NLP.  Finally, in the sections that follow, the subscript "$:$" denotes either a row or a column, where the ``$:$'' notation is analogous to the syntax used in the MATLAB programming language.

\subsubsection{NLP Variables}\label{subsubsect:nlp-variables}

For a continuous optimal control problem transcribed into $P$ phases, the NLP decision vector, $\m{z}$, has the following form:
\begin{equation}\label{eq:NLP-Dec-Vec}
\m{z} = \begin{bmatrix}
\m{z}^{(1)}\\
\vdots\\
\m{z}^{(P)}\\
s_1\\
\vdots\\
s_{n_s}
\end{bmatrix}
~,\quad
\text{where }
\m{z}^{(p)} = \begin{bmatrix}
\m{Y}_{(:,1)}^{(p)}\\
\vdots\\
\m{Y}_{(:,n_y^{(p)})}^{(p)} \vspace{3pt}\\
\m{U}_{(:,1)}^{(p)}\\
\vdots\\
\m{U}_{(:,n_u^{(p)})}^{(p)} \vspace{3pt}\\
(\m{Q}^{(p)})^{\tr}\\
t_0^{(p)} \\
t_f^{(p)}\\
\end{bmatrix}~,
\end{equation}
$\m{Y}^{(p)} \in \mathbb{R}^{N^{(p)}~\times~n_y^{(p)}}$ is the state approximation matrix [see Eq.~\eqref{eq:NLP-Var-Mat}], $\m{U}^{(p)} \in \mathbb{R}^{N^{(p)}~\times~n_u^{(p)}}$ is the control parameterization matrix [see Eq.~\eqref{eq:NLP-Var-Mat}], $\m{Q}^{(p)} \in \mathbb{R}^{1~\times~n_q^{(p)}}$ is the integral approximation vector, and $t_0^{(p)}$ and $t_f^{(p)}$ are scalars of the initial and final time, respectively, for phase $p \in \{1,\ldots,P\}$, and $s_i$ for $i \in \{1,\dots,n_s\}$ are the static parameters appearing throughout the entire problem.

\subsubsection{NLP Objective and Constraint Functions}\label{subsubsect:nlp-objective-constraints}

The NLP objective function, $f(\m{z})$, is given in the form
\begin{equation}\label{eq:NLP-Obj-Func}
f(\m{z}) = \phi^{(p)}(\m{E}^{(1)},\ldots,\m{E}^{(P)},\m{s})~,
\end{equation}
where $\m{E}^{(p)},~p\in\{1,\ldots,P\}$, is the endpoint approximation vector defined in Eq. \eqref{eq:NLP-endpoint-vec}, and the typical cost functional of a general multiple-phase optimal control problem has been turned simply into a Mayer cost function by using the integral variables, $\m{Q}^{(p)}$, to approximate the Lagrange cost in each phase~$p$.  The NLP constraint vector, $\m{H}(\m{z})$, is given in the form
\begin{equation}\label{eq:NLP-Con-Vec}
\m{H}(\m{z}) = \begin{bmatrix}
\m{h}^{(1)}\\
\vdots\\
\m{h}^{(P)}\\
\m{b}
\end{bmatrix}
~,~
\text{where }
\m{h}^{(p)} = \begin{bmatrix}
\m{\Delta}_{(:,1)}^{(p)}\\
\vdots\\
\m{\Delta}_{(:,n_y^{(p)})}^{(p)} \vspace{3pt} \\
\m{C}_{(:,1)}^{(p)}\\
\vdots\\
\m{C}_{(:,n_c^{(p)})}^{(p)} \\
(\g{\rho}^{(p)})^{\tr}\\
\end{bmatrix}~,\quad p = \{1,\dots,P\}~,
\end{equation}
$\m{\Delta}^{(p)} \in \mathbb{R}^{N^{(p)}~\times~n_y^{(p)}}$, $\g{\rho}^{(p)}\in\mathbb{R}^{1~\times~n_q^{(p)}}$, and $\m{C}^{(p)} \in \mathbb{R}^{N^{(p)}~\times~n_c^{(p)}}$, are, respectively, the defect constraint matrix, the integral approximation constraint vector, and the path constraint matrix in phase $p \in \{1,\dots,P \}$, and $\m{b} \in \mathbb{R}^{n_b~\times~1}$ is the event constraint vector for the entire problem.  The defect constraint matrix, integral approximation constraint vector, and path constraint matrix in phase $p$ are defined by Eqs.~\eqref{eq:NLP-defect}, \eqref{eq:NLP-approx}, and \eqref{eq:NLP-C-def}, respectively.  It is noted that the constraints are divided into the equality defect and integral constraints
\begin{equation}\label{eq:NLP-defect-integral-bounds}
\begin{array}{rcl}
\g{\Delta}^{(p)} & = & \m{0}~,\\
\g{\rho}^{(p)} & = & \m{0}~,
\end{array}
\quad p\in\{1,\dots,P\}~,
\end{equation}
and the inequality discretized path and event constraints
\begin{equation}
\begin{array}{rcccll}
\m{c}^{(p)}_{\min} & \leq & \m{C}_i^{(p)} & \leq & \m{c}^{(p)}_{\max}~, & i\in\{1,\dots,n_c^{(p)}\}~,~p\in\{1,\dots,P\}~, \\
\m{b}_{\min} & \leq & \m{b} & \leq & \m{b}_{\max}~.
\end{array}
\end{equation}


\subsection{Sparse Structure of NLP Derivative Functions}\label{subsect:sparse-structure}


The structure of the NLP created by the LGR collocation method has been described in detail in Refs.~\cite{Patterson2012} and \cite{Agamawi2018}.  Specifically, Refs.~\cite{Patterson2012} and \cite{Agamawi2018} describe the sparse structure of the NLP for the differential form of the LGR collocation method for the single and multiple phase optimal control problem, respectively.  As described in Section~\ref{subsubsect:nlp-variables}, the values of the state approximation coefficients at the discretization points, the control parameters at the collocation points, the initial time, the final time, and the integral vector of each phase, as well as any static parameters of the problem make up the NLP decision vector.  The NLP constraints vector consists of the defect constraints and path constraints applied at each of the collocation points, as well as any integral approximation constraints, for each phase, and event constraints, as described in Section~\ref{subsubsect:nlp-objective-constraints}.  The derivation of the NLP derivative matrices in terms of the original continuous optimal control problem functions is described in detail in Refs.~\cite{Patterson2012,Patterson2014,Agamawi2018} and is beyond the scope of this paper.  It is noted that the sparsity exploitation derived in Refs.~\cite{Patterson2012,Patterson2014,Agamawi2018} requires computing partial derivatives of the continuous optimal control problem functions on the first- and second-order derivative levels.

\begin{figure}
\centering

\subfloat[NLP Constraint Jacobian]{\hspace{-0.1in}\includegraphics[height=3.1in]{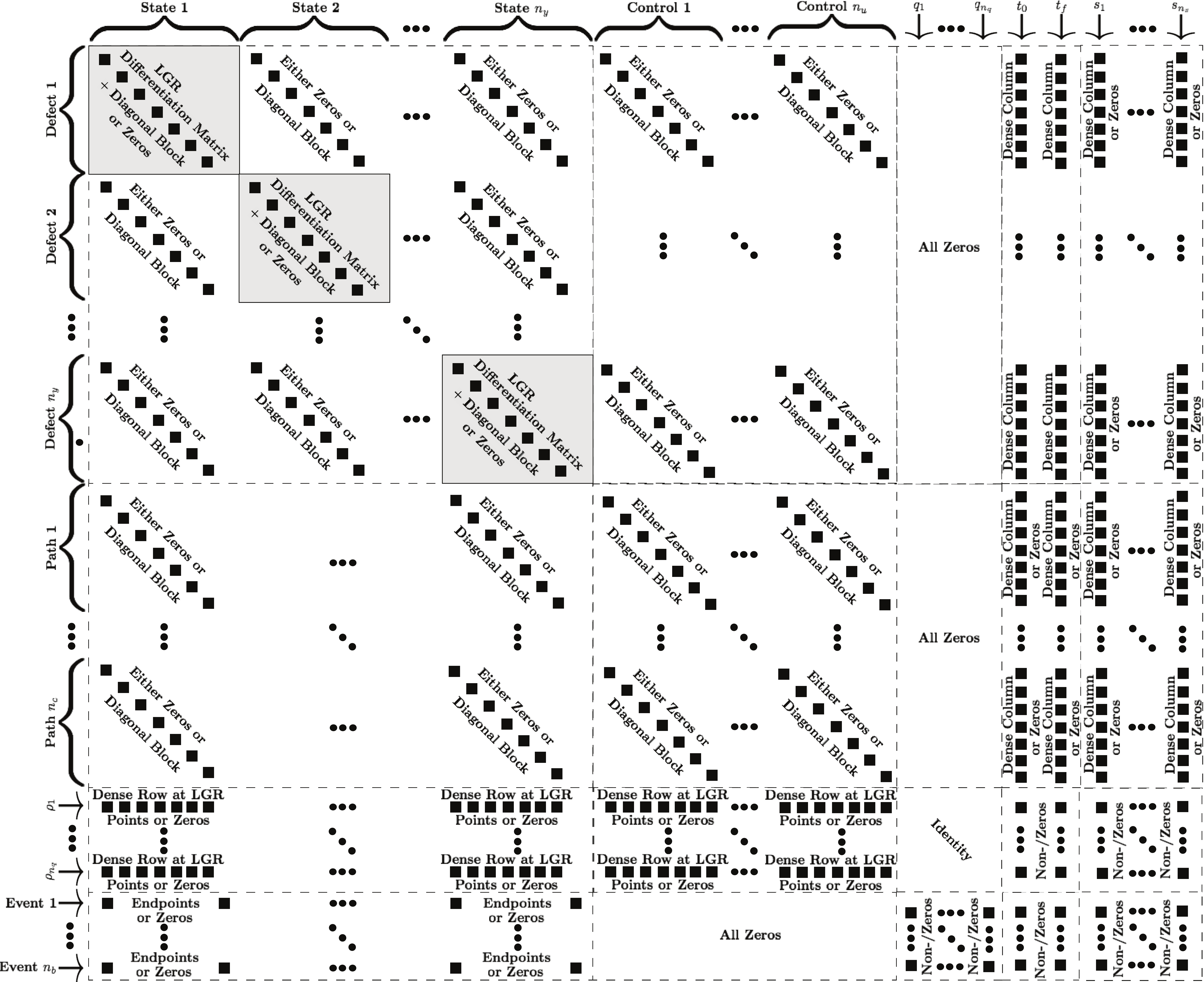}\label{fig:ExampleNLPConJac}}

\subfloat[NLP Lagrangian Hessian]{\includegraphics[height=3.7in]{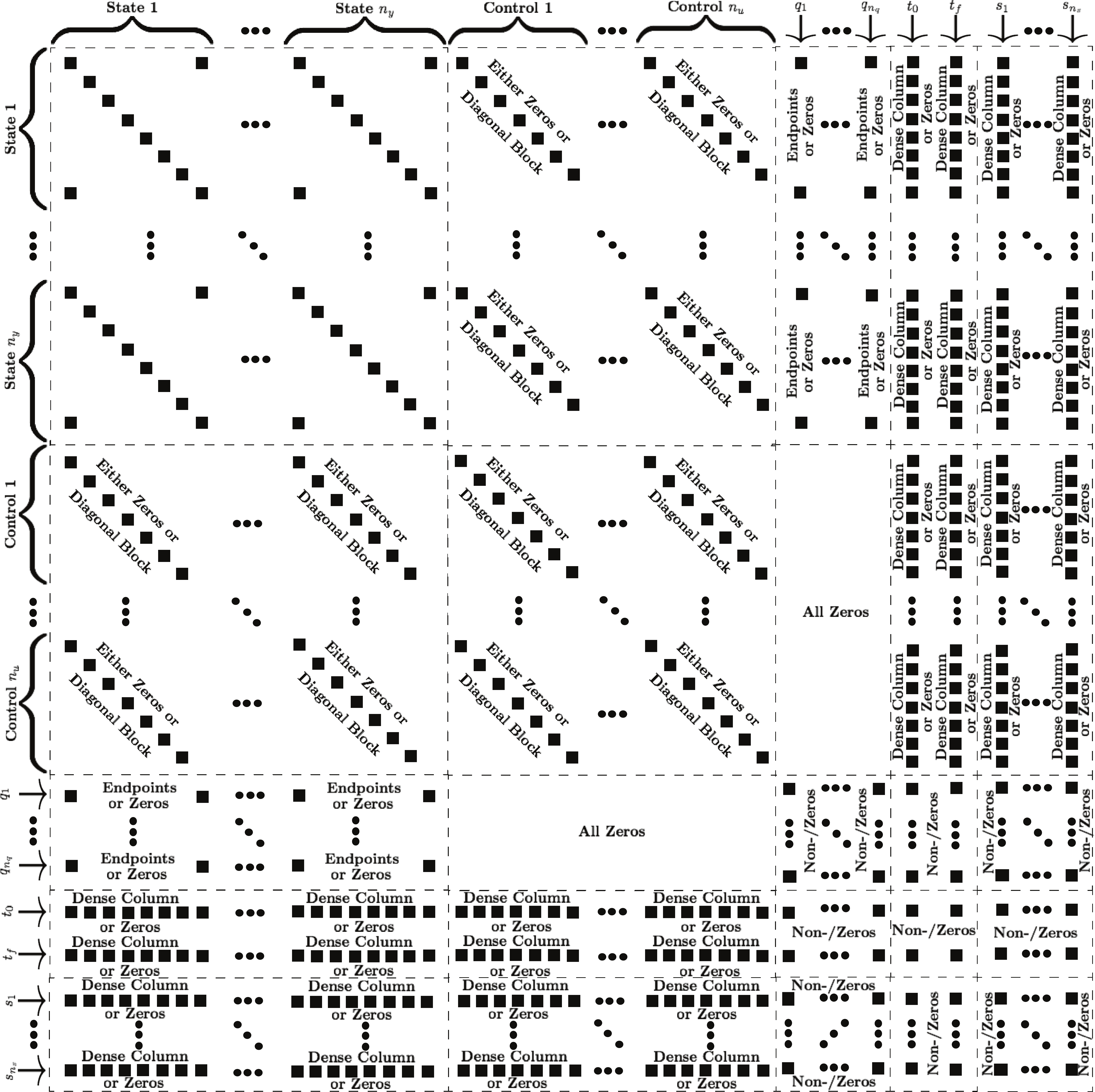}\label{fig:ExampleNLPLagHess}}
  
\renewcommand{\baselinestretch}{1}\normalsize\normalfont
\caption{Example Sparsity Patterns for Single Phase Optimal Control Problem Containing $n_y$ State Components, $n_u$ Control Components, $n_q$ Integral Components, and $n_c$ Path Constraints, $n_s$ Static Parameters, and $n_b$ Event Constraints. \label{fig:ExampleNLPSparsePat}}

\end{figure}

Examples of the sparsity patterns of the NLP constraint Jacobian and Lagrangian Hessian are shown, respectively, in Figs.~\ref{fig:ExampleNLPConJac} and \ref{fig:ExampleNLPLagHess} for a single-phase optimal control problem.  It is noted that for the NLP constraint Jacobian, all of the off-diagonal phase blocks relating constraints in phase $i$ to variables in phase $j$ for $i\neq j$ are all zeros.  Similarly, for the NLP Lagrangian Hessian, all of the off-diagonal phase blocks relating variables in phase $i$ to variables in phase $j$ for $i\neq j$ are all zeros except for the variables making up the endpoint vectors which may be related via the objective function or event constraints.  The sparsity patterns shown in Fig.~\ref{fig:ExampleNLPSparsePat} are determined explicitly by identifying the derivative dependencies of the NLP objective and constraints functions with respect to the NLP decision vector variables.  It is noted that the phases are connected using the initial and terminal values of the time and state in each phase along with the static parameters.

\subsection{Scaling of Optimal Control Problem for NLP\label{subsect:scaling}}


The NLP described in Section~\ref{subsect:nlp} must be well scaled in order for the NLP solver to obtain a solution.  $\CGPOPS$ includes the option for the NLP to be scaled automatically by scaling the continuous optimal control problem.  The approach to automatic scaling is to scale the variables and the first derivatives of the optimal control functions to be $\approx\C{O}(1)$.  First, the optimal control variables are scaled to lie on the unit interval $[-1/2, 1/2]$ and is accomplished as follows.  Suppose it is desired to scale an arbitrary variable $x\in [a,b]$ to $\tilde{x}$ such that $\tilde{x}\in[-1/2,1/2]$.  This variable scaling is accomplished via the affine transformation
\begin{equation}\label{variable-scaling}
 \tilde{x} = v_xx + r_x,
\end{equation}
where $v_x$ and $r_x$ are the variable scale and shift, respectively, defined as
\begin{equation}\label{variable-scale-and-shift}
  \begin{array}{rcl}
    v_x & = & \displaystyle \frac{1}{b-a}~,\\
    r_x & = & \displaystyle \frac{1}{2} - \frac{b}{b-a}~.
  \end{array}
\end{equation}
Every variable in the continuous optimal control problem is scaled using Eqs.~\eqref{variable-scaling} and \eqref{variable-scale-and-shift}.  Next, the Jacobian of the NLP constraints can be made $\approx\C{O}(1)$ by scaling the derivatives of the optimal control functions to be approximately unity.  First, using the approach derived in \cite{Betts3}, in $\CGPOPS$ the defect constraints are scaled using the same scale factors as were used to scale the state.  Next, the objective function, event constraints, and path constraints scale factors are obtained by sampling the gradient of each constraint at a variety of sample points within the bounds of the unscaled optimal control problem and taking the average norm of each gradient across all sample points.


\subsection{Computation Derivatives of NLP Functions\label{subsect:derivative-computation}}

The NLP derivative functions are obtained by exploiting the sparse structure of the NLP arising from the $hp$ LGR collocation method.  Specifically, in Refs.~\cite{Patterson2012,Agamawi2018} it has been shown that by using the derivative form of the LGR collocation method, the NLP derivatives can be obtained by computing the derivatives of the optimal control problem functions at the LGR points and inserting these derivatives into the appropriate locations in the NLP derivative functions.  In $\CGPOPS$, the optimal control derivative functions are approximated using one of four types of derivative estimation methods: sparse central finite-differencing, bicomplex-step derivative approximations, hyper-dual derivative approximations, and automatic differentiation.

\subsubsection{Central Finite Difference}\label{subsubsect:central-finite-difference}

To see how the central finite-difference derivative approximation works in practice, consider the function $\m{f}(\m{x})$, where $\m{f}:\mathbb{R}^n\rightarrow\mathbb{R}^m$ is one of the {\em optimal control functions} (that is, $n$ and $m$ are, respectively, the size of an optimal control variable and an optimal control function).  Then $\partial\m{f}/\partial\m{x}$ is approximated using a central finite difference as
\begin{equation}\label{eq:Central-Finite-Difference-First-Order}
\pd{\m{f}}{x_i} \approx \frac{\m{f}(\m{x}+\m{h}_i)-\m{f}(\m{x}-\m{h}_i)}{2h},
\end{equation}
where $\m{h}_i$ arises from perturbing the $i^{th}$ component of $\m{x}$.
The vector $\m{h}_i$ is computed as
\begin{equation}\label{eq:vector-perturbation}
  \m{h}_i = h_i\m{e}_i
\end{equation}
where $\m{e}_i$ is the $i^{th}$ row of the $n\times n$ identity matrix and $h_i$ is the perturbation size associated with $x_i$.  The perturbation $h_i$ is computed using the equation
\begin{equation}\label{eq:perturbation-size}
  h_i = h(1+|x_i|), 
\end{equation}
where the base perturbation size $h$ is chosen to be the optimal step size for a function whose input and output are $\approx\C{O}(1)$ as described in Ref.~\cite{Gill3}.  Second derivative approximations are computed in a manner similar to that used for first derivative approximations with the key difference being that perturbations in two variables are performed.  For example, $\partial^2\m{f}/\partial x_i\partial x_j$ can be approximated using a central finite-difference approximation as
\begin{equation}\label{eq:Central-Finite-Difference-Second-Order-Mixed}
\frac{\partial^2 \m{f}(\m{x})}{\partial x_i \partial x_j} \approx \frac{\m{f}(\m{x}+\m{h}_i+\m{h}_j)-\m{f}(\m{x}+\m{h}_i-\m{h}_j)-\m{f}(\m{x}-\m{h}_i+\m{h}_j)+\m{f}(\m{x}-\m{h}_i-\m{h}_j)}{4h_ih_j},
\end{equation}
where $\m{h}_i$, $\m{h}_j$, $h_i$, and $h_j$ are as defined in Eqs.~\eqref{eq:vector-perturbation} and \eqref{eq:perturbation-size}.  The base perturbation size is chosen to minimize round-off error in the finite-difference approximation.  Furthermore, it is noted that $h_i\rightarrow h$ as $|x_i|\rightarrow 0$.  

\subsubsection{Bicomplex-step}\label{subsubsect:bicomplex-step}

To see how the bicomplex-step derivative approximation works in practice, consider the function $\m{f}(\m{x})$, where $\m{f}:\mathbb{R}^n\rightarrow\mathbb{R}^m$ is one of the {\em optimal control functions} (that is, $n$ and $m$ are, respectively, the size of an optimal control variable and an optimal control function).  Then $\partial\m{f}/\partial\m{x}$ is approximated using a bicomplex-step derivative approximation as
\begin{equation}\label{eq:Bicomplex-Step-First-Order-Derivative-Approximation-I1}
\frac{\partial \m{f}(\m{x})}{\partial x_i} \approx \frac{\text{Im}_1\left[f(\m{x}+i_1h\m{e}_i)\right]}{h},
\end{equation}
where $\text{Im}_{1}[\cdot]$ denotes the imaginary $i_1$ component of the function evaluated with the perturbed bicomplex input, $\m{e}_i$ is the $i^{th}$ row of the $n\times n$ identity matrix, and the base perturbation size $h$ is chosen to be a step size that will minimize truncation error while refraining from encountering roundoff error due to bicomplex arithmetic, which is described in detail in Ref.~\cite{Lantoine1}, and is beyond the scope of this paper.  It is noted that the imaginary component $i_1$ has the property $i_1^2 = -1$.  Second derivative approximations are computed in a manner similar to that used for first derivative approximations with the key difference being that perturbations in two variables are performed in two separate imaginary directions.  For example, $\partial^2\m{f}/\partial x_i\partial x_j$ can be approximated using a bicomplex-step derivative approximation as
\begin{equation}\label{eq:Bicomplex-Step-Second-Order-Derivative-Approximation}
\frac{\partial^2 \m{f}(\m{x})}{\partial x_i \partial x_j} \approx \frac{\text{Im}_{1,2}\left[f(\m{x}+i_1h\m{e}_i+i_2h\m{e}_j)\right]}{h^2},
\end{equation}
where $\text{Im}_{1,2}[\cdot]$ denotes the imaginary $i_1i_2$ component of the function evaluated with the perturbed bicomplex input, where it is noted that $i_2^2=-1$, and $i_1i_2$ is a bi-imaginary direction distinct from either the $i_1$ or $i_2$ imaginary directions (i.e. $i_1i_2=i_2i_1$).

\subsubsection{Hyper-Dual}\label{subsubsect:hyper-dual}

To see how the hyper-dual derivative approximation works in practice, consider the function $\m{f}(\m{x})$, where $\m{f}:\mathbb{R}^n\rightarrow\mathbb{R}^m$ is one of the {\em optimal control functions} (that is, $n$ and $m$ are, respectively, the size of an optimal control variable and an optimal control function).  Then $\partial\m{f}/\partial\m{x}$ is approximated using a hyper-dual derivative approximation as
\begin{equation}\label{eq:Hyper-Dual-First-Order-Derivative-Approximation-I1}
\frac{\partial \m{f}(\m{x})}{\partial x_i} = \frac{\text{Ep}_1\left[f(\m{x}+\epsilon_1h\m{e}_i)\right]}{h},
\end{equation}
where $\text{Ep}_{1}[\cdot]$ denotes the imaginary $\epsilon_1$ component of the function evaluated with the perturbed hyper-dual input, $\m{e}_i$ is the $i^{th}$ row of the $n\times n$ identity matrix, and the base perturbation size $h$ is chosen to be unity because for first- and second-derivatives the hyper-dual arithmetic does not suffer from either truncation or roundoff error (described in detail in Ref.~\cite{Fike2011} and beyond the scope of this paper).  It is noted that the imaginary component $\epsilon_1$ has the property of being nilpotent (that is, $\epsilon_1^2 = 0$).  Second derivative approximations are computed in a manner similar to that used for first derivative approximations with the key difference being that perturbations in two variables are performed in two separate imaginary directions.  For example, $\partial^2\m{f}/\partial x_i\partial x_j$ can be approximated using a hyper-dual derivative approximation as
\begin{equation}\label{eq:Hyper-Dual-Second-Order-Derivative-Approximation}
\frac{\partial^2 \m{f}(\m{x})}{\partial x_i \partial x_j} = \frac{\text{Ep}_{1,2}\left[f(\m{x}+\epsilon_1h\m{e}_i+\epsilon_2h\m{e}_j)\right]}{h^2},
\end{equation}
where $\text{Ep}_{1,2}[\cdot]$ denotes the imaginary $\epsilon_1\epsilon_2$ component of the function evaluated with the perturbed hyper-dual input, where it is noted that $\epsilon_2$ also has the property of being nilpotent (i.e. $\epsilon_2^2=0$), and $\epsilon_1\epsilon_2$ is a bi-imaginary direction distinct from either the $\epsilon_1$ or $\epsilon_2$ imaginary directions (i.e. $\epsilon_1\epsilon_2=\epsilon_2\epsilon_1$).


\subsubsection{Automatic Differentiation \label{subsubsect:AlgorithmicDiff}}


In this section, the basis of automatic differentiation is discussed.  As described in Ref.~\cite{Martins2}, automatic (algorithmic) differentiation may be derived from the unifying chain rule and supplies numerical evaluations of the derivative for a defined computer program by decomposing the program into a sequence of elementary function operations and applying the calculus chain rule algorithmically through the computer \cite{Griewank0}.  The process of automatic differentiation is described in detail in Ref.~\cite{Griewank0}, and is beyond the scope of this paper.  It is noted, however, that the first- and second-order partial derivatives obtained using the Taylor series-based derivative approximation methods described in Sections~\ref{subsubsect:central-finite-difference} -- \ref{subsubsect:hyper-dual} may be computed to machine precision using automatic differentiation.  Specifically, $\CGPOPS$ employs the well-known open source software ADOL-C \cite{Griewank1,Walther1} to compute derivatives using automatic differentiation.


\subsection{Method for Determining the Optimal Control Function Dependencies\label{subsect:dependencies}}


It can be seen from Section~\ref{subsect:sparse-structure} that the NLP associated with the LGR collocation method has a sparse structure where the blocks of the constraint Jacobian and Lagrangian Hessian are dependent upon whether a particular NLP function depends upon a particular NLP variable, as was shown in Refs.~\cite{Patterson2012,Agamawi2018}.  The method for identifying the optimal control function derivative dependencies in $\CGPOPS$ utilizes the independent nature of the hyper-dual derivative approximations.  Specifically, since the imaginary directions used for hyper-dual derivative approximations are completely independent of one another, second-order derivative approximations only appear nonzero if the partial actually exists (same for first-order derivative approximations).  For example, suppose that $\m{f}(\m{x})$ is a function where $\m{f}:\mathbb{R}^n\rightarrow\mathbb{R}^m$ and $\m{x}=\left[ x_1 \ldots x_n\right]$.  The hyper-dual derivative approximation of $\partial^2 \m{f}(\m{x})/\partial x_i \partial x_j$ will only be nonzero if the actual $\partial^2 \m{f}(\m{x})/\partial x_i \partial x_j$ exists and is nonzero.  Given this knowledge of the exact correspondence of hyper-dual derivative approximations to the actual derivative evaluations, identifying derivative dependencies of optimal control problem functions with respect to optimal control problem variables becomes simple, as existing partial derivatives will have nonzero outputs when approximated by the hyper-dual derivative approximations, while non-existing partial derivatives will simply be zero always.  In order to ensure that derivative dependencies aren't mistakenly missed due to a derivative approximation happening to equal zero at the point its being evaluated at for an existing nonzero partial derivative, the hyper-dual derivative approximations are evaluated at multiple sample points within the variable bounds.  In this manner, the derivative dependencies of the optimal control problem functions can be easily identified exactly for the first- and second-order derivative levels.  The computational expense of identifying the derivative dependencies in this manner is minimal, while the exact second-order derivative sparsity pattern that is obtained can significantly reduce the cost of computing the NLP Lagrangian Hessian when compared to using an over-estimated sparsity pattern as done in $\mathbb{GPOPS-II}$ \cite{Patterson2014}.


\subsection{Adaptive Mesh Refinement\label{subsect:mesh-refinement}}


In the past few years, the subject of adaptive mesh refinement has been of considerable study in the efficient implementation of Gaussian quadrature collocation methods.  The work on adaptive Gaussian quadrature mesh refinement has led to several articles in the literature including those found in Refs.~\cite{Gong3,Darby2,Darby3,Patterson2015,Liu2015,Liu2018}.  $\CGPOPS$ employs the recently developed mesh refinement methods described in \cite{Darby2,Patterson2015,Liu2015,Liu2018,Agamawi2019}.  The mesh refinement methods of Refs.~\cite{Patterson2015}, \cite{Darby2}, \cite{Liu2015}, \cite{Liu2018}, and \cite{Agamawi2019} are referred to, respectively, as the $hp$-I, $hp$-II , $hp$-III, $hp$-IV, and $hp$-BB methods.  In all five of the $hp$-adaptive mesh refinement methods, the number of mesh intervals, width of each mesh interval, and the degree of the approximating polynomial can be varied until a user-specified accuracy tolerance has been achieved.  When using any of the methods in $\CGPOPS$, the terminology $hp$-Method$(N_{\min},N_{\max})$ refers to a method whose minimum and maximum allowable polynomial degrees within a mesh interval are $N_{\min}$ and $N_{\max}$, respectively.  All five methods estimate the solution error using a relative difference between the state estimate and the integral of the dynamics at a modified set of LGR points.  The key difference between the five methods lies in the manner in which the decision is made to either increase the number of collocation points in a mesh interval or to refine the mesh.  In Ref.~\cite{Darby2} the degree of the approximating polynomial is increased if the ratio of the maximum curvature over the mean curvature of the state in a particular mesh interval is below a user-specified threshold.  On the other hand, Ref.~\cite{Patterson2015} uses the exponential convergence property of the LGR collocation method and increases the polynomial degree within a mesh interval if the estimate of the required polynomial degree is less than a user-specified upper limit.  Similarly, Refs.~\cite{Liu2015} and \cite{Liu2018} employ nonsmoothness criterion to determine whether an $h$ or $p$ method should be used for a given mesh interval, while also utilizing mesh reduction techniques in order to minimize the size of the transcribed NLP in regions of the solution where such high resolution is not required.  If a $p$ method refinement is prescribed for a given mesh interval and the estimate of the polynomial degree exceeds the allowed upper limit, the mesh interval is divided into more mesh intervals (i.e. $h$ method employed).  Lastly, the mesh refinement method developed in Ref.~\cite{Agamawi2019} is designed for bang-bang optimal control problems and employs estimates of the switching functions of the Hamiltonian in order to obtain the solution profile.  In $\CGPOPS$, the user can choose between these five mesh refinement methods.  Finally, it is noted that $\CGPOPS$ has been designed in a modular way, making it possible to add a new mesh refinement method in a relatively straightforward way if it is so desired.


\subsection{Algorithmic Flow of $\mathbb{CGPOPS}$\label{subsect:flow}} 

\begin{figure}
 \centering
 \epsfig{file=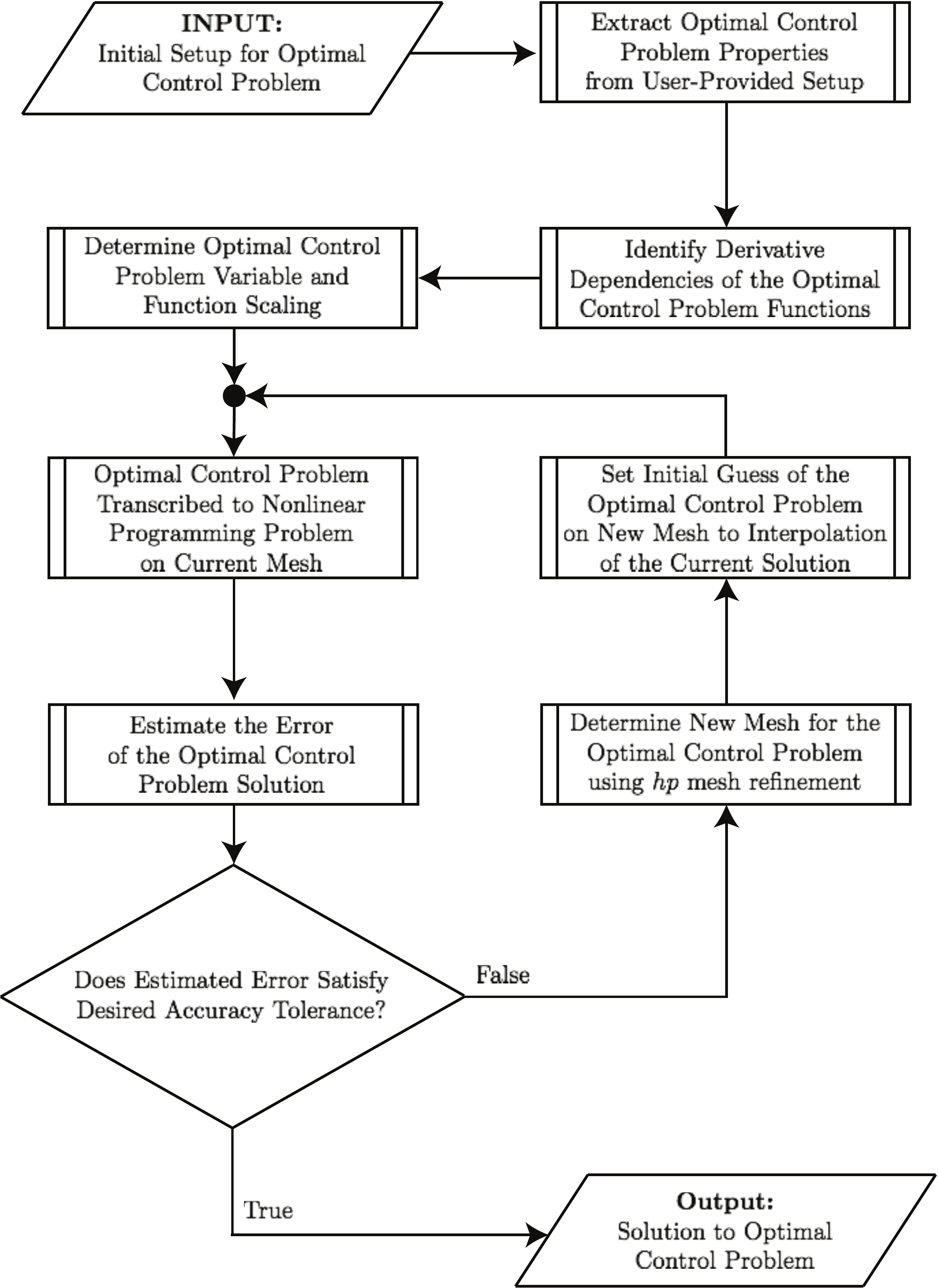,height=4in}
 \caption{Flowchart of the $\CGPOPS$ Algorithm. \label{fig:algorithmic_flow}}
\end{figure}

In this section we describe the operational flow of $\CGPOPS$ with the aid of Fig.~\ref{fig:algorithmic_flow}.  First, the user provides a description of the optimal control problem that is to be solved.  The properties of the optimal control problem are then extracted from the user description from which the state, control, time, and parameter dependencies of the optimal control problem functions are identified.  Subsequently, assuming that the user has specified that the optimal control problem be scaled automatically, the optimal control problem scaling algorithm is called and these scale factors are determined and used to scale the NLP.  The optimal control problem is then transcribed to a large sparse NLP and the NLP is solved on the initial mesh, where the initial mesh is either user-supplied or is determined by the default settings in $\CGPOPS$.  Once the NLP is solved, the NLP solution is analyzed as a discrete approximation of the optimal control problem and the error in the discrete approximation for the current mesh is estimated.  If the user-specified accuracy tolerance is met, the software terminates and outputs the solution.  Otherwise, a new mesh is determined using one of the supplied mesh refinement algorithms and the resulting NLP is solved on the new mesh.  


\section{Examples\label{sect:examples}}


$\CGPOPS$ is now demonstrated on five examples taken from the open literature.  The first example is the hyper-sensitive optimal control problem from Ref.~\cite{Rao4} and demonstrates the ability of $\CGPOPS$ to efficiently solve problems that have rapid changes in dynamics in particular regions of the solution.  The second example is the reusable launch vehicle entry problem taken from Ref.~\cite{Betts3} and demonstrates the efficiency of $\CGPOPS$ on a more realistic problem.  The third example is the space station attitude optimal control problem taken from Refs.~\cite{Pietz2003} and \cite{Betts3} and demonstrates the efficiency of $\CGPOPS$ on a problem whose solution is highly non-intuitive.  The fourth example is a free-flying robot problem taken from Ref.~\cite{Betts3} and demonstrates the ability of $\CGPOPS$ to solve a bang-bang optimal control problem using discontinuity detection.  The fifth example is a multiple-stage launch vehicle ascent problem taken from Refs.~\cite{Benson1,Rao8,Betts3} and demonstrates the ability of $\CGPOPS$ to solve a problem with multiple-phases.

All five examples were solved using the open-source NLP solver IPOPT \cite{Biegler2} in second derivative (full Newton) mode with the publicly available multifrontal massively parallel sparse direct linear solver MA57 \cite{Duff1}.  All results were obtained using the differential form of the LGR collocation method and various forms of the aforementioned $hp$ mesh refinement method using default NLP solver settings and the automatic scaling routine in $\CGPOPS$.  For $\CGPOPS$, all first- and second-order derivatives for the NLP solver were obtained using hyper-dual derivative approximations as described in Section~\ref{subsubsect:hyper-dual} with a perturbation step size of $h=1$.  All solutions obtained by $\CGPOPS$ are compared against the solutions obtained using the previously developed MATLAB software $\mathbb{GPOPS-II}$ \cite{Patterson2014} which also employs $hp$ LGR collocation methods.  For $\mathbb{GPOPS-II}$, the first- and second-order derivatives for the NLP solver were obtained using the automatic differentiation software ADiGator \cite{Weinstein2014} for all examples except the fifth example which used sparse central finite-differences.  All computations were performed on a 2.9 GHz Intel Core i7 MacBook Pro running MAC OS-X version 10.13.6 (High Sierra) with 16GB 2133MHz LPDDR3 of RAM.  C++ files were compiled using Apple LLVM version 9.1.0 (clang-1000.10.44.2).  All m-scripts were executed using MATLAB Version R2016a (build 9.0.0.341360).  All plots were created using MATLAB Version R2016a (build 9.0.0.341360).



\subsection{Example 1: Hyper-Sensitive Problem\label{subsect:hyperSensitive}}


Consider the following {\em hyper-sensitive} optimal control problem taken from Ref.~\cite{Rao4}.  The objective is to minimize the cost functional
\begin{equation}\label{eq:hyperSensitiveProblem-cost}
  \C{J} = \frac{1}{2} \int_{0}^{t_f} (x^2 + u^2) ~dt~,
\end{equation}
subject to the dynamic constraints
\begin{equation}\label{eq:hyperSensitiveProblem-dyn}
\dot{x} = -x^3+u~,
\end{equation}
and the boundary conditions
\begin{equation}\label{eq:hyperSensitiveProblem-BC}
\begin{array}{lcl}
x(0) = 1~, & & x(t_f) = 1.5~,
\end{array}
\end{equation}
where $t_f=10000$.  It is known for a sufficiently large value of $t_f$ the interesting behavior in the solution for the optimal control problem defined by Eqs.~\eqref{eq:hyperSensitiveProblem-cost} -- \eqref{eq:hyperSensitiveProblem-BC} occurs near $t=0$ and $t=t_f$ (see Ref.~\cite{Rao4} for details), while  the vast majority of the solution is a constant.  Given the structure of the solution, a majority of collocation points need to be placed near $t=0$ and $t=t_f$.  

\begin{figure}[htp]
\centering

\subfloat[$x(t)$ vs.~$t$.\label{fig:hyperSensitiveState}]{\epsfig{figure=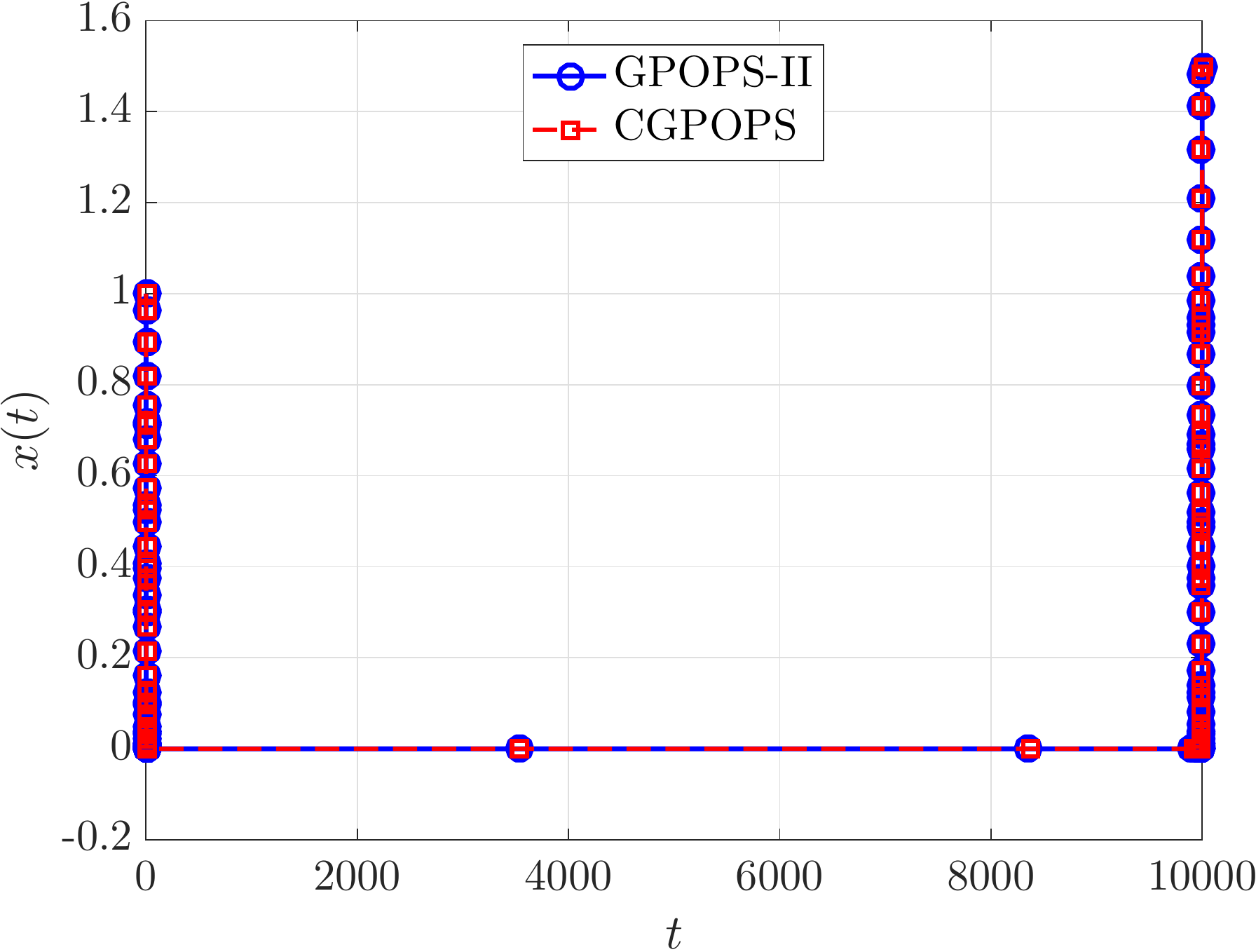,height=1.85in}}~~~\subfloat[$u(t)$ vs.~$t$. \label{fig:hyperSensitiveControl}]{\epsfig{figure=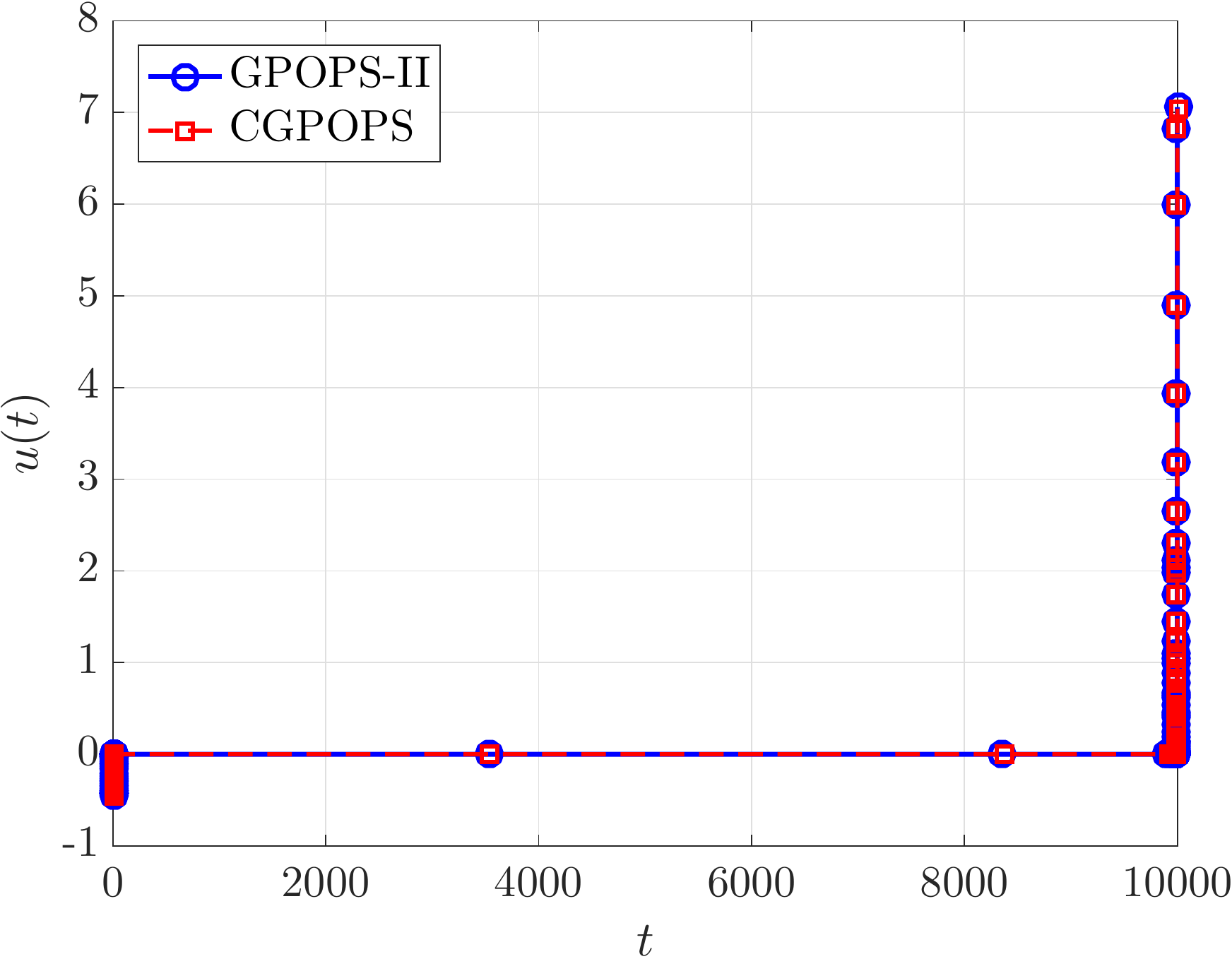,height=1.85in}}

\caption{$\CGPOPS$ and $\mathbb{GPOPS-II}$ Solutions to Example 1 using $hp$-IV(3,10). \label{fig:hyperSensitiveSolution}}
\end{figure}

\begin{figure}[h!]
\centering

\subfloat[$x(t)$ vs.~$t$ Near $t=0$.\label{fig:hyperSensitiveStateInitialLayer}]{\epsfig{figure=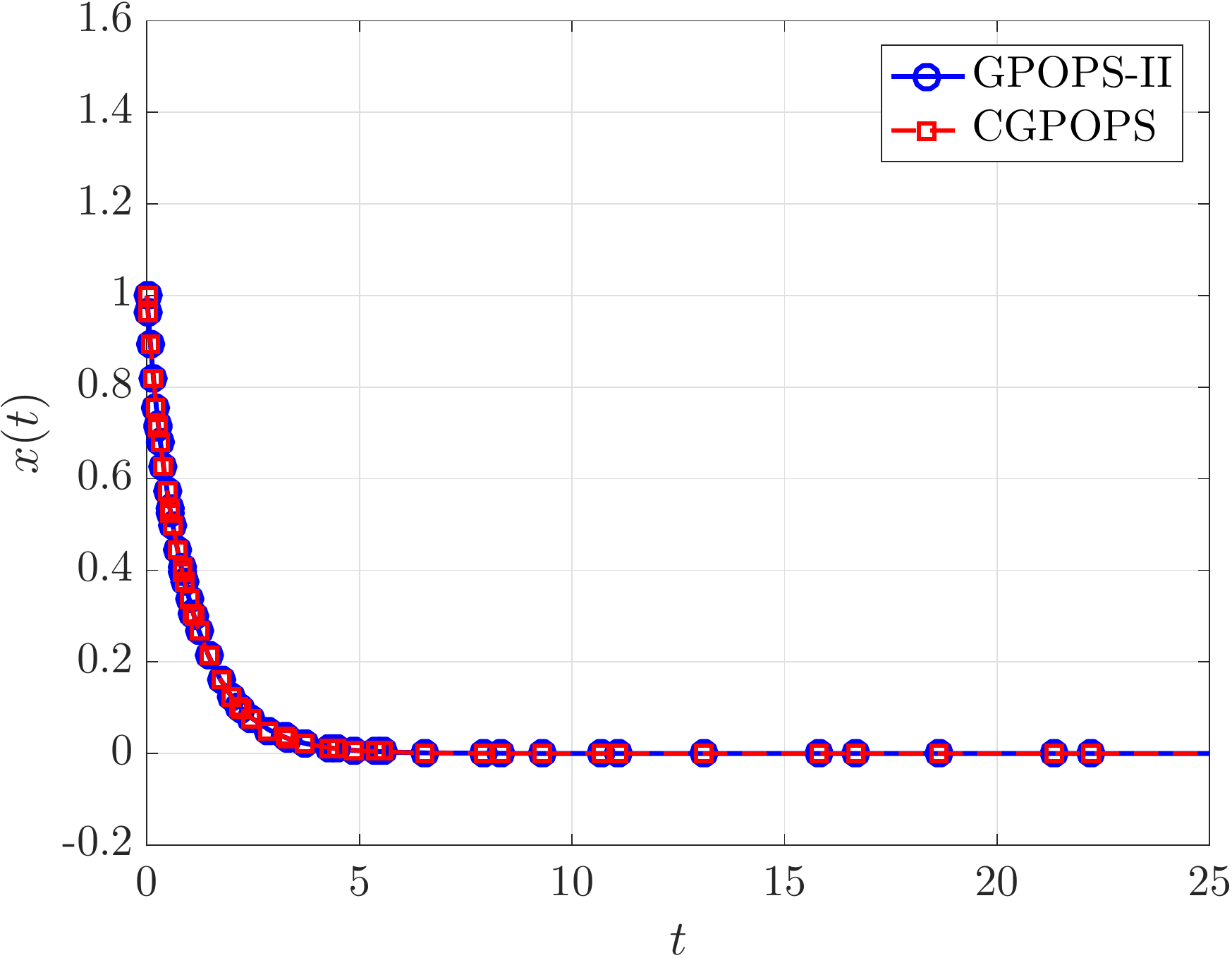,height=1.85in}}~~~\subfloat[$x(t)$ vs.~$t$ Near $t=t_f$. \label{fig:hyperSensitiveStateFinalLayer}]{\epsfig{figure=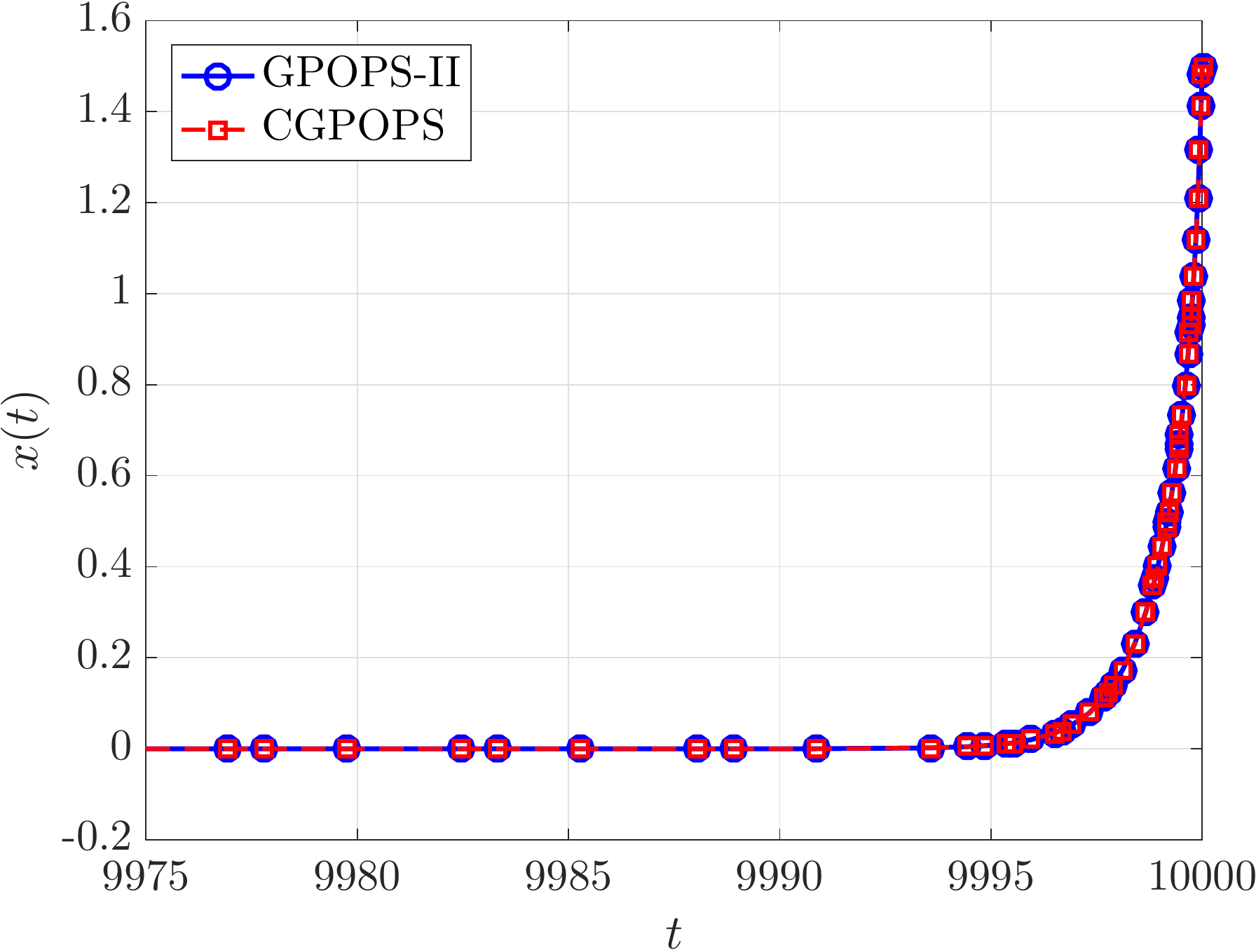,height=1.85in}}

\subfloat[$u(t)$ vs.~$t$ Near $t=0$.\label{fig:hyperSensitiveControlInitialLayer}]{\epsfig{figure=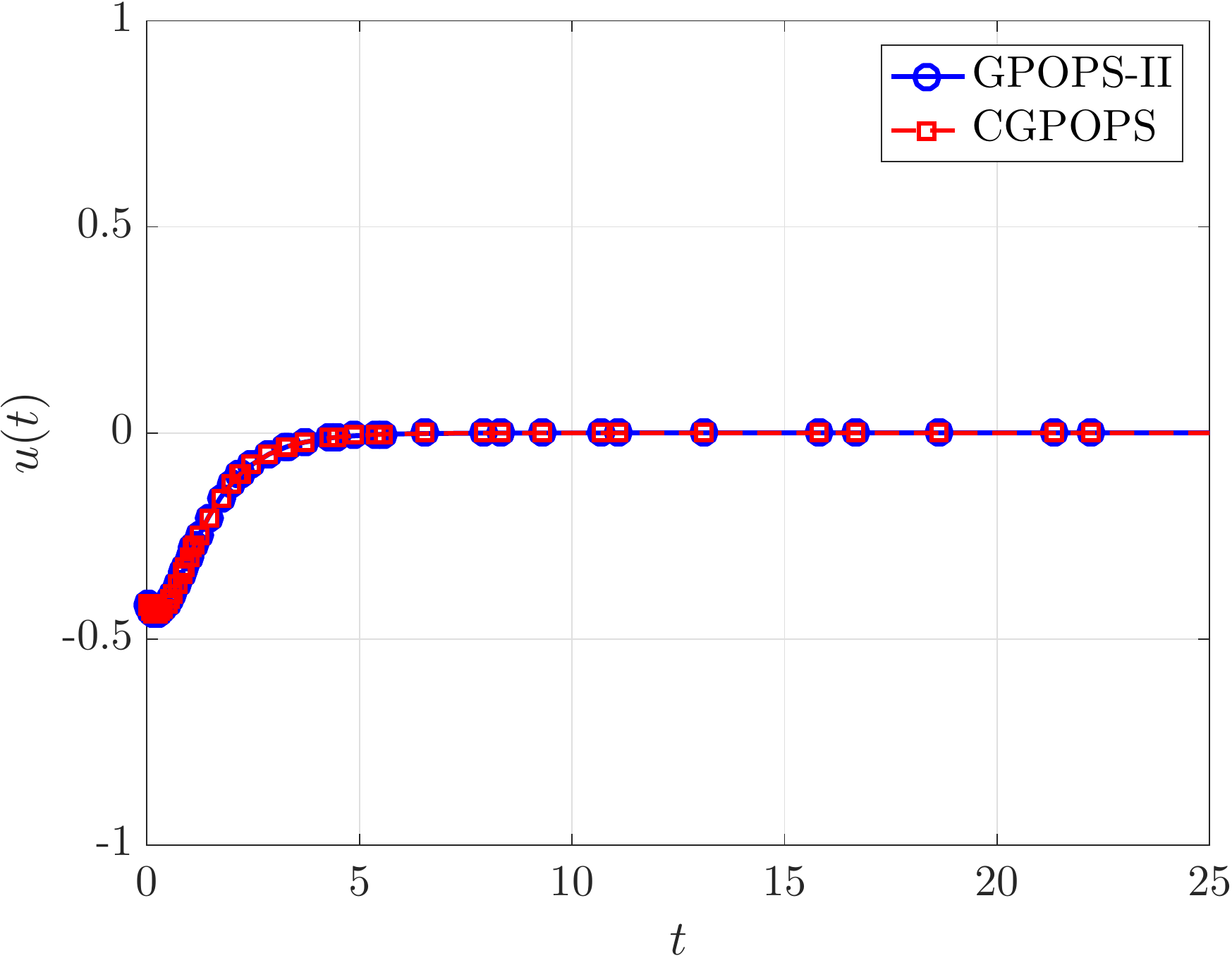,height=1.85in}}~~~\subfloat[$u(t)$ vs.~$t$ Near $t=t_f$. \label{fig:hyperSensitiveControlFinalLayer}]{\epsfig{figure=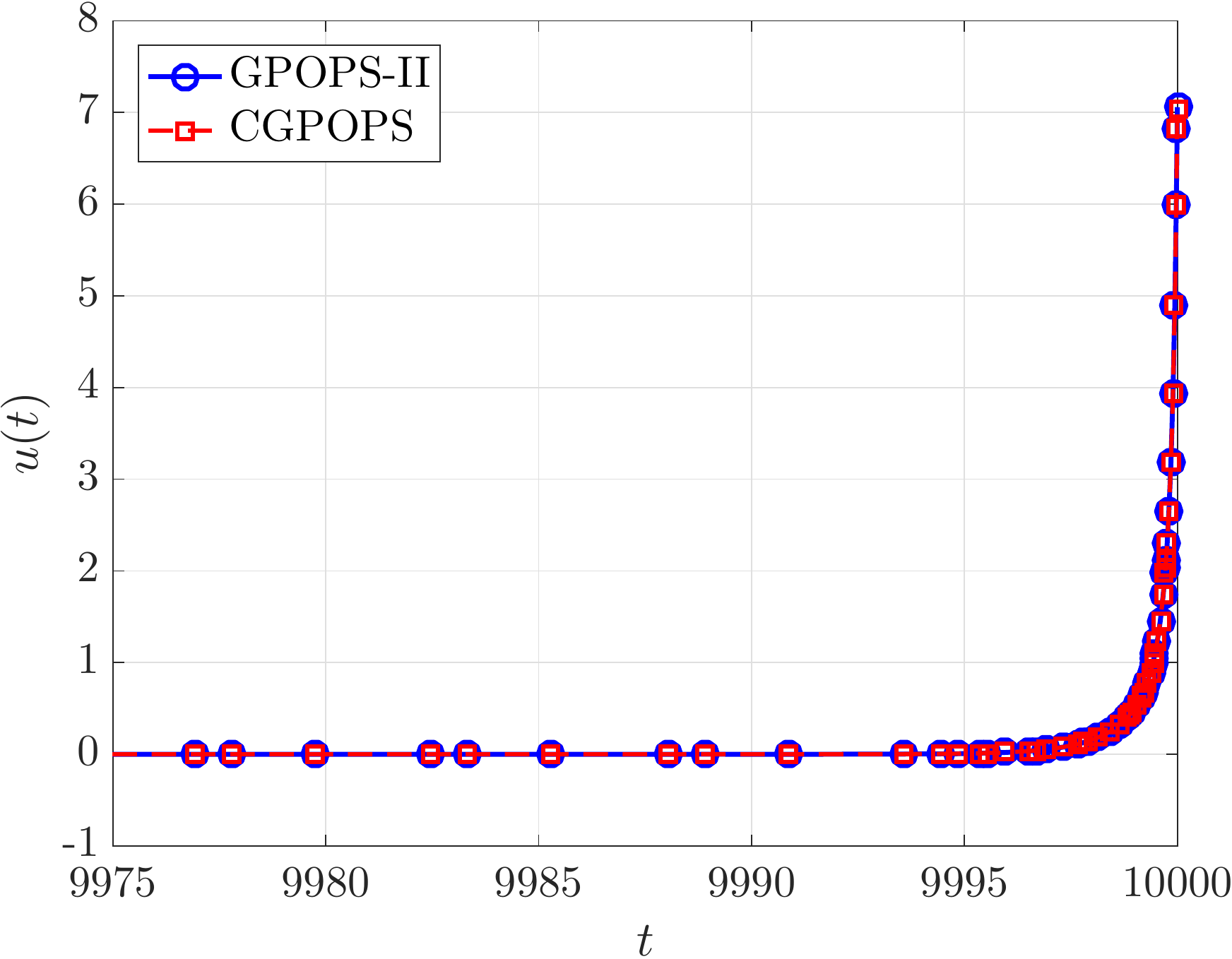,height=1.85in}}
\caption{$\CGPOPS$ and $\mathbb{GPOPS-II}$
  Solutions to Example 1 Near $t=0$ and $t=t_f$ using $hp$-IV(3,10). \label{fig:hyperSensitiveSolutionLayers}}
\end{figure}

\begin{figure}[htp]
\centering

\epsfig{file=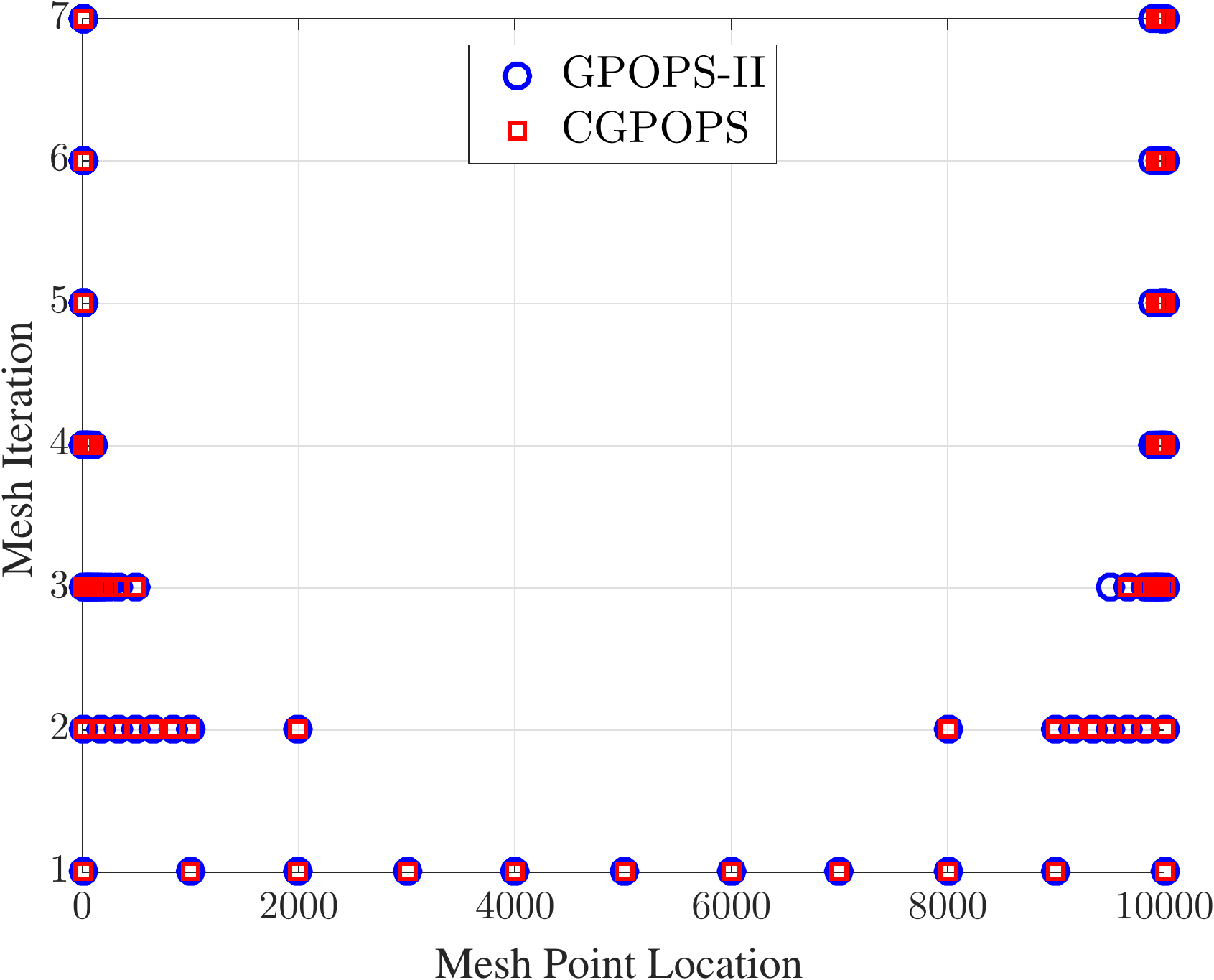,height=2in}
  
\caption{$\CGPOPS$ and $\mathbb{GPOPS-II}$ Mesh Refinement History for Example 1 Using $hp$-IV(3,10). \label{fig:hyperSensitiveMeshRefinement}}
\end{figure}

The hyper-sensitive optimal control problem was solved using $\CGPOPS$ with the mesh refinement methods $hp$-I(3,10), $hp$-II(3,10), $hp$-III(3,10), and $hp$-IV(3,10) on an initial mesh of ten evenly spaced mesh intervals with three LGR points per mesh interval.  The NLP solver and mesh refinement accuracy tolerances were set to $10^{-7}$ and $10^{-6}$, respectively.  The solution obtained using $\CGPOPS$ with the $hp$-IV(3,10) method is shown in Fig.~\ref{fig:hyperSensitiveSolution} alongside the solution obtained using $\mathbb{GPOPS-II}$ \cite{Patterson2014} with the $hp$-IV(3,10) method.  It is seen that the $\CGPOPS$ and $\mathbb{GPOPS-II}$ solutions are in excellent agreement.  Moreover, the optimal cost obtained using $\CGPOPS$ and $\mathbb{GPOPS-II}$ are extremely close, with values of $3.3620559$ and $3.3620559$, respectively, agreeing to the seventh decimal place.  Additionally, the computational time required by $\CGPOPS$ and $\mathbb{GPOPS-II}$ to solve the optimal control problem was $0.2153$ seconds and $1.5230$ seconds, respectively.  In order to demonstrate how $\CGPOPS$ is capable of capturing the interesting features of the optimal solution, Fig.~\ref{fig:hyperSensitiveSolutionLayers} shows the solution on the intervals $t\in [0,25]$ (near the initial time) and $t\in[9975,10000]$ (near the final time).  It is seen that $\CGPOPS$ accurately captures the rapid decay from $x(0)=1$ and the rapid growth to meet the terminal condition $x(t_f)=1.5$, with the density of the mesh points near $t=0$ and $t=t_f$ increasing as the mesh refinement progresses.  Additionally, Fig.~\ref{fig:hyperSensitiveMeshRefinement} shows the mesh refinement history.  Finally, Tables~\ref{tab:hyperSensitiveMeshRefinement-hpPatterson} -- \ref{tab:hyperSensitiveMeshRefinement-hpLegendre} show the estimated error on each mesh, where it is seen that the solution error decreases steadily with each mesh refinement iteration for all $hp$ methods employed.
 
\begin{table}[htp]
\centering
\caption{Performance of $\CGPOPS$ on Example 1 using $hp$-I(3,10). \label{tab:hyperSensitiveMeshRefinement-hpPatterson}}
\footnotesize
\begin{tabular}{|c|c|c|c|c|} \hline
{\bf Mesh} & {\bf Estimated}&{\bf Number of} &{\bf Estimated }&{\bf Number of} \\
{\bf Iteration} & {\bf Error ($\CGPOPS$)}&{\bf Collocation}&{\bf Error ($\mathbb{GPOPS-II}$)}&{\bf Collocation} \\
{\bf Number} & {\bf $hp$-I(3,10)}&{\bf Points}&{\bf $hp$-I(3,10)}&{\bf Points} \\\hline\hline
$1$	& $28.27$ 						& $31$ &	$28.27$ 						& $31$  \\ \hline
$2$	& $4.090$ 						& $67$ &	$4.090$ 						& $67$  \\ \hline
$3$	& $7.060\times 10^{-1}$	& $101$&	$7.060\times 10^{-1}$	& $101$ \\ \hline
$4$	& $1.661\times 10^{-1}$	&$134$ &	$1.661\times 10^{-1}$	&$134$ \\ \hline
$5$	& $1.476\times 10^{-2}$	& $158$ & $1.476\times 10^{-2}$	& $158$  \\ \hline
$6$	& $1.139\times 10^{-3}$	& $191$ & $1.139\times 10^{-3}$	& $191$  \\ \hline
$7$	& $7.557\times 10^{-7}$	& $218$ & $7.557\times 10^{-7}$	& $218$  \\ \hline
\end{tabular}
\normalsize
\end{table}

\begin{table}[htp]
\centering
\caption{Performance of $\CGPOPS$ on Example 1 using $hp$-II(3,10). \label{tab:hyperSensitiveMeshRefinement-hpDarby}}
\footnotesize
\begin{tabular}{|c|c|c|c|c|} \hline
{\bf Mesh} & {\bf Estimated}&{\bf Number of} &{\bf Estimated }&{\bf Number of} \\
{\bf Iteration} & {\bf Error ($\CGPOPS$)}&{\bf Collocation}&{\bf Error ($\mathbb{GPOPS-II}$)}&{\bf Collocation} \\
{\bf Number} & {\bf $hp$-II(3,10)}&{\bf Points}&{\bf $hp$-II(3,10)}&{\bf Points} \\\hline\hline
$1$	& $28.27$ 						& $31$	& $28.27$ 						& $31$  \\ \hline
$2$	& $1.667$ 						& $65$ & $1.667$ 						& $65$ \\ \hline
$3$	& $3.193$							& $106$& $3.193$						& $106$ \\ \hline
$4$	& $1.557\times 10^{-1}$	&$140$ &	$1.557\times 10^{-1}$	&$140$ \\ \hline
$5$	& $4.142\times 10^{-1}$	& $165$ & $4.142\times 10^{-1}$	& $165$ \\ \hline
$6$	& $1.261\times 10^{-2}$	& $185$ & $1.261\times 10^{-2}$	& $185$  \\ \hline
$7$	& $4.423\times 10^{-2}$	& $204$ & $4.423\times 10^{-2}$	& $204$  \\ \hline
$8$	& $4.707\times 10^{-4}$	&$209$ &	$4.707\times 10^{-4}$	&$209$ \\ \hline
$9$	& $1.090\times 10^{-3}$	& $226$ & $1.090\times 10^{-3}$	& $226$ \\ \hline
$10$	& $7.742\times 10^{-6}$	& $247$ & $7.742\times 10^{-6}$	& $247$  \\ \hline
$11$	& $7.470\times 10^{-7}$	& $250$ & $7.470\times 10^{-7}$	& $250$  \\ \hline
\end{tabular}
\normalsize
\end{table}

\begin{table}[htp]
\centering
\caption{Performance of $\CGPOPS$ on Example 1 using $hp$-III(3,10). \label{tab:hyperSensitiveMeshRefinement-hpLiu}}
\footnotesize
\begin{tabular}{|c|c|c|c|c|} \hline
{\bf Mesh} & {\bf Estimated}&{\bf Number of} &{\bf Estimated }&{\bf Number of} \\
{\bf Iteration} & {\bf Error ($\CGPOPS$)}&{\bf Collocation}&{\bf Error ($\mathbb{GPOPS-II}$)}&{\bf Collocation} \\
{\bf Number} & {\bf $hp$-III(3,10)}&{\bf Points}&{\bf $hp$-III(3,10)}&{\bf Points} \\\hline\hline
$1$	& $28.27$ 						& $31$ &	$28.27$ 						& $31$  \\ \hline
$2$	& $5.207$ 						& $22$ &	$5.207$ 						& $22$  \\ \hline
$3$	& $5.848\times 10^{-1}$	& $112$&	$5.848\times 10^{-1}$	& $112$ \\ \hline
$4$	& $9.156\times 10^{-2}$	&$139$ &	$9.156\times 10^{-2}$	&$142$ \\ \hline
$5$	& $5.732\times 10^{-3}$	& $115$ & $5.732\times 10^{-3}$	& $112$  \\ \hline
$6$	& $9.927\times 10^{-5}$	& $146$ & $9.927\times 10^{-5}$	& $146$ \\ \hline
$7$	& $2.451\times 10^{-5}$	& $153$ & $2.451\times 10^{-5}$	& $153$ \\ \hline
$8$	& $8.237\times 10^{-7}$	& $160$ & $8.237\times 10^{-7}$	& $160$ \\ \hline
\end{tabular}
\normalsize
\end{table}
   
\begin{table}[htp]
\centering
\caption{Performance of $\CGPOPS$ on Example 1 using $hp$-IV(3,10). \label{tab:hyperSensitiveMeshRefinement-hpLegendre}}
\footnotesize
\begin{tabular}{|c|c|c|c|c|} \hline
{\bf Mesh} & {\bf Estimated}&{\bf Number of} &{\bf Estimated }&{\bf Number of} \\
{\bf Iteration} & {\bf Error ($\CGPOPS$)}&{\bf Collocation}&{\bf Error ($\mathbb{GPOPS-II}$)}&{\bf Collocation} \\
{\bf Number} & {\bf $hp$-IV(3,10)}&{\bf Points}&{\bf $hp$-IV(3,10)}&{\bf Points} \\\hline\hline
$1$	& $28.27$ 						& $31$ &	$28.27$ 						& $31$  \\ \hline
$2$	& $4.763$ 						& $46$ &	$4.763$ 						& $46$  \\ \hline
$3$	& $8.214\times 10^{-1}$	& $52$&	$8.214\times 10^{-1}$	& $55$ \\ \hline
$4$	& $1.813\times 10^{-1}$	&$55$ &	$1.813\times 10^{-1}$	&$58$ \\ \hline
$5$	& $2.114\times 10^{-2}$	& $61$ & $2.114\times 10^{-2}$	& $61$  \\ \hline
$6$	& $1.688\times 10^{-3}$	& $87$ & $1.688\times 10^{-3}$	& $87$  \\ \hline
$7$	& $8.991\times 10^{-7}$	& $106$ & $8.991\times 10^{-7}$	& $106$  \\ \hline
\end{tabular}
\normalsize
\end{table}

\clearpage

\clearpage

\subsection{Example 2: Reusable Launch Vehicle Entry\label{subsect:rlvEntry}}


Consider the following optimal control problem of maximizing the crossrange during the atmospheric entry of a reusable  launch vehicle taken from Ref.~\cite{Betts3} where the numerical values in Ref.~\cite{Betts3} are converted from English units to SI units.  Maximize the cost functional
\begin{equation}\label{eq:rlv-cost}
  \cal{J} = \phi(t_f)~,
\end{equation}
subject to the dynamic constraints
\begin{equation}\label{eq:rlv-eom}
  \begin{array}{lclclcl}
    \dot{r} & = & v \sin \gamma~, & & 
    \dot{\theta} & = & \displaystyle \frac{v \cos \gamma \sin \psi}{r \cos
      \phi}~, \\
    \dot{\phi} & = & \displaystyle \frac{v \cos \gamma \cos \psi}{r}~, & &
    \dot{v} & = & \displaystyle -\frac{D}{m}-g\sin \gamma~, \\
    \dot{\gamma} & = & \displaystyle \frac{L \cos \sigma}{ m v} -\left(\frac{g}{v}-\frac{v}{r}\right)\cos\gamma~, & & 
    \dot{\psi}  & = & \displaystyle \frac{L \sin \sigma}{m v \cos \gamma}+\frac{v \cos \gamma \sin \psi \tan \phi }{r}~,  
  \end{array}
\end{equation}
and the boundary conditions
\begin{equation} \label{eq:rlv-bcs}
\begin{array}{lclclcl}
h(0) & = & 79248 \textrm{ km}~, & & h(t_f) & = &  24384 \textrm{ km}~, \\
\theta(0) & = & 0 \textrm{ deg}~, & & \theta(t_f) & = & \textrm{Free}~, \\
\phi(0) & = & 0 \textrm{ deg}~, & & \phi(t_f) & = & \textrm{Free}~, \\
v(0) & = & 7.803 \textrm{ km/s}~, & & v(t_f) & = & 0.762 \textrm{ km/s}~, \\
\gamma(0) & = & -1 \textrm{ deg}~, & & \gamma(t_f) & = & -5 \textrm{ deg}~, \\
\psi(0) & = & 90 \textrm{ deg}~, & & \psi(t_f) & = & \textrm{Free}~,
\end{array}
\end{equation}
where $r=h+R_e$ is the geocentric radius, $h$ is the altitude, $R_e$ is the polar radius of the Earth, $\theta$ is the longitude, $\phi$ is the latitude, $v$ is the speed, $\gamma$ is the flight path angle, $\psi$ is the azimuth angle and $m$ is the mass of the vehicle.  Furthermore, the aerodynamic and gravitational forces are computed as 
\begin{equation}\label{eq:rlv-aux}
  \begin{array}{lclclclclcl}
    D & = & \rho v^2 S C_D/2~, & & L & = & \rho v^2 S C_L/2~, & & g & = & \mu/r^2~,
  \end{array}
\end{equation}
where $\rho = \rho_0\exp(-h/H)$ is the atmospheric density, $\rho_0$ is the density at sea level, $H$ is the density scale height, $S$ is the vehicle reference area, $C_D$ is the coefficient of drag, $C_L$ is the coefficient of lift, and $\mu$ is the gravitational parameter.

\begin{figure}[ht!]
\centering

\subfloat[$h(t)$ vs.~$t$.\label{fig:rlvEntryAltitude}]{\epsfig{figure=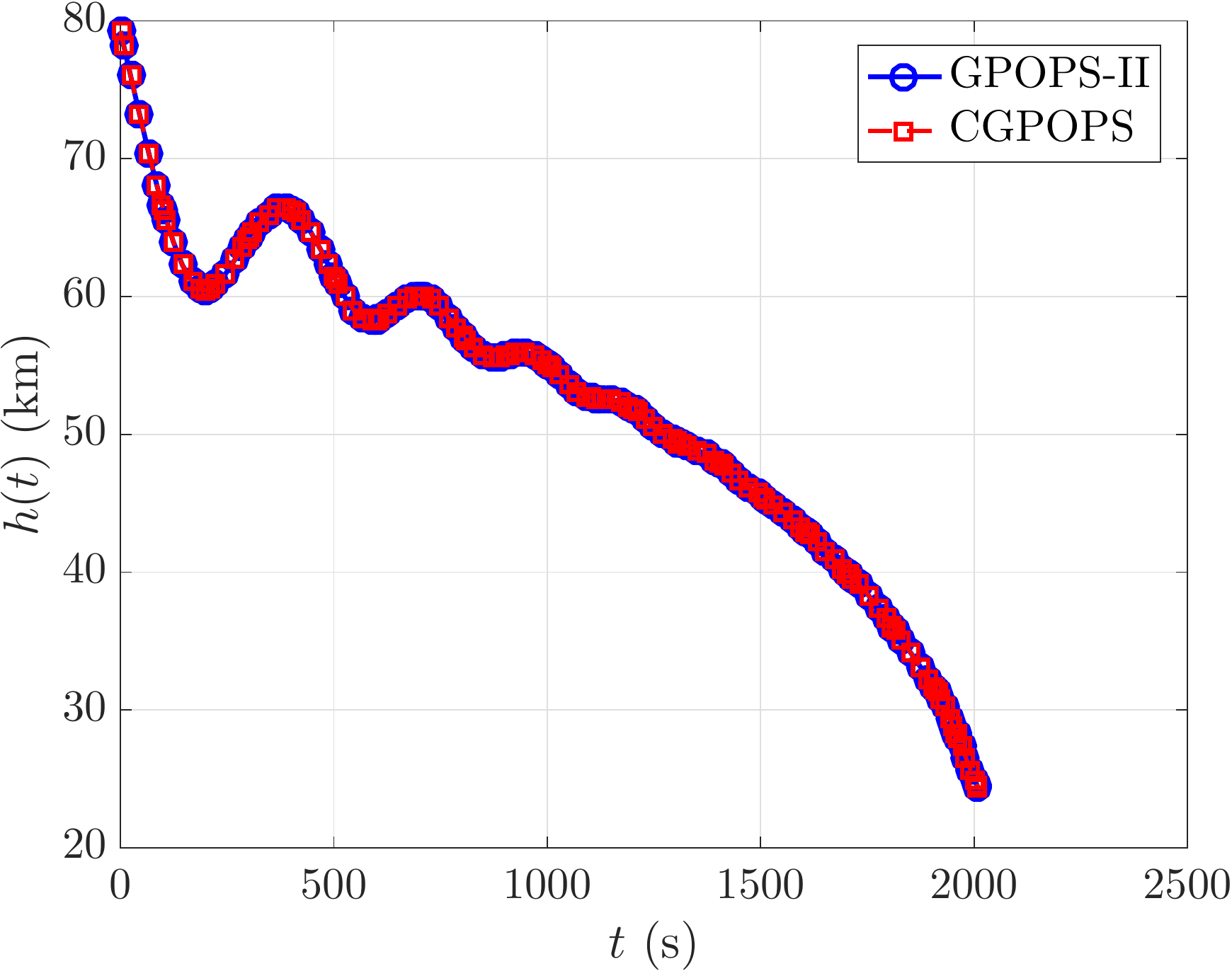,height=2in}}~~~\subfloat[$v(t)$ vs.~$t$. \label{fig:rlvEntrySpeed}]{\epsfig{figure=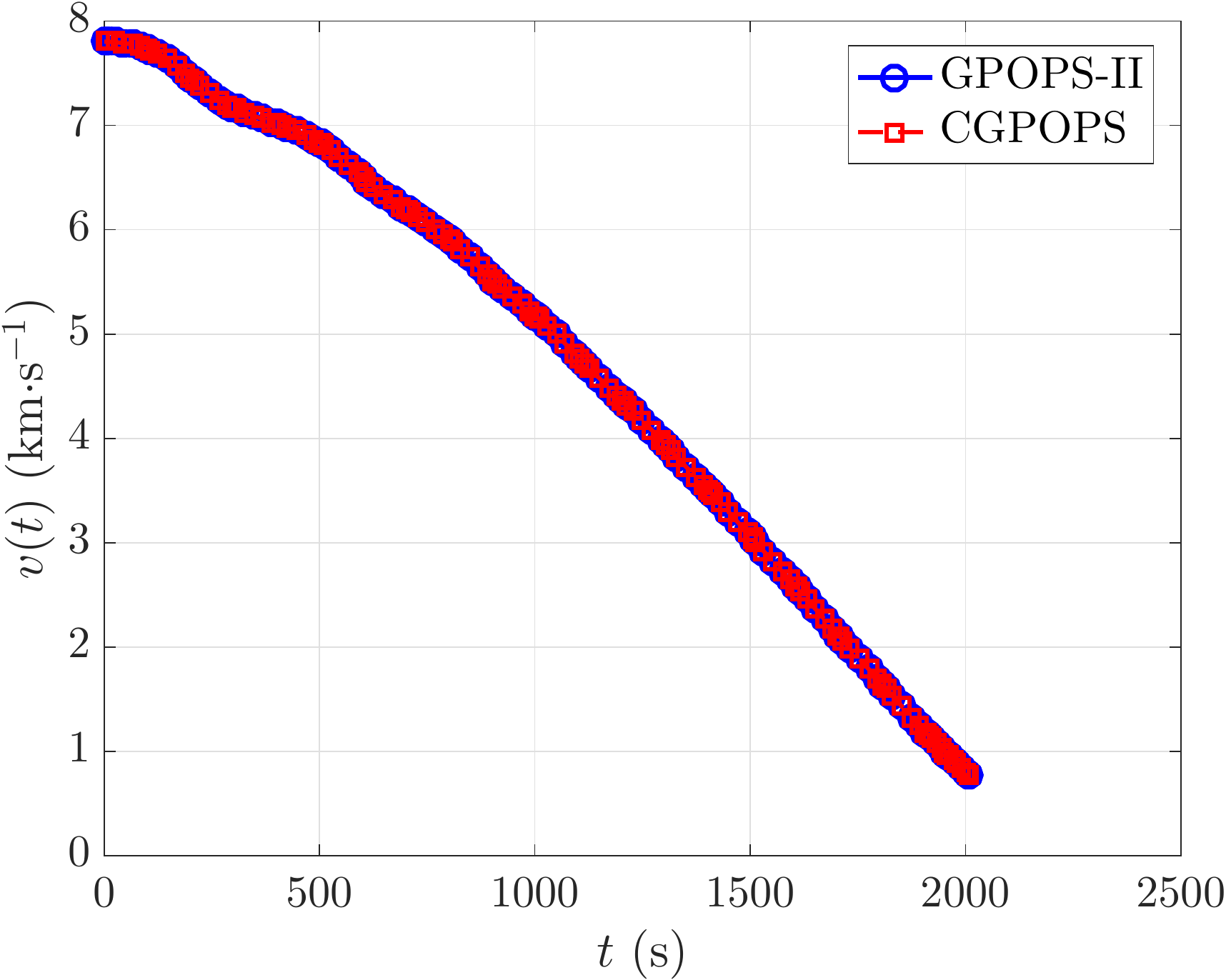,height=2in}}

\subfloat[$\phi(t)$ vs.~$t$.\label{fig:rlvEntryLat}]{\epsfig{figure=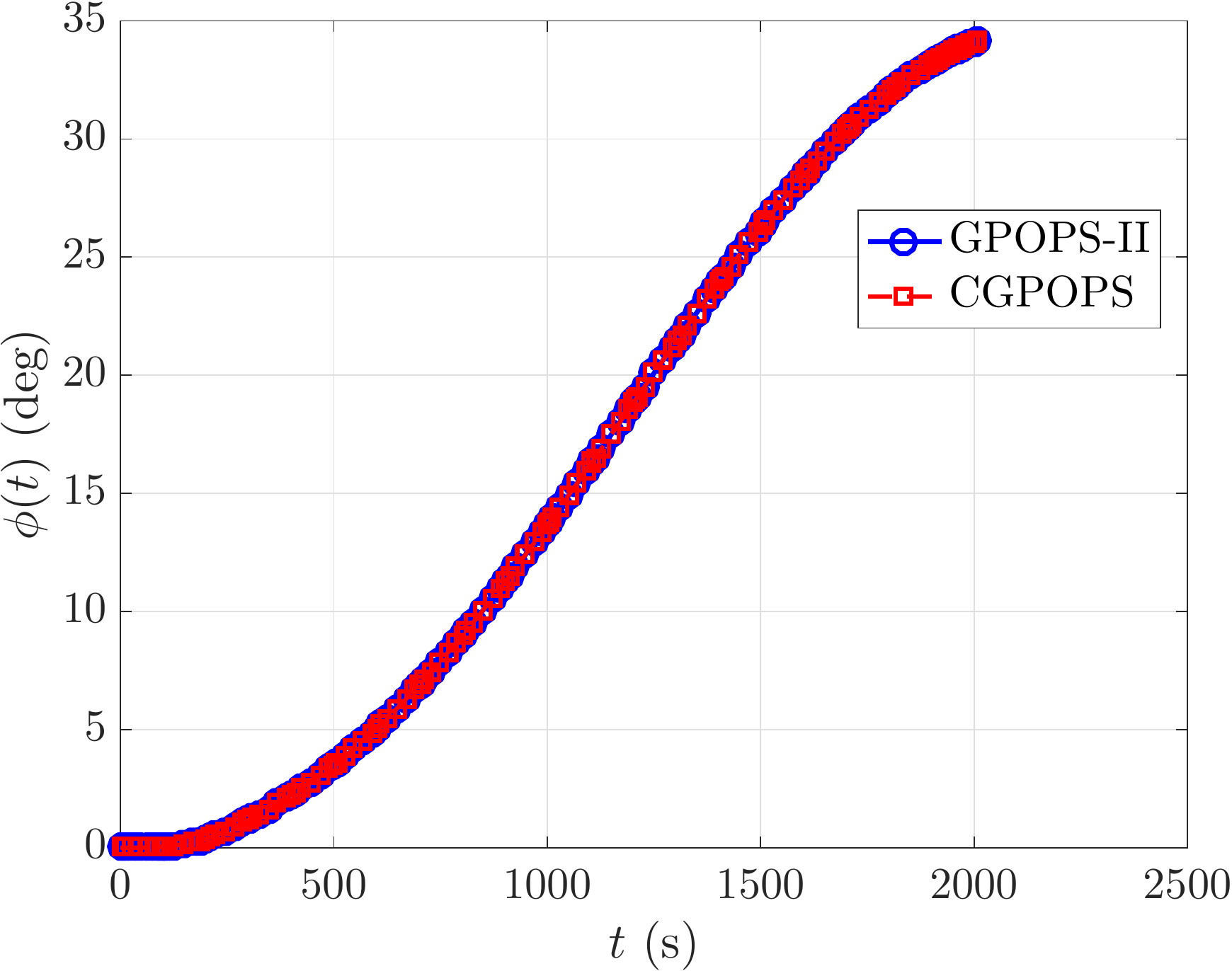,height=2in}}~~~\subfloat[$\gamma(t)$ vs.~$t$. \label{fig:rlvEntryFlightPathAngle}]{\epsfig{figure=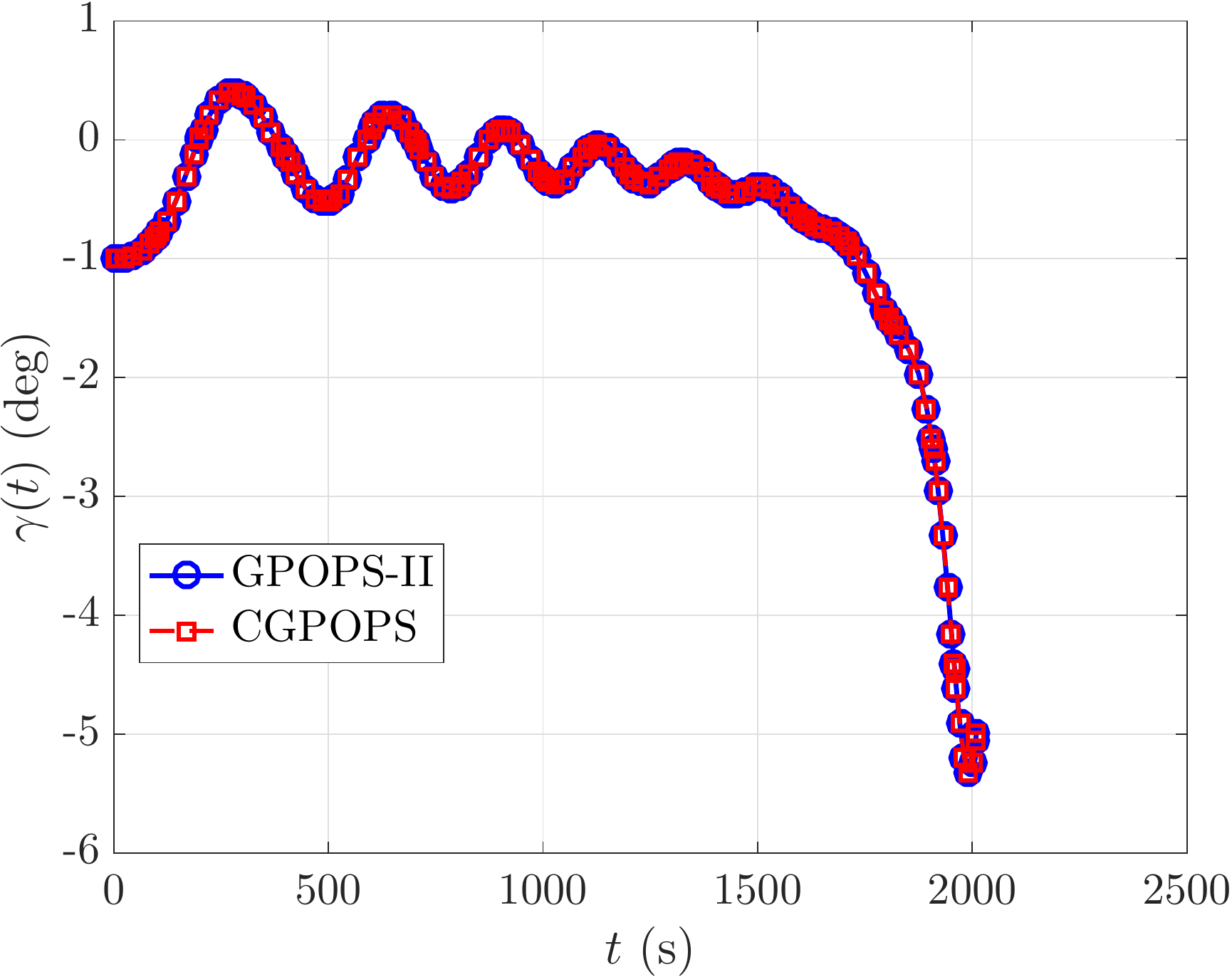,height=2in}}
  
\caption{$\CGPOPS$ and $\mathbb{GPOPS-II}$ State Solutions to Example 2 Using $hp$-III(4,10). \label{fig:rlvEntryStateSolution}}
\end{figure}

\begin{figure}[ht!]
\centering

\subfloat[$\alpha(t)$ vs.~$t$.\label{fig:rlvEntryAngleofAttack}]{\epsfig{figure=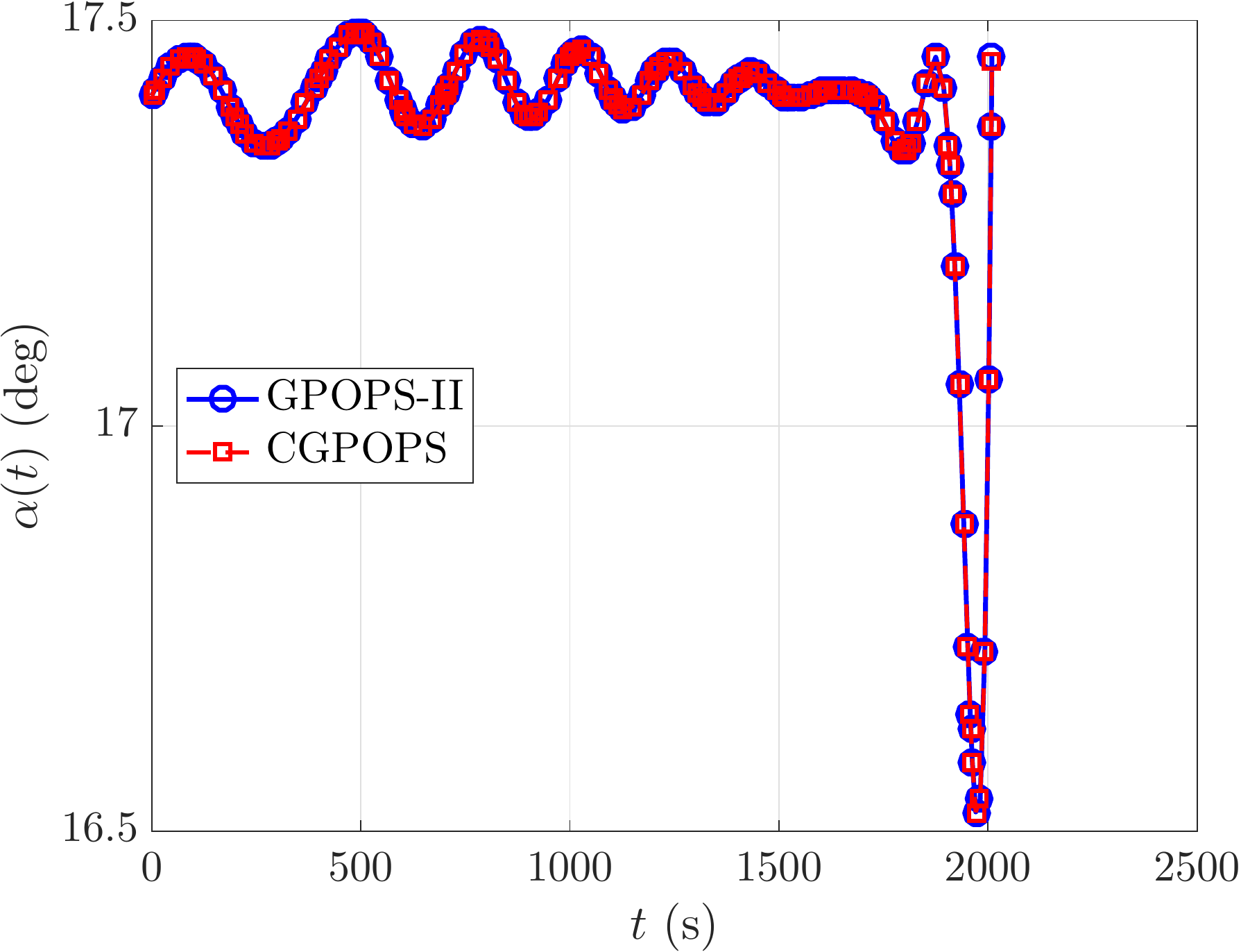,height=2in}}~~~\subfloat[$\sigma(t)$ vs.~$t$. \label{fig:rlvEntryBankAngle}]{\epsfig{figure=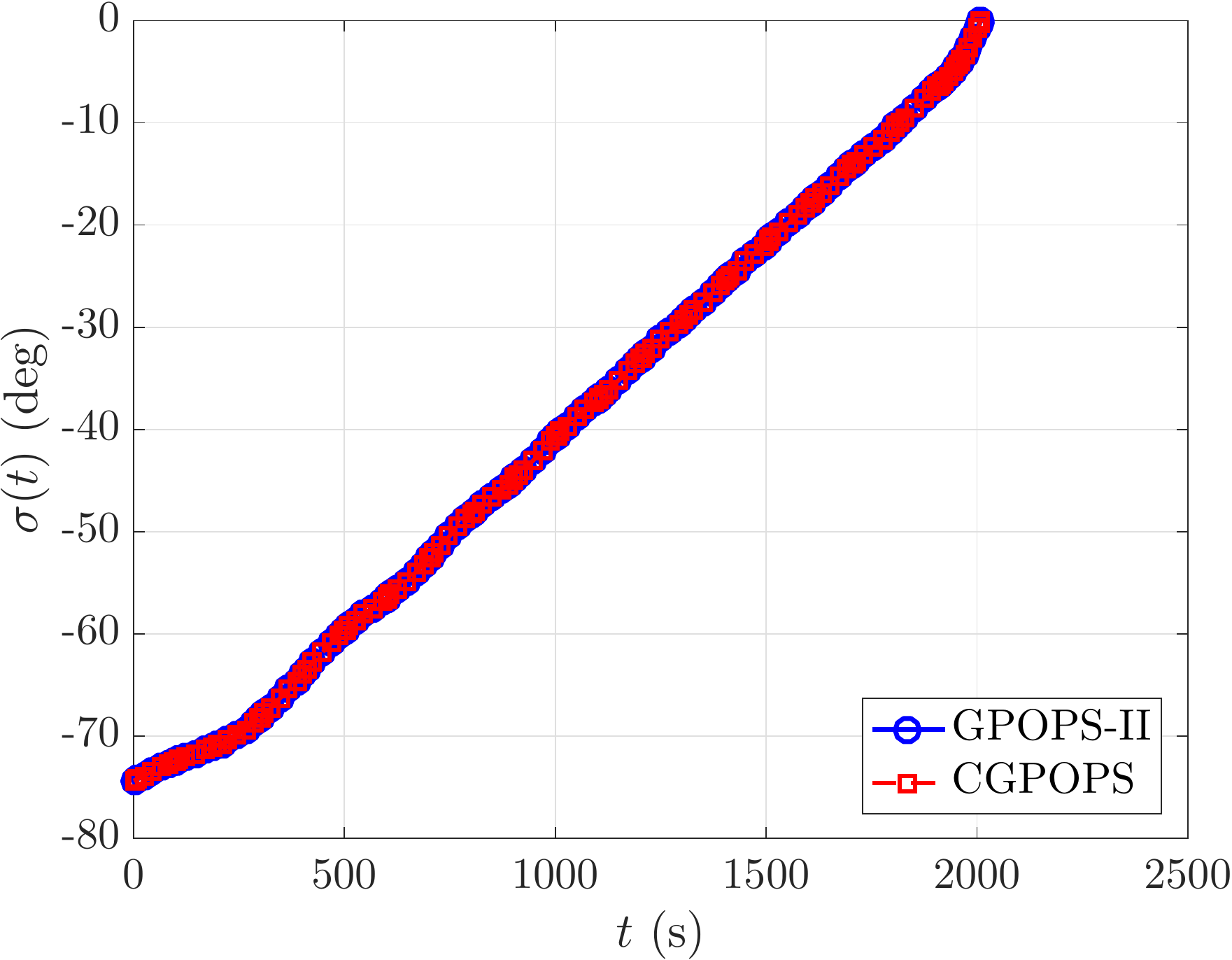,height=2in}}
  
\caption{$\CGPOPS$ and $\mathbb{GPOPS-II}$ Control Solutions to Example 2 Using $hp$-III(4,10). \label{fig:rlvEntryControlSolution}}
\end{figure}

The reusable launch vehicle entry optimal control problem defined by Eqs.~\eqref{eq:rlv-cost} -- \eqref{eq:rlv-aux} was solved using $\CGPOPS$ with the $hp$-I(4,10), $hp$-II(4,10), $hp$-III(4,10), and $hp$-IV(4,10) methods on an initial mesh consisting of ten evenly spaced mesh intervals with four LGR points per mesh interval.  The NLP solver and mesh refinement accuracy tolerances were both set to $10^{-7}$.  The initial guess of the state was a straight line over the duration $t\in[0,1000]$ between the known initial and final components of the state or a constant at the initial values of the components of the state whose terminal values are not specified, while the initial guess of both controls was zero.  Tables~\ref{tab:rlvEntryPerformance-hpPatterson} -- \ref{tab:rlvEntryPerformance-hpLegendre} show the performance of both $\CGPOPS$ and $\mathbb{GPOPS-II}$ on this example for the four $hp$ methods, where the mesh refinement history is nearly identical for each of the respective methods.  The solution obtained using $\CGPOPS$ with the $hp$-III(4,10) method is shown in Figs.~\ref{fig:rlvEntryAltitude} -- \ref{fig:rlvEntryBankAngle} alongside the solution obtained using the software $\mathbb{GPOPS-II}$ \cite{Patterson2014} with the $hp$-III(4,10) method, where it is seen that the two solutions obtained are essentially identical.  It is noted that the optimal cost obtained by $\CGPOPS$ and $\mathbb{GPOPS-II}$ are also identical at $0.59627639$ and $0.59627639$, respectively.  Finally, the computational time used by $\CGPOPS$ is approximately half the amount of time required by $\mathbb{GPOPS-II}$ to solve the optimal control problem, taking 0.9105~seconds and 1.9323 seconds, respectively.

\begin{table}[htp]
\centering
\caption{Performance of $\CGPOPS$ on Example 2 using $hp$-I(4,10). \label{tab:rlvEntryPerformance-hpPatterson}}
\footnotesize
\begin{tabular}{|c|c|c|c|c|} \hline
{\bf Mesh} & {\bf Estimated}&{\bf Number of} &{\bf Estimated }&{\bf Number of} \\
{\bf Iteration} & {\bf Error ($\CGPOPS$)}&{\bf Collocation}&{\bf Error ($\mathbb{GPOPS-II}$)}&{\bf Collocation} \\
{\bf Number} & {\bf $hp$-I(4,10)}&{\bf Points}&{\bf $hp$-I(4,10)}&{\bf Points} \\\hline\hline
$1$& $2.463\times 10^{-3}$ & $41$ &$2.463\times 10^{-3}$ & $41$  \\ \hline
$2$ & $9.891\times 10^{-5}$ & $103$ &$9.896\times 10^{-5}$ & $103$ \\ \hline
$3$ & $3.559\times 10^{-6}$ & $118$&  $3.559\times 10^{-6}$ & $118$ \\ \hline
$4$ & $3.287\times 10^{-7}$ &$133$ & $3.287\times 10^{-7}$ &$133$ \\ \hline
$5$ & $8.706\times 10^{-8}$ & $134$ & $8.706\times 10^{-8}$ & $134$   \\ \hline
\end{tabular}
\normalsize
\end{table}

\begin{table}[htp]
\centering
\caption{Performance of $\CGPOPS$ on Example 2 using $hp$-II(4,10). \label{tab:rlvEntryPerformance-hpDarby}}
\footnotesize
\begin{tabular}{|c|c|c|c|c|} \hline
{\bf Mesh} & {\bf Estimated}&{\bf Number of} &{\bf Estimated }&{\bf Number of} \\
{\bf Iteration} & {\bf Error ($\CGPOPS$)}&{\bf Collocation}&{\bf Error ($\mathbb{GPOPS-II}$)}&{\bf Collocation} \\
{\bf Number} & {\bf $hp$-II(4,10)}&{\bf Points}&{\bf $hp$-II(4,10)}&{\bf Points} \\\hline\hline
$1$& $2.463\times 10^{-3}$ & $41$ &$2.463\times 10^{-3}$ & $41$  \\ \hline
$2$ & $6.026\times 10^{-6}$ & $193$ &$6.023\times 10^{-6}$ & $193$ \\ \hline
$3$ & $8.227\times 10^{-8}$ & $261$&  $8.227\times 10^{-8}$ & $261$ \\ \hline
\end{tabular}
\normalsize
\end{table}

\begin{table}[htp]
\centering
\caption{Performance of $\CGPOPS$ on Example 2 using $hp$-III(4,10). \label{tab:rlvEntryPerformance-hpLiu}}
\footnotesize
\begin{tabular}{|c|c|c|c|c|} \hline
{\bf Mesh} & {\bf Estimated}&{\bf Number of} &{\bf Estimated }&{\bf Number of} \\
{\bf Iteration} & {\bf Error ($\CGPOPS$)}&{\bf Collocation}&{\bf Error ($\mathbb{GPOPS-II}$)}&{\bf Collocation} \\
{\bf Number} & {\bf $hp$-III(4,10)}&{\bf Points}&{\bf $hp$-III(4,10)}&{\bf Points} \\\hline\hline
$1$& $2.463\times 10^{-3}$ & $41$ &$2.463\times 10^{-3}$ & $41$  \\ \hline
$2$ & $2.850\times 10^{-5}$ & $71$ &$2.850\times 10^{-5}$ & $71$ \\ \hline
$3$ & $2.065\times 10^{-6}$ & $141$&  $2.065\times 10^{-6}$ & $141$ \\ \hline
$4$ & $8.887\times 10^{-8}$ &$148$ & $8.887\times 10^{-8}$ &$148$ \\ \hline
\end{tabular}
\normalsize
\end{table}

\begin{table}[htp]
\centering
\caption{Performance of $\CGPOPS$ on Example 2 using $hp$-IV(4,10). \label{tab:rlvEntryPerformance-hpLegendre}}
\footnotesize
\begin{tabular}{|c|c|c|c|c|} \hline
{\bf Mesh} & {\bf Estimated}&{\bf Number of} &{\bf Estimated }&{\bf Number of} \\
{\bf Iteration} & {\bf Error ($\CGPOPS$)}&{\bf Collocation}&{\bf Error ($\mathbb{GPOPS-II}$)}&{\bf Collocation} \\
{\bf Number} & {\bf $hp$-IV(4,10)}&{\bf Points}&{\bf $hp$-IV(4,10)}&{\bf Points} \\\hline\hline
$1$& $2.463\times 10^{-3}$ & $41$ &$2.463\times 10^{-3}$ & $41$  \\ \hline
$2$ & $2.364\times 10^{-5}$ & $122$ &$3.364\times 10^{-5}$ & $122$ \\ \hline
$3$ & $3.286\times 10^{-7}$ & $200$&  $3.286\times 10^{-7}$ & $192$ \\ \hline
$4$ & $9.561\times 10^{-8}$ &$203$ & $1.285\times 10^{-7}$ &$194$ \\ \hline
$5$ & -- & -- & $9.561\times 10^{-8}$ & $195$   \\ \hline
\end{tabular}
\normalsize
\end{table}

\clearpage

\begin{figure}[ht!]
\centering

\epsfig{file=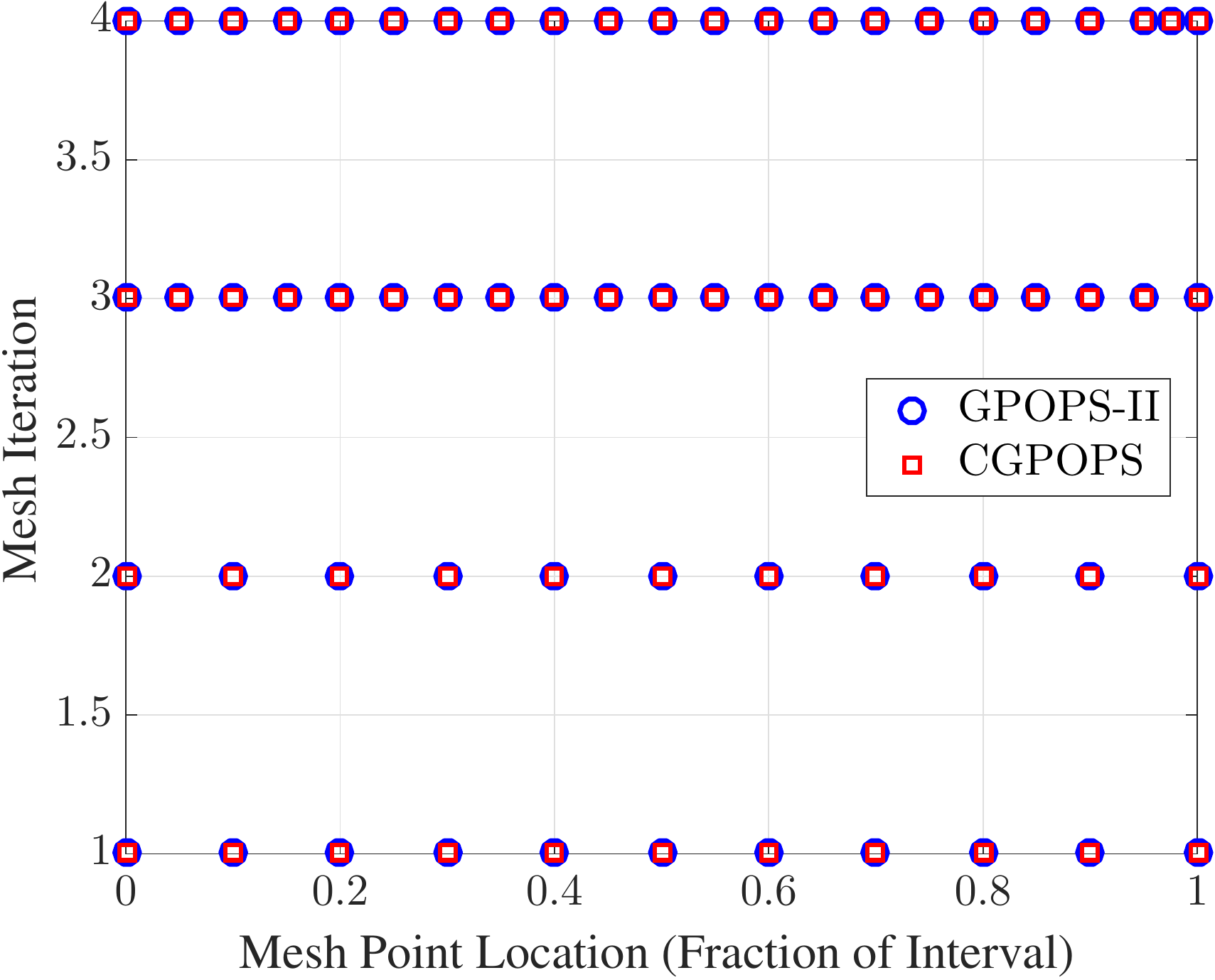,height=2in}
  
\caption{$\CGPOPS$ and $\mathbb{GPOPS-II}$ Mesh Refinement History for Example 2 Using $hp$-III(4,10). \label{fig:rlvEntryMeshRefinement}}
\end{figure}


\subsection{Example 3: Space Station Attitude Control\label{subsect:spaceStation}}


Consider the following space station attitude control optimal control problem taken from Refs.~\cite{Pietz2003} and \cite{Betts3}.  Minimize the cost functional
\begin{equation}\label{eq:space-station-cost}
  \cal{J} = \displaystyle \frac{1}{2}\int_{t_0}^{t_f} \m{u}\tr\m{u}~dt~,
\end{equation}
subject to the dynamic constraints
\begin{equation}\label{eq:space-station-dynamics}
  \begin{array}{lcl}
    \g{\omega} & = & \displaystyle \m{J}^{-1}\left\{\g{\tau}_{gg}(\m{r})-\g{\omega}^{\otimes}\left[\m{J}\g{\omega}+\m{h}\right]-\m{u}\right\}~, \\
    \dot{\m{r}} & = &   \displaystyle \frac{1}{2}\left[\m{r}\m{r}\tr+\m{I}+\m{r}\right]\left[\g{\omega}-\g{\omega}(\m{r})\right]~, \\
    \dot{\m{h}} & = & \m{u}~,
  \end{array}
\end{equation}
the inequality path constraint
\begin{equation}
  \left\|\m{h}\right\| \leq h_{\max}~,
\end{equation}
and the boundary conditions
\begin{equation}\label{eq:space-station-boundary-conditions}
\begin{array}{l}
  \begin{array}{rclclcl} 
    \hspace{0.16in} t_0 & = & 0~, & &
    \hspace{0.2in} t_f & = & 1800~, \vspace{3pt}\\
  \end{array} \vspace{3pt}\\
  \begin{array}{rclclclclcl} 
    \g{\omega}(0) & = & \bar{\g{\omega}}_0~, & &
    \m{r}(0) & = & \bar{\m{r}}_0~, & &
    \m{h}(0) & = & \bar{\m{h}}_0~, \vspace{3pt}\\
  \end{array} \\ 
  \begin{array}{rcl} 
    \hspace{0.19in} \mathbf{0} & = & \displaystyle \mathbf{J}^{-1}\left\{\boldsymbol{\tau}_{gg}(\mathbf{r}(t_f)) - \boldsymbol{\omega}^{\otimes}(t_f)\left[\mathbf{J}\boldsymbol{\omega}(t_f)+\mathbf{h}(t_f)\right]\right\}~, \vspace{3pt}\\
    \mathbf{0} & = & \displaystyle \frac{1}{2}\left[\mathbf{r}(t_f)\mathbf{r}^{\textsf{T}}(t_f)+\mathbf{I}+\mathbf{r}(t_f)\right]\left[\boldsymbol{\omega}(t_f)-\boldsymbol{\omega}_0(\mathbf{r}(t_f))\right]~,
\end{array}
\end{array}
\end{equation}
where $(\g{\omega},\m{r},\m{h})$ is the state and $\m{u}$ is the control.  In this formulation $\g{\omega}$ is the angular velocity, $\m{r}$ is the Euler-Rodrigues parameter vector, $\m{h}$ is the angular momentum, and $\m{u}$ is the input moment (and is the control).   Furthermore,   
\begin{equation}\label{eq:space-station-om_orb-defs}
  \begin{array}{lclclcl}
    \g{\omega}_0(\m{r}) & = & -\omega_{\textrm{orb}}\m{C}_2~, & & \g{\tau}_{gg} & = & 3\omega_{\textrm{orb}}^2 \m{C}_3^{\otimes}\m{J}\m{C}_3~, \vspace{3pt}\\
    \omega_{\textrm{orb}} & = & \displaystyle 0.6511\frac{\pi}{180}~, & &  h_{\max} & = & 10000~, \vspace{3pt}\\
  \end{array}
\end{equation}
and $\m{C}_2$ and $\m{C}_3$ are the second and third column, respectively, of the matrix
\begin{equation}
 \m{C} = \m{I}+\frac{2}{1+\m{r}\tr\m{r}}\left(\m{r}^{\otimes}\m{r}^{\otimes}-\m{r}^{\otimes}\right)~.
\end{equation}
In this example the matrix $\m{J}$ is given as
\begin{equation}\label{eq:space-station-inertia-def}
  \m{J} = \left[
    \begin{array}{ccc} 2.80701911616\times 10^{7} & 4.822509936\times 10^{5} & -1.71675094448\times 10^{7}\\
      4.822509936\times 10^{5} & 9.5144639344\times 10^{7}  & 6.02604448\times 10^{4} \\
      -1.71675094448\times 10^{7} & 6.02604448\times 10^{4} & 7.6594401336\times 10^{7}
    \end{array}
  \right]~,
\end{equation}
while the initial conditions $\bar{\g{\omega}}_0$, $\bar{\m{r}}_0$, and
$\bar{\m{h}}_0$ are
\begin{equation}\label{eq:space-station-IC}
  \begin{array}{l}
  \begin{array}{lcl}
    \bar{\g{\omega}}_0 & = & \left[\begin{array}{c} -9.5380685844896\times 10^{-6}\\ -1.1363312657036\times 10^{-3} \\+5.3472801108427\times 10^{-6}\end{array}\right]~,
  \end{array} \\
  \\
  \begin{array}{lclclcl}
    ~\bar{\m{r}}_0 & = & \left[\begin{array}{c} 2.9963689649816\times10^{-3}  \\1.5334477761054\times 10^{-1} \\ 3.8359805613992\times 10^{-3}\end{array}\right]~, & &
    \bar{\m{h}}_0 & = & \left[\begin{array}{c} 5000 \\ 5000 \\ 5000 \end{array}\right]~.
  \end{array}
  \end{array}
\end{equation}
A more detailed description of this problem, including all of the constants $\m{J}$, $\bar{\g{\omega}}_0$, $\bar{\m{r}}_0$, and $\bar{\m{h}}_0$, can be found in \cite{Pietz2003} or \cite{Betts3}.

\begin{figure}[ht!]
\centering

\subfloat[$\omega_1(t)$ vs.~$t$.\label{fig:spaceStationState1}]{\epsfig{figure=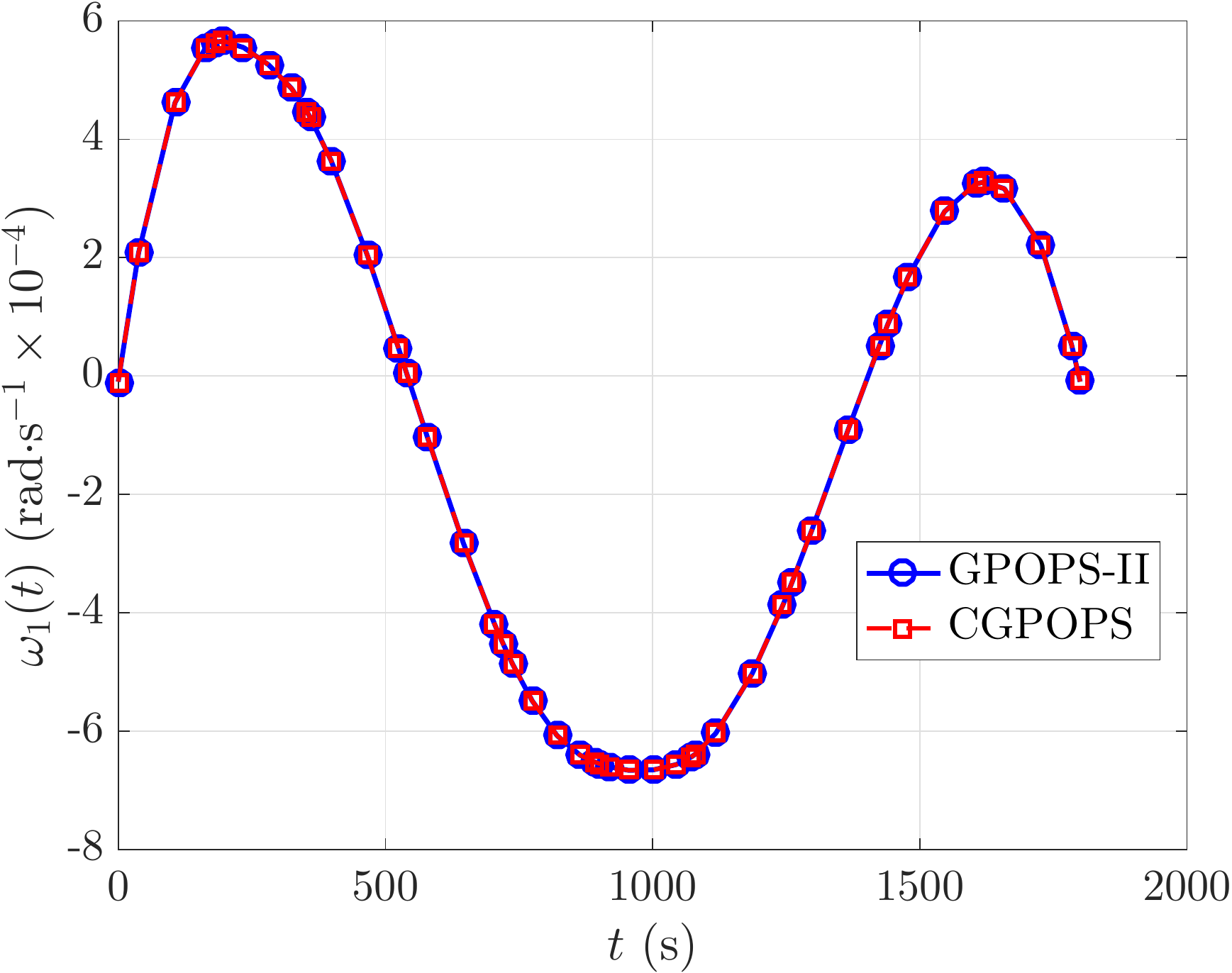,height=2in}}~\subfloat[$\omega_2(t)$ vs.~$t$. \label{fig:spaceStationState2}]{\epsfig{figure=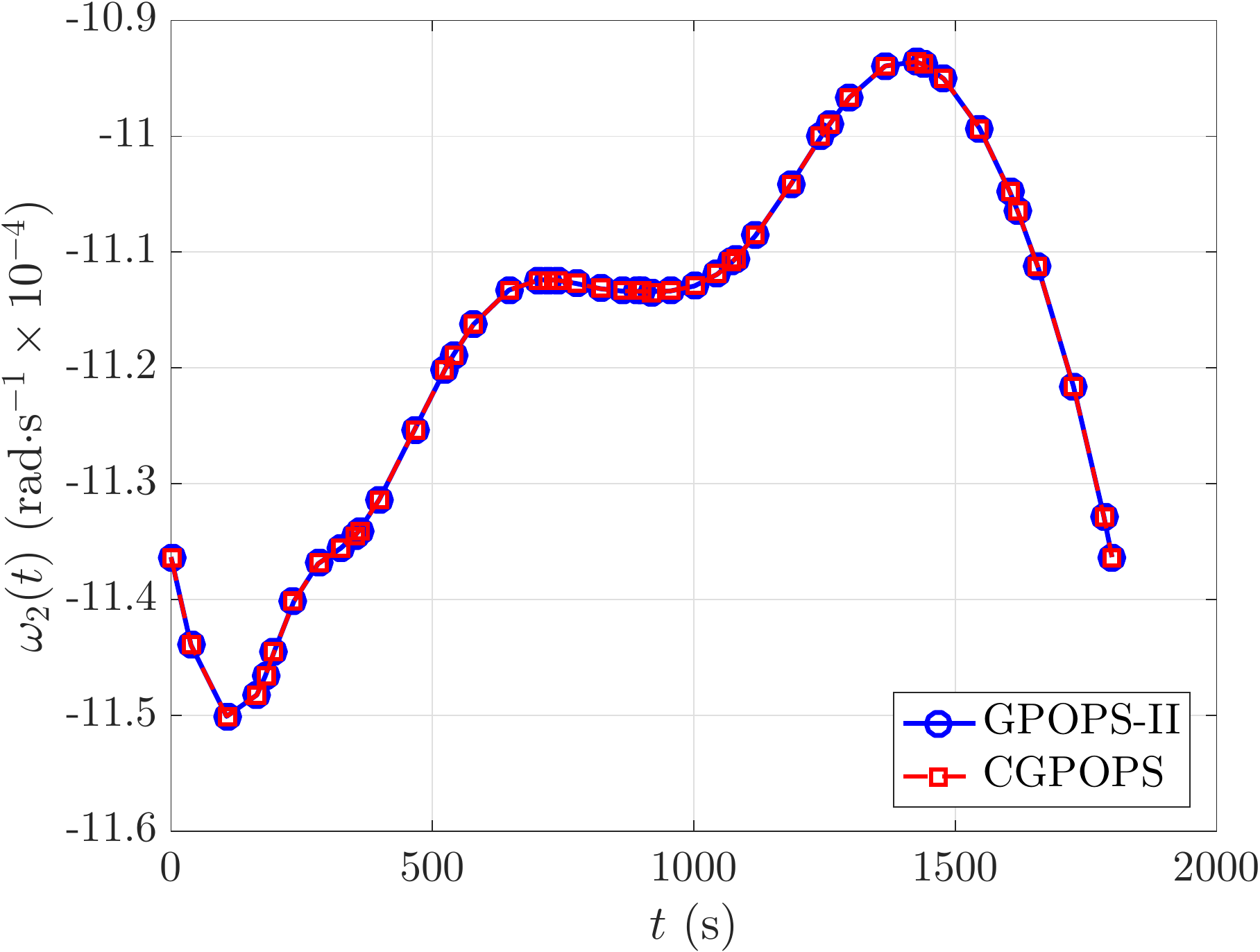,height=2in}}

\subfloat[$\omega_3(t)$ vs.~$t$.\label{fig:spaceStationState3}]{\epsfig{figure=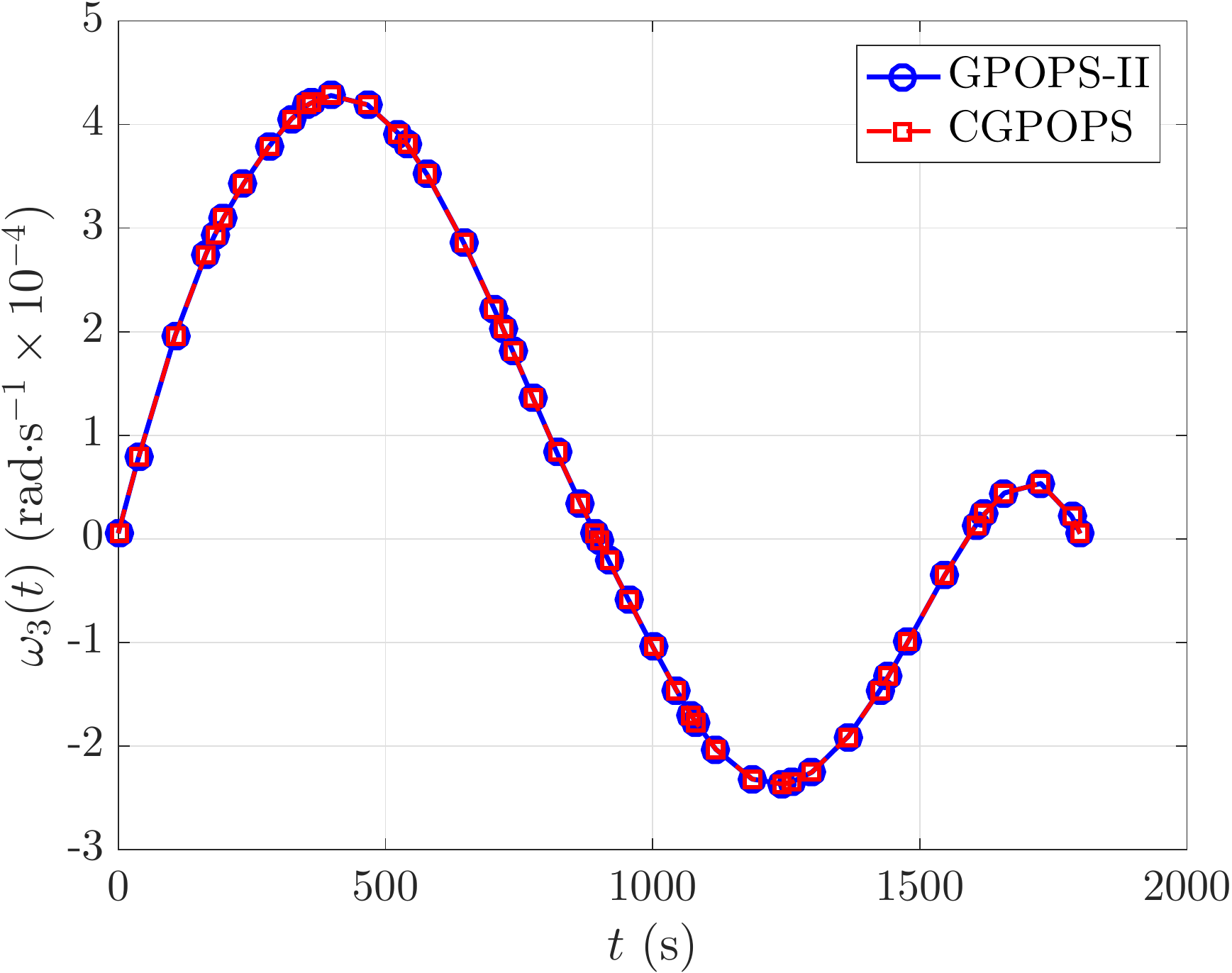,height=2in}}~\subfloat[$r_1(t)$ vs.~$t$. \label{fig:spaceStationState4}]{\epsfig{figure=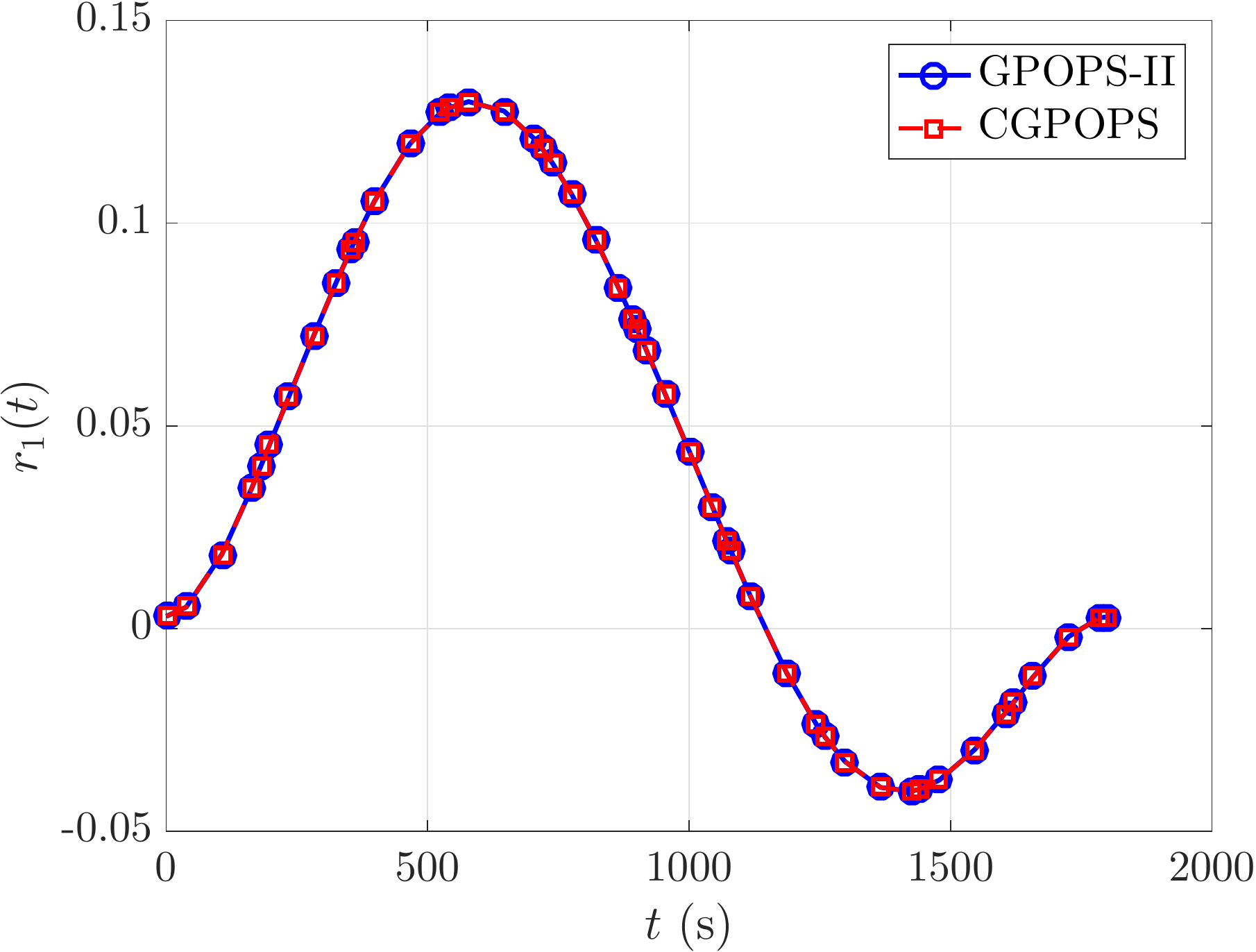,height=2in}}

\subfloat[$r_2(t)$ vs.~$t$.\label{fig:spaceStationState5}]{\epsfig{figure=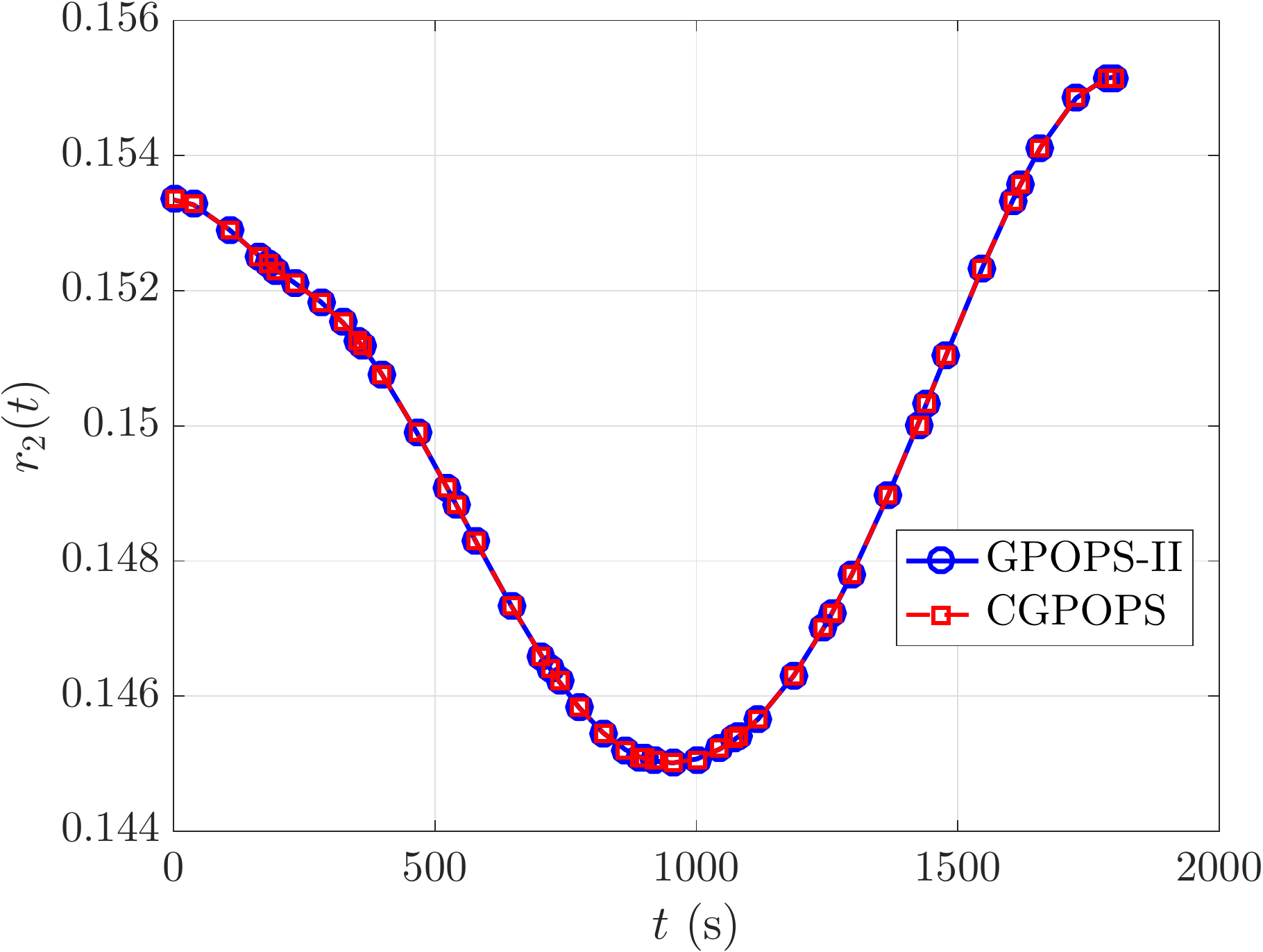,height=2in}}~\subfloat[$r_3(t)$ vs.~$t$. \label{fig:spaceStationState6}]{\epsfig{figure=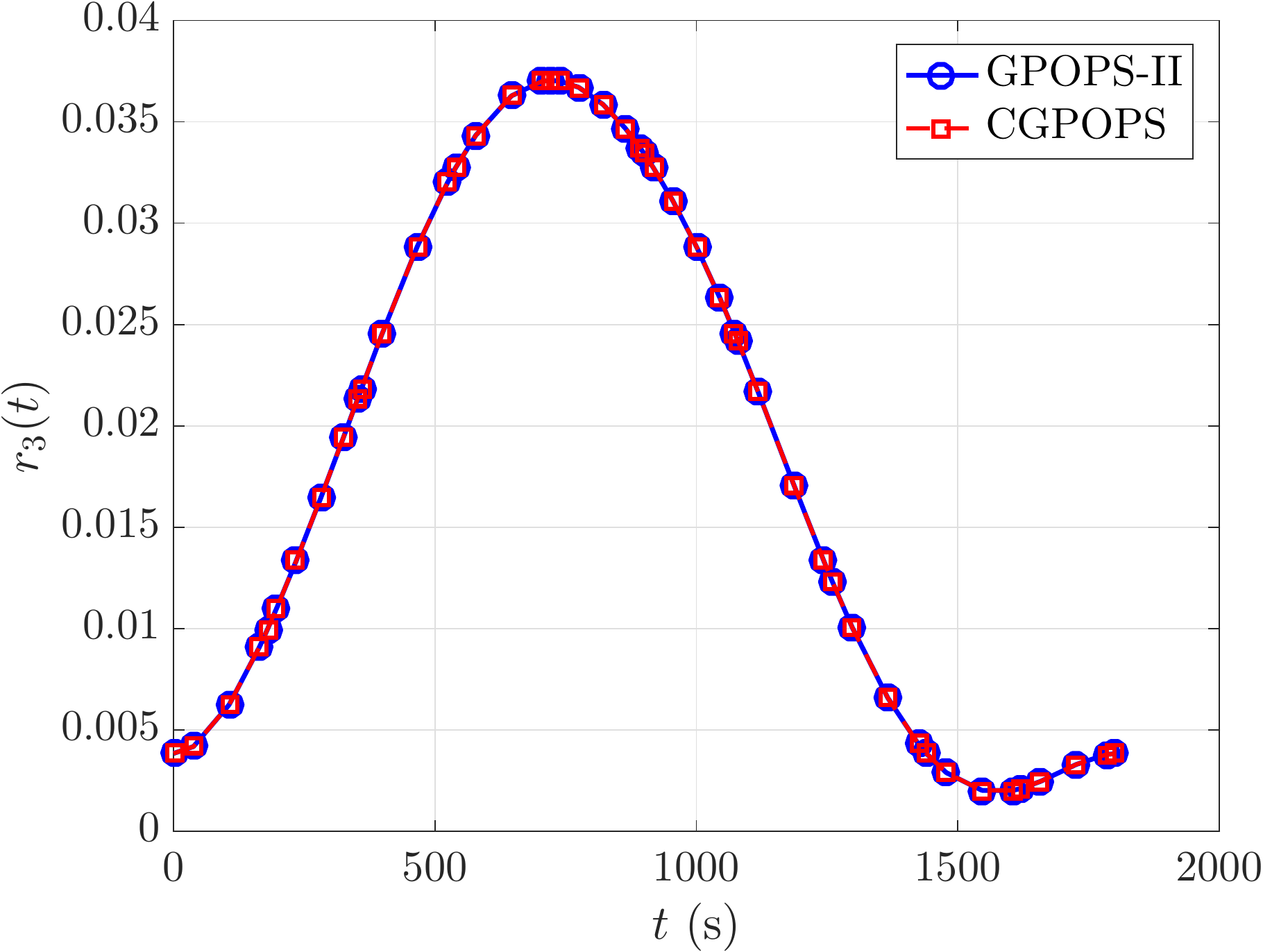,height=2in}}
  
  \caption{$\CGPOPS$ and $\mathbb{GPOPS-II}$ Solutions to Example 3 Using with the NLP Solver IPOPT and a Mesh  Refinement Tolerance of $10^{-6}$ using $hp$-I(4,10). \label{fig:spaceStationStateSolution}}
\end{figure}

\begin{figure}[ht!]
\centering

\subfloat[$h_1(t)$ vs.~$t$.\label{fig:spaceStationState7}]{\epsfig{figure=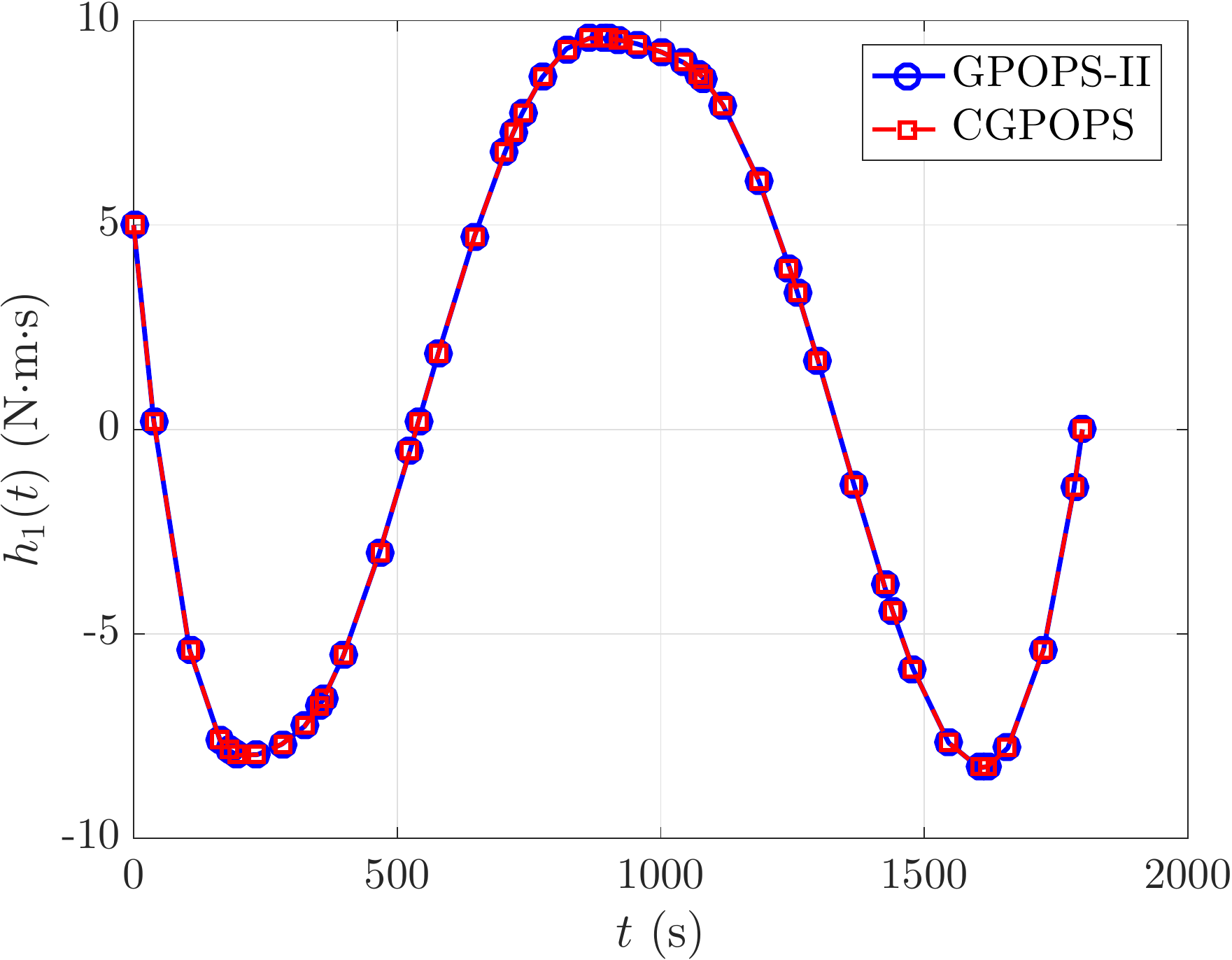,height=2in}}~\subfloat[$h_2(t)$ vs.~$t$. \label{fig:spaceStationState8}]{\epsfig{figure=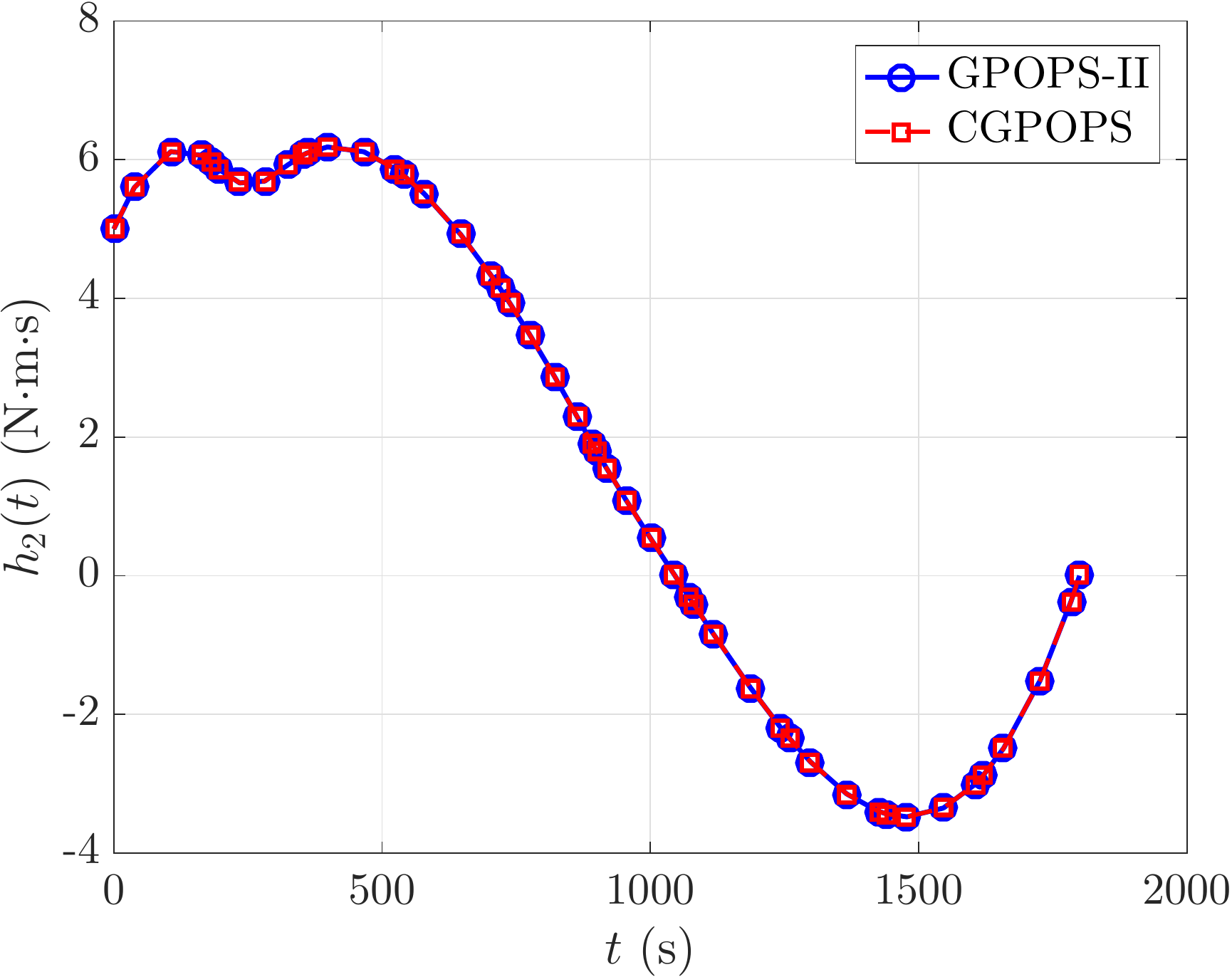,height=2in}}

\subfloat[$h_3(t)$ vs.~$t$.\label{fig:spaceStationState9}]{\epsfig{figure=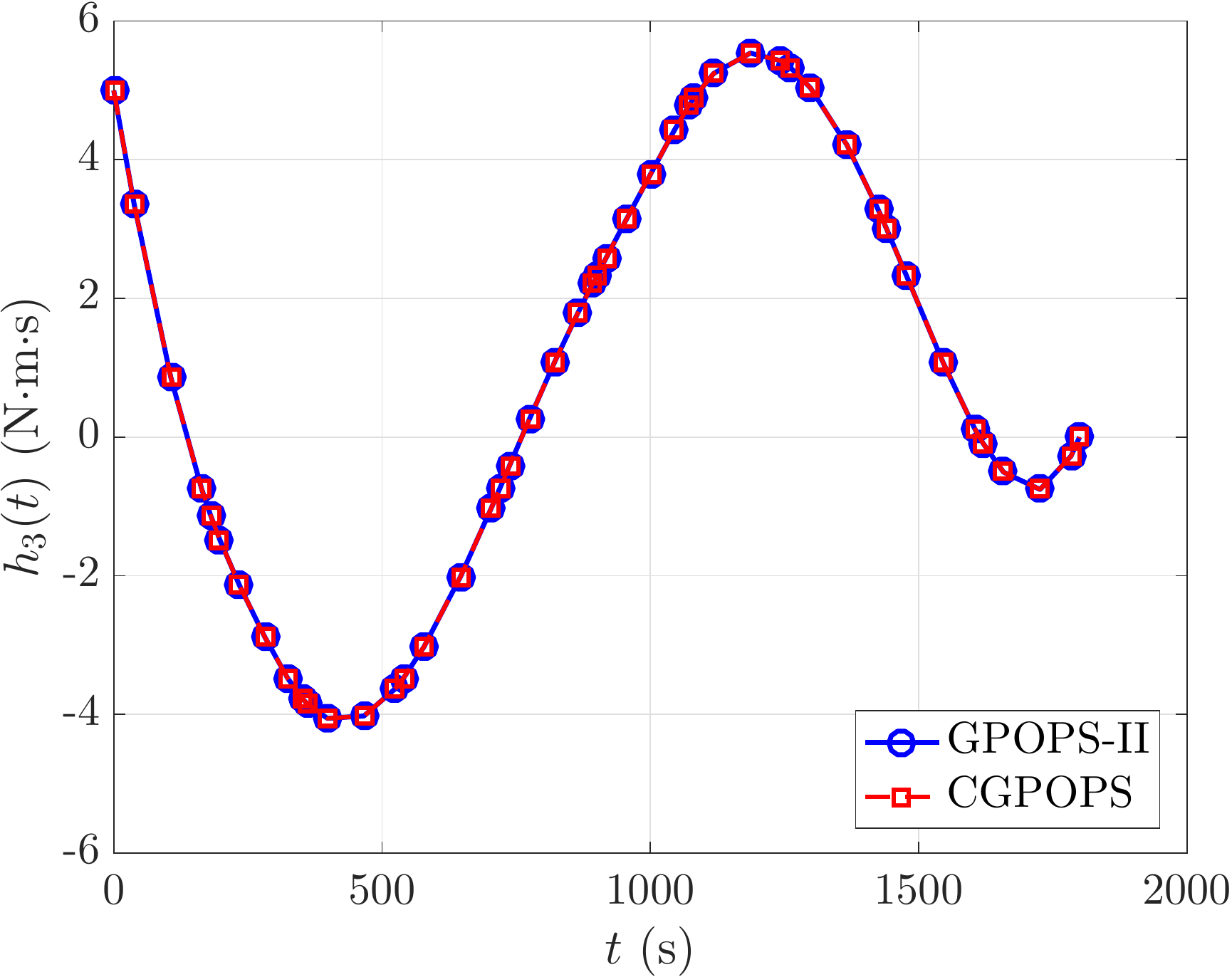,height=2in}}~\subfloat[$u_1(t)$ vs.~$t$. \label{fig:spaceStationControl1}]{\epsfig{figure=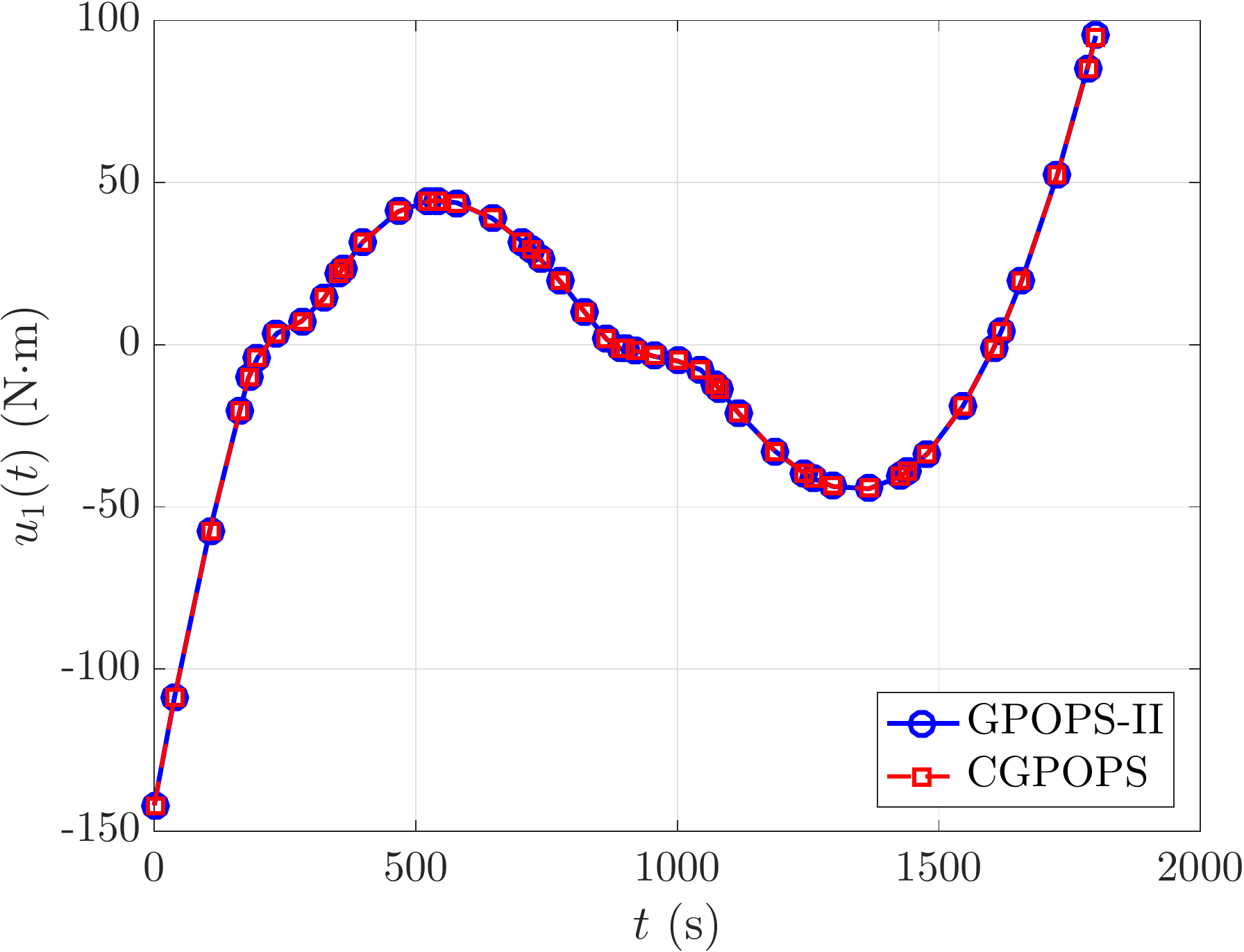,height=2in}}

\subfloat[$u_2(t)$ vs.~$t$.\label{fig:spaceStationControl2}]{\epsfig{figure=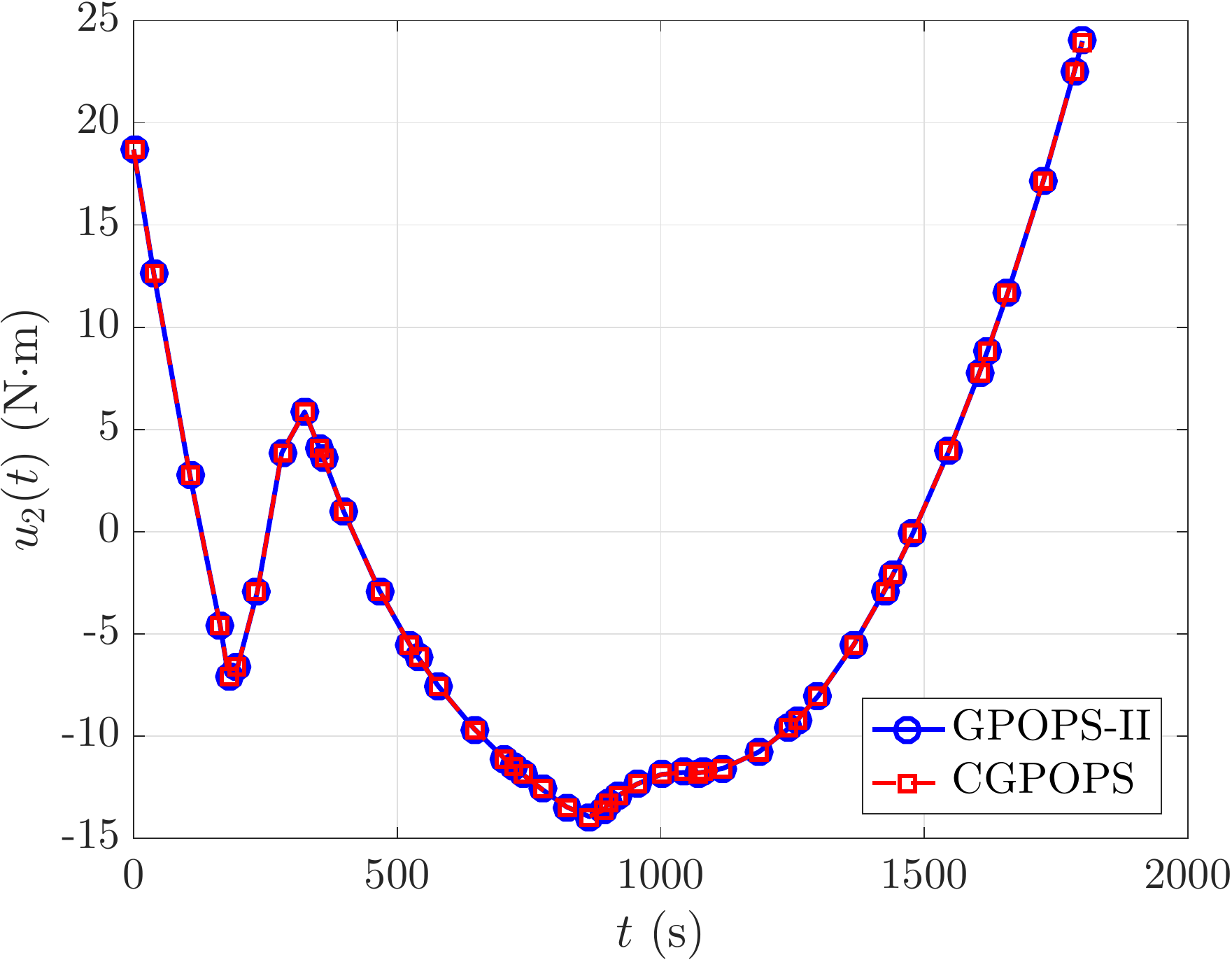,height=2in}}~\subfloat[$u_3(t)$ vs.~$t$. \label{fig:spaceStationControl3}]{\epsfig{figure=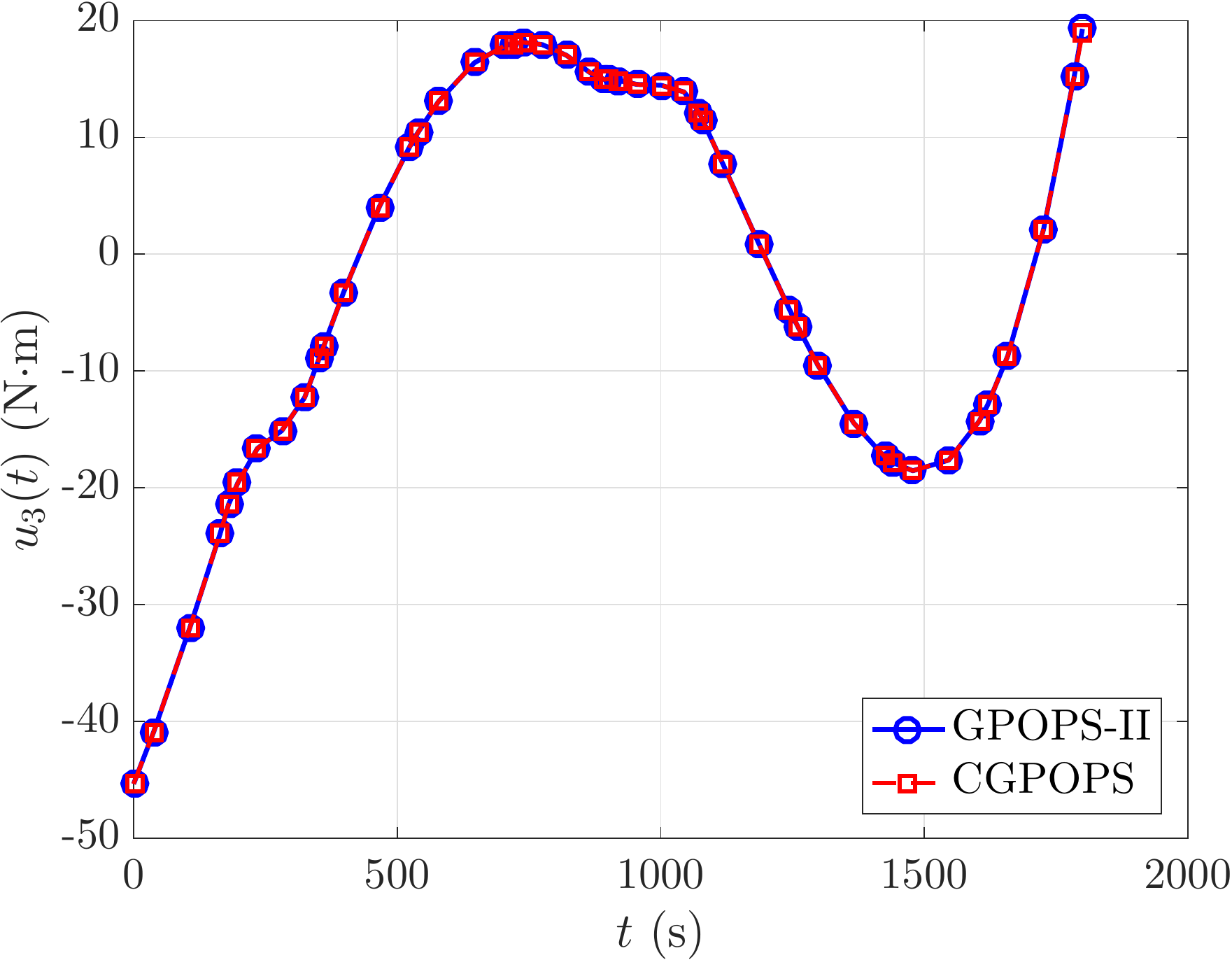,height=2in}}
  
  \caption{$\CGPOPS$ and $\mathbb{GPOPS-II}$ Solutions to Example 3 Using with the NLP Solver IPOPT and a Mesh  Refinement Tolerance of $10^{-6}$ using $hp$-I(4,10). \label{fig:spaceStationControlSolution}}
\end{figure}

The optimal control problem defined by Eqs.~\eqref{eq:space-station-cost} -- \eqref{eq:space-station-IC} was solved using $\CGPOPS$ with the $hp$-I(4,10), $hp$-II(4,10), $hp$-III(4,10), and $hp$-IV(4,10) methods on an initial mesh consisting of ten uniformly spaced mesh intervals and four LGR points per mesh interval.  The NLP solver and mesh refinement accuracy tolerances were set to $10^{-7}$ and $10^{-6}$, respectively.  The initial guess was a constant over the time interval $t\in[0,1800]$, where the constant was $(\bar{\g{\omega}}_0,\bar{\m{r}}_0,\bar{\m{h}}_0)$ for the state and zero for the control.  The state and control solutions obtained using $\CGPOPS$ are shown, respectively, in Fig.~\ref{fig:spaceStationStateSolution} and \ref{fig:spaceStationControlSolution} alongside the solution obtained using the optimal control software $\mathbb{GPOPS-II}$ \cite{Patterson2014} with the $hp$-I(4,10).  It is seen that the $\CGPOPS$ solution matches extremely well with the $\mathbb{GPOPS-II}$ solution.  Moreover, the optimal cost obtained using $\CGPOPS$ and $\mathbb{GPOPS-II}$ are essentially the same, with values of $3.5867511\times10^{-6}$ and $3.5867511\times10^{-6}$, respectively.  Finally, the computational time required by $\CGPOPS$ and $\mathbb{GPOPS-II}$ to solve the optimal control problem was $0.5338$ seconds and $2.7696$ seconds, respectively.

\begin{table}[htp]
\centering
\caption{Performance of $\CGPOPS$ on Example 3 using $hp$-I(4,10). \label{tab:spaceStation-hpPatterson}}
\footnotesize
\begin{tabular}{|c|c|c|c|c|} \hline
{\bf Mesh} & {\bf Estimated}&{\bf Number of} &{\bf Estimated }&{\bf Number of} \\
{\bf Iteration} & {\bf Error ($\CGPOPS$)}&{\bf Collocation}&{\bf Error ($\mathbb{GPOPS-II}$)}&{\bf Collocation} \\
{\bf Number} & {\bf $hp$-I(4,10)}&{\bf Points}&{\bf $hp$-I(4,10)}&{\bf Points} \\\hline\hline
$1$& $9.409\times 10^{-6}$ & $41$ &$9.409\times 10^{-6}$ & $41$  \\ \hline
$2$ & $6.496\times 10^{-7}$ & $47$ &$6.496\times 10^{-7}$ & $47$ \\ \hline
\end{tabular}
\normalsize
\end{table}

\begin{table}[htp]
\centering
\caption{Performance of $\CGPOPS$ on Example 3 using $hp$-II(4,10). \label{tab:spaceStation-hpDarby}}
\footnotesize
\begin{tabular}{|c|c|c|c|c|} \hline
{\bf Mesh} & {\bf Estimated}&{\bf Number of} &{\bf Estimated }&{\bf Number of} \\
{\bf Iteration} & {\bf Error ($\CGPOPS$)}&{\bf Collocation}&{\bf Error ($\mathbb{GPOPS-II}$)}&{\bf Collocation} \\
{\bf Number} & {\bf $hp$-II(4,10)}&{\bf Points}&{\bf $hp$-II(4,10)}&{\bf Points} \\\hline\hline
$1$& $9.409\times 10^{-6}$ & $41$ &$9.409\times 10^{-6}$ & $41$  \\ \hline
$2$ & $2.389\times 10^{-6}$ & $50$ &$2.387\times 10^{-6}$ & $50$ \\ \hline
$3$ & $7.125\times 10^{-7}$ & $55$ &$7.130\times 10^{-7}$ & $55$ \\ \hline
\end{tabular}
\normalsize
\end{table}
    
\begin{table}[htp]
\centering
\caption{Performance of $\CGPOPS$ on Example 3 using $hp$-III(4,10). \label{tab:spaceStation-hpLiu}}
\footnotesize
\begin{tabular}{|c|c|c|c|c|} \hline
{\bf Mesh} & {\bf Estimated}&{\bf Number of} &{\bf Estimated }&{\bf Number of} \\
{\bf Iteration} & {\bf Error ($\CGPOPS$)}&{\bf Collocation}&{\bf Error ($\mathbb{GPOPS-II}$)}&{\bf Collocation} \\
{\bf Number} & {\bf $hp$-III(4,10)}&{\bf Points}&{\bf $hp$-III(4,10)}&{\bf Points} \\\hline\hline
$1$& $9.409\times 10^{-6}$ & $41$ &$9.409\times 10^{-6}$ & $41$  \\ \hline
$2$ & $9.542\times 10^{-7}$ & $50$ &$9.559\times 10^{-7}$ & $50$ \\ \hline
\end{tabular}
\normalsize
\end{table}
    
\begin{table}[htp]
\centering
\caption{Performance of $\CGPOPS$ on Example 3 using $hp$-IV(4,10). \label{tab:spaceStation-hpLegendre}}
\footnotesize
\begin{tabular}{|c|c|c|c|c|} \hline
{\bf Mesh} & {\bf Estimated}&{\bf Number of} &{\bf Estimated }&{\bf Number of} \\
{\bf Iteration} & {\bf Error ($\CGPOPS$)}&{\bf Collocation}&{\bf Error ($\mathbb{GPOPS-II}$)}&{\bf Collocation} \\
{\bf Number} & {\bf $hp$-IV(4,10)}&{\bf Points}&{\bf $hp$-IV(4,10)}&{\bf Points} \\\hline\hline
$1$& $9.409\times 10^{-6}$ & $41$ &$9.409\times 10^{-6}$ & $41$  \\ \hline
$2$ & $1.049\times 10^{-7}$ & $53$ &$1.046\times 10^{-7}$ & $53$ \\ \hline
$3$ & $7.125\times 10^{-7}$ & $57$ &$7.130\times 10^{-7}$ & $57$ \\ \hline
\end{tabular}
\normalsize
\end{table}

\clearpage


\subsection{Example 4: Free-Flying Robot Problem\label{subsect:freeFlyingRobot}}


Consider the following {\em free-flying robot} optimal control problem taken from Refs.~\cite{Betts3} and \cite{Sakawa}.  The objective is to minimize the cost functional  
\begin{equation}\label{eq:freeFlyingRobot-Cost}
  \cal{J} = \int_{0}^{t_f}(u_1+u_2+u_3+u_4)dt~,
\end{equation}
subject to the dynamic constraints
\begin{equation}\label{eq:freeFlyingRobot-Dynamics}
 \begin{array}{lclclcl}
   \dot{x} & = & v_x~, & & \dot{y} & = & v_y~, \\
   \dot{v_x} & = & (F_1+F_2)\cos(\theta)~, & & \dot{v_y} & = & (F_1+F_2)\sin(\theta)~, \\
   \dot{\theta} & = & \omega~, & & \dot{\omega} & = & \alpha F_1-\beta F_2~, \\
\end{array}
\end{equation}
the control inequality constraints
\begin{equation}\label{eq:freeFlyingRobot-Constraint}
\begin{array}{lcl}
 0 \leq u_i \leq 1~, \quad (i=1,2,3,4)~, & & \quad F_i \leq 1~,~\quad (i=1,2)~,
\end{array}  
\end{equation}
and the boundary conditions
\begin{equation} \label{eq:freeFlyingRobot-BCs}
\begin{array}{lclclcl}
 x(0) = -10~, & & x(t_f) = 0~, & &
 y(0) = -10~, & & y(t_f) = 0~, \\
 v_x(0) = 0~, & & v_x(t_f) = 0~, & &
 v_y(0) = 0~, & & v_y(t_f) = 0~, \\
 \theta(0) = \frac{\pi}{2}~, & & \theta(t_f) = 0~, & &
 \omega(0) = 0~, & & \omega(t_f) = 0~,
\end{array}
\end{equation} 
where 
\begin{equation}
\label{eq:freeFlyingRobot-Aux}
 \begin{array}{lclclcl}
 F_1 = u_1-u_2~, & & F_2 = u_3-u_4~,& &
 \alpha = 0.2~, & & \beta = 0.2~.
 \end{array}
\end{equation}
It is known that the optimal control problem defined by Eqs.~\eqref{eq:freeFlyingRobot-Cost} -- \eqref{eq:freeFlyingRobot-Aux} is a bang-bang optimal control problem, as the Hamiltonian is linear with respect to the control.  Given the structure of the solution, the $hp$-BB(3,10) mesh refinement method \cite{Agamawi2019} is also employed to solve this example.
 
\begin{figure}[htp]
\centering

\epsfig{file=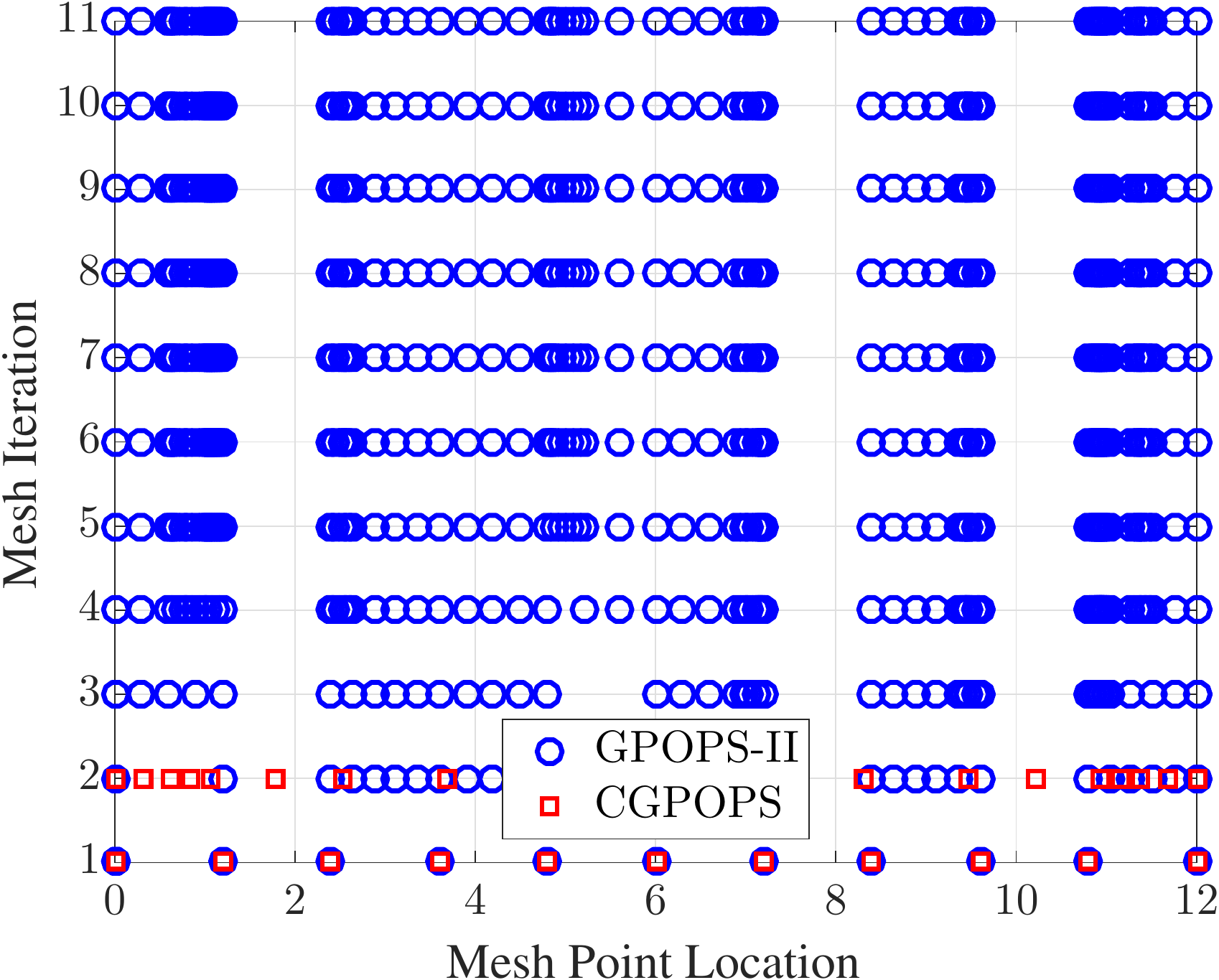,height=2in}
  
\caption{$\CGPOPS$ and $\mathbb{GPOPS-II}$ Mesh Refinement History for Example 4 Using $hp$-BB(3,10) and $hp$-II(3,10), respectively. \label{fig:freeFlyingRobotMeshRefinement}}
\end{figure}

\begin{figure}[htp]
\centering

\subfloat[$x(t)$ vs.~$t$.\label{fig:freeFlyingRobotStateX}]{\epsfig{figure=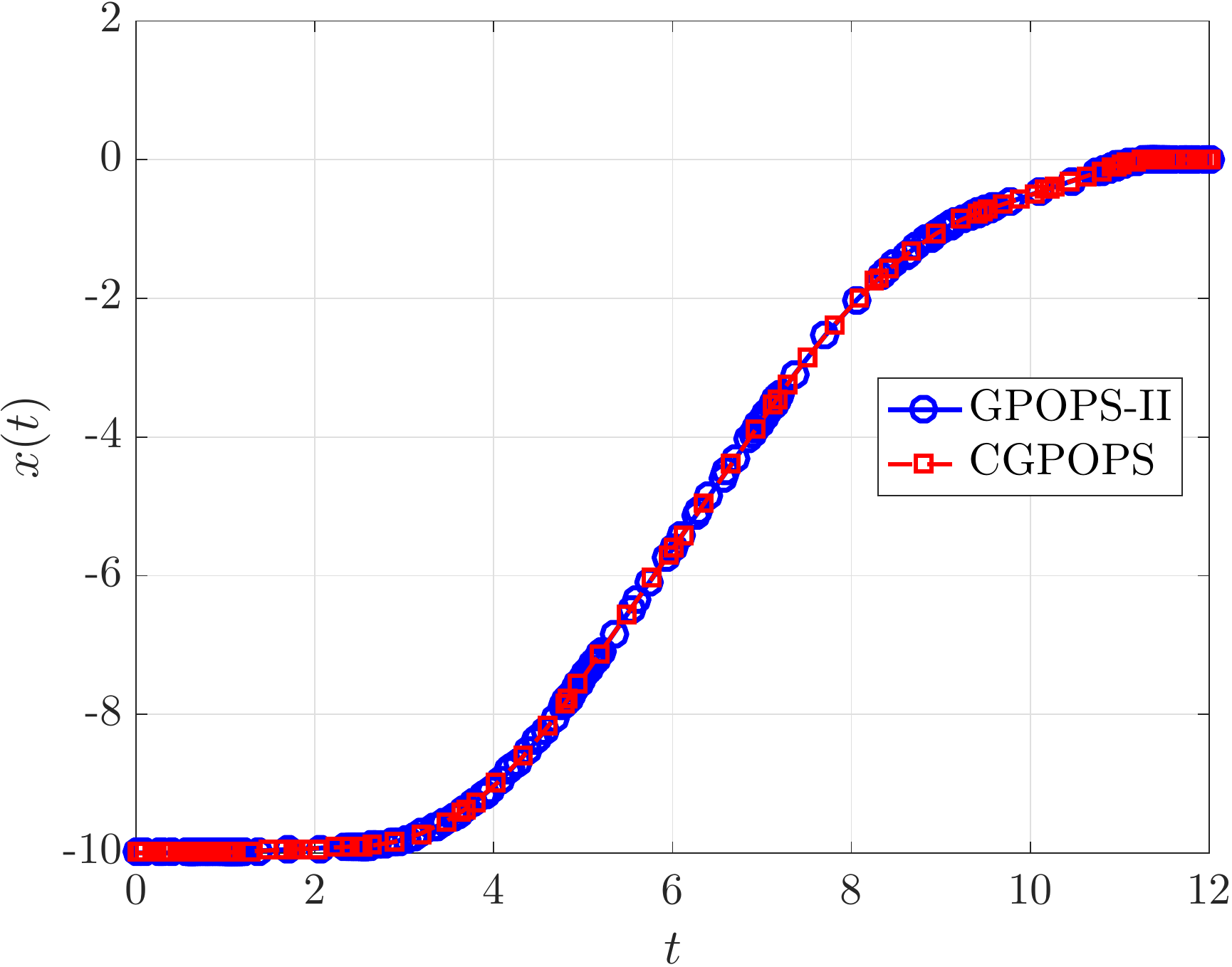,height=2in}}~~~\subfloat[$y(t)$ vs.~$t$. \label{fig:freeFlyingRobotStateY}]{\epsfig{figure=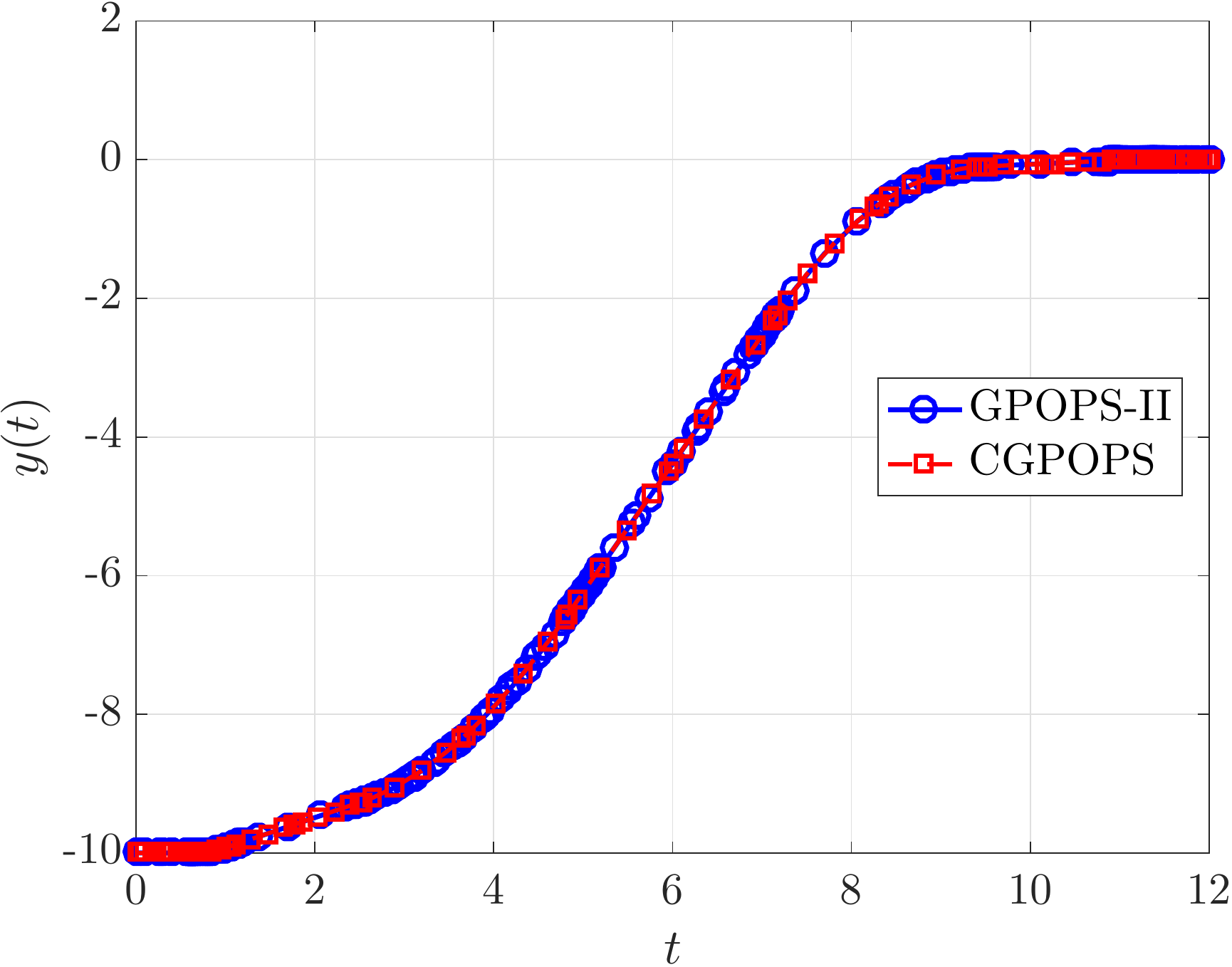,height=2in}}

\subfloat[$v_x(t)$ vs.~$t$.\label{fig:freeFlyingRobotStateVx}]{\epsfig{figure=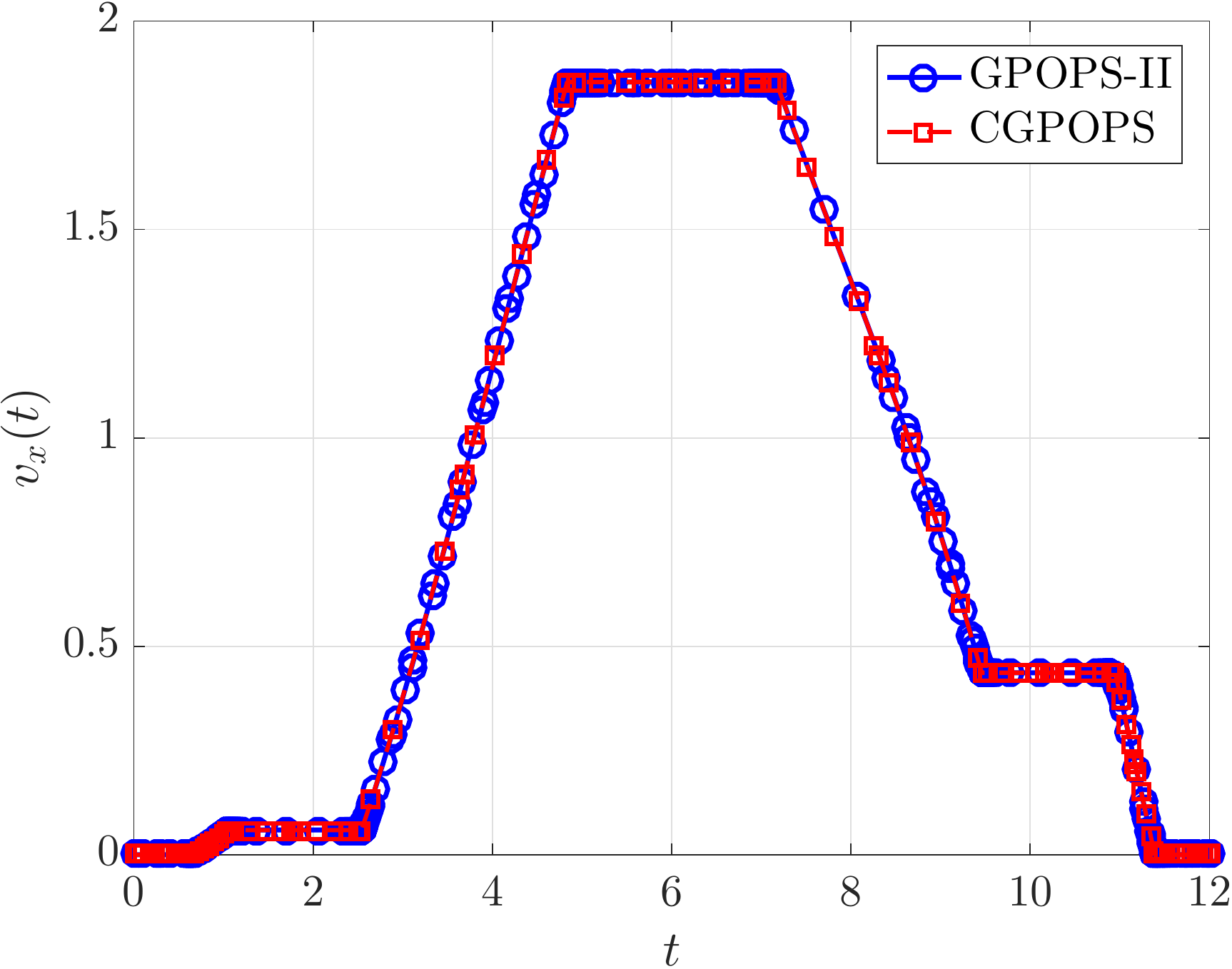,height=2in}}~~~\subfloat[$v_y(t)$ vs.~$t$. \label{fig:freeFlyingRobotStateVy}]{\epsfig{figure=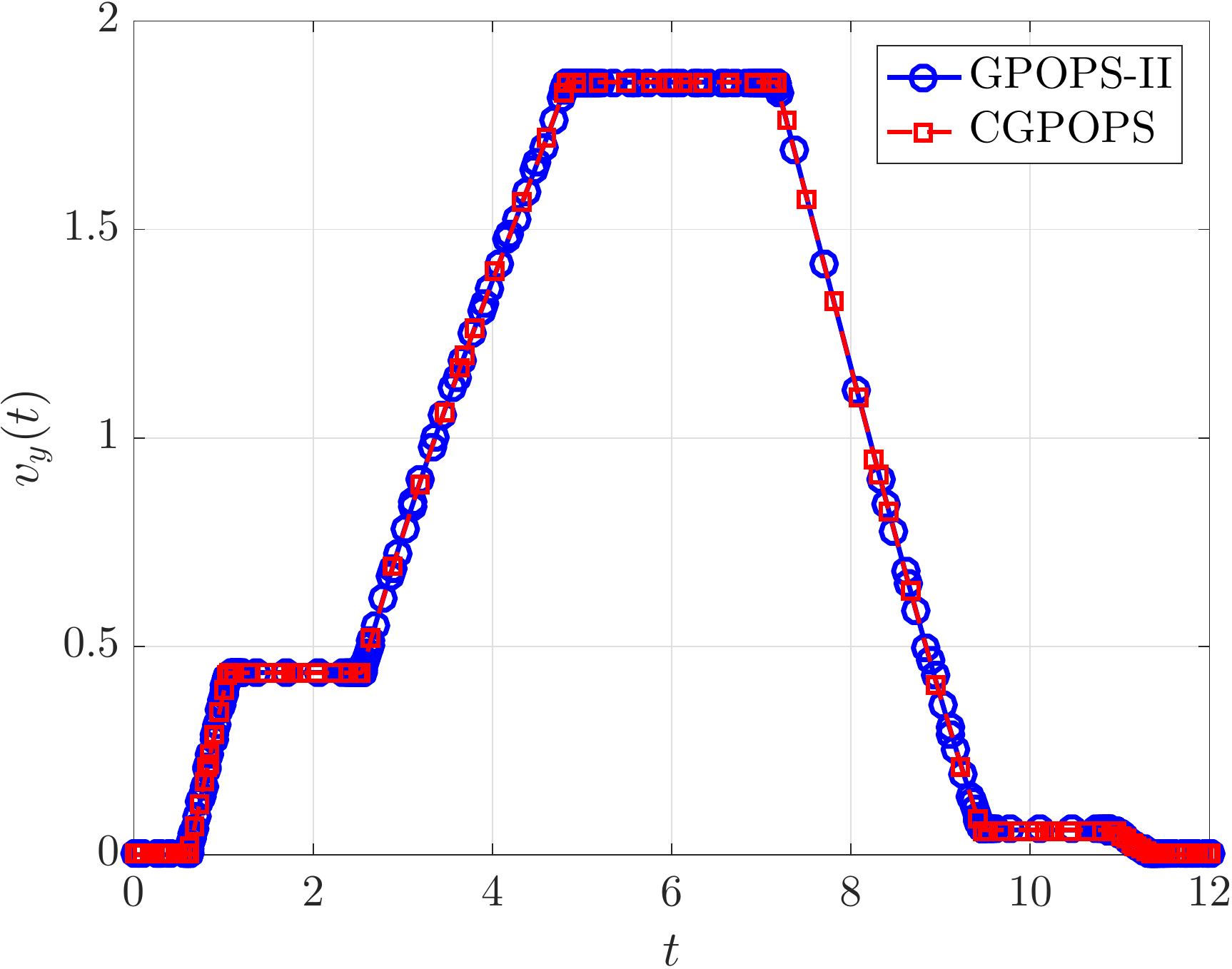,height=2in}}

\subfloat[$\theta(t)$ vs.~$t$.\label{fig:freeFlyingRobotStateTheta}]{\epsfig{figure=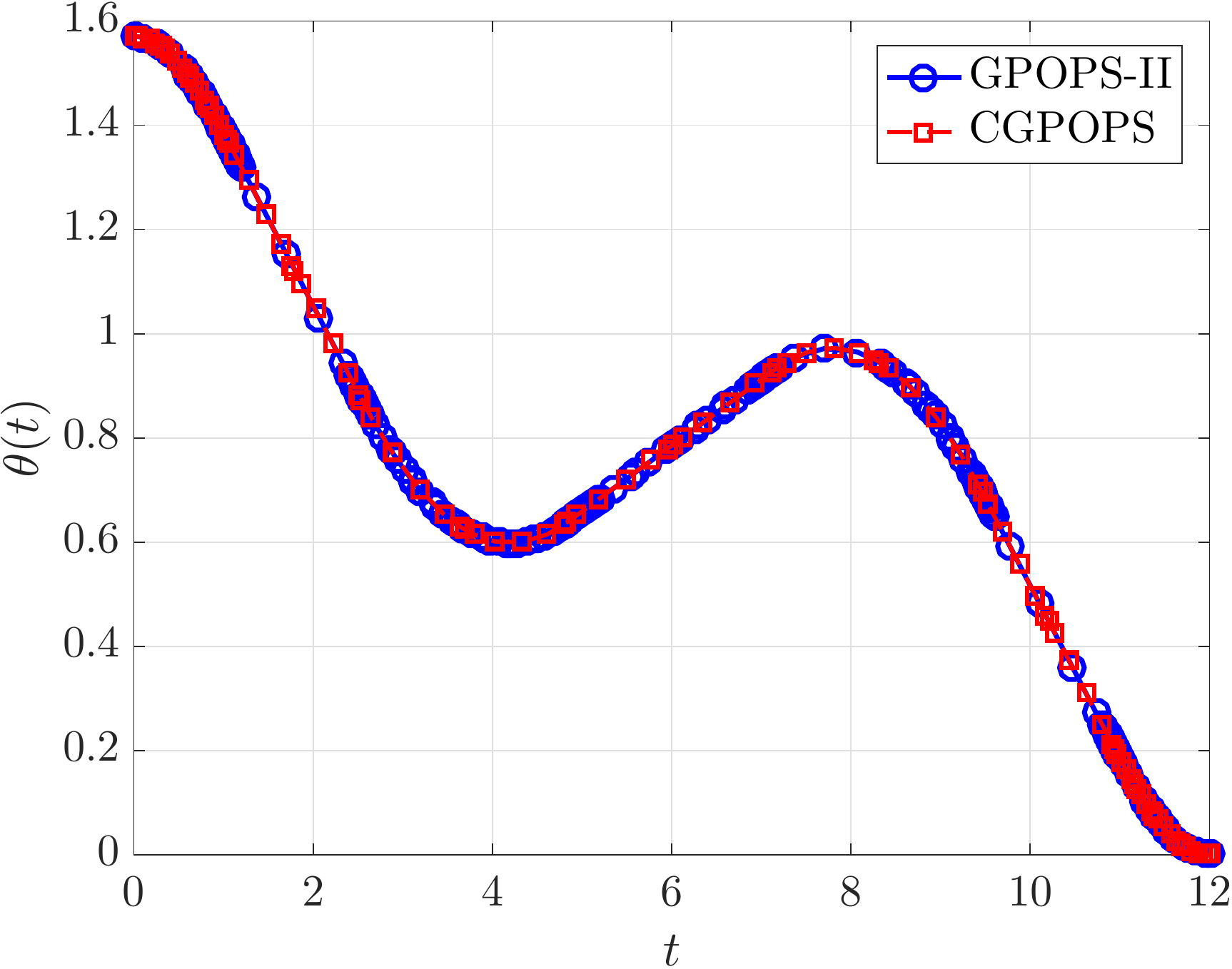,height=2in}}~~~\subfloat[$\omega(t)$ vs.~$t$. \label{fig:freeFlyingRobotStateOmega}]{\epsfig{figure=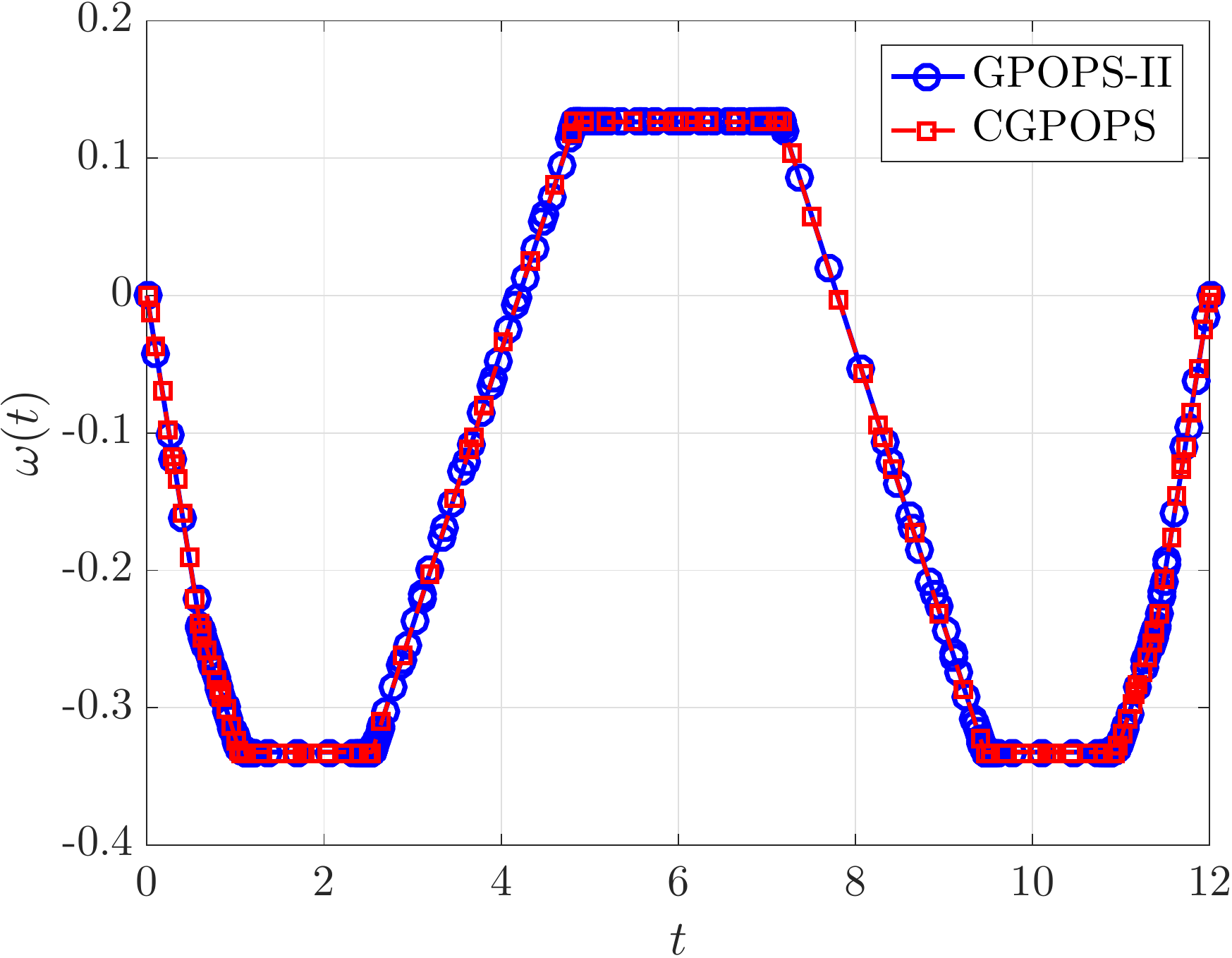,height=2in}}
\caption{$\CGPOPS$ and $\mathbb{GPOPS-II}$ State Solutions to Example 4 using $hp$-BB(3,10) and $hp$-II(3,10), respectively. \label{fig:freeFlyingRobotStateSolution}}
\end{figure}

\begin{figure}[htp]
\centering

\subfloat[$u_1(t)$ vs.~$t$.\label{fig:freeFlyingRobotControlU1}]{\epsfig{figure=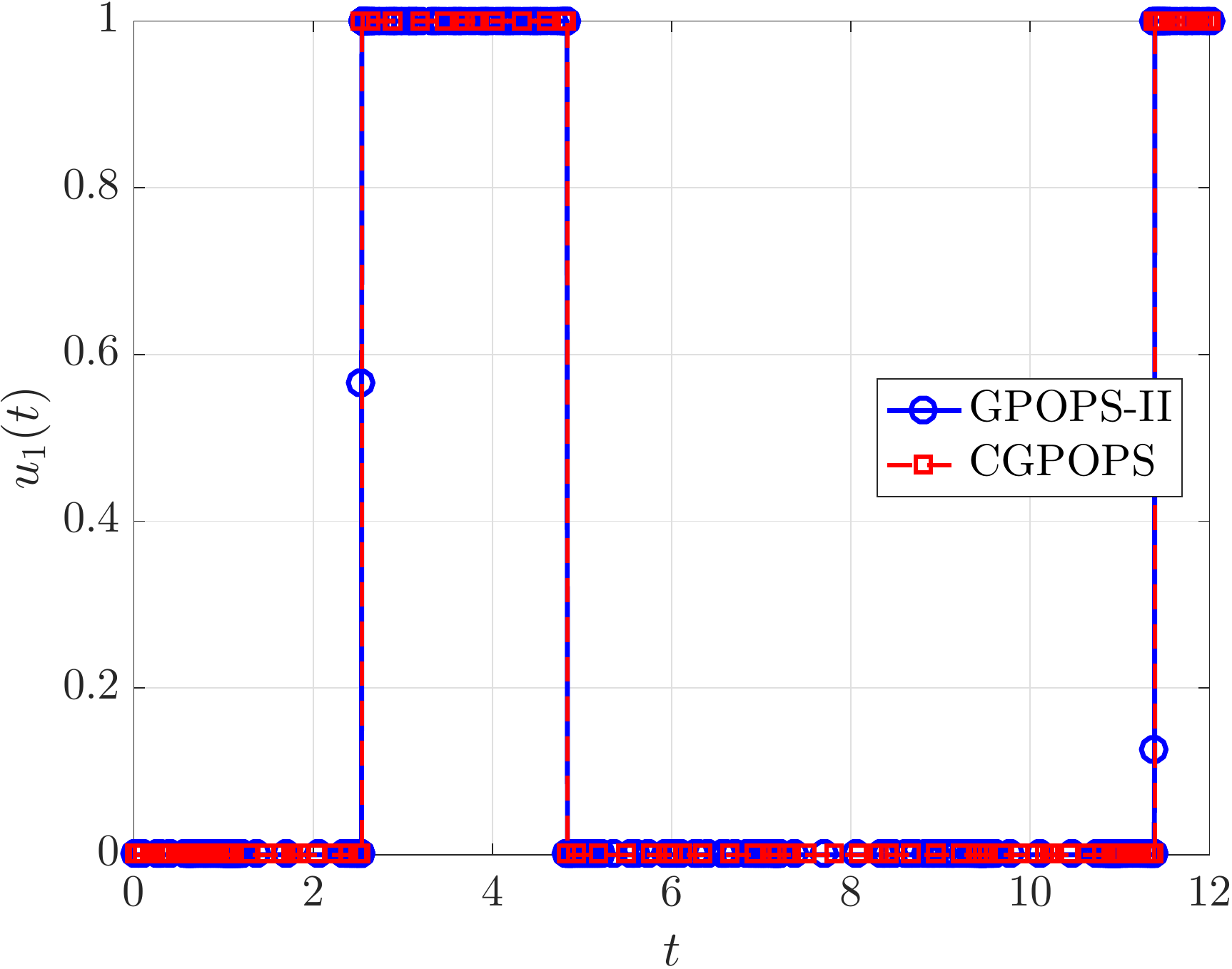,height=2in}}~~~\subfloat[$u_2(t)$ vs.~$t$. \label{fig:freeFlyingRobotControlU2}]{\epsfig{figure=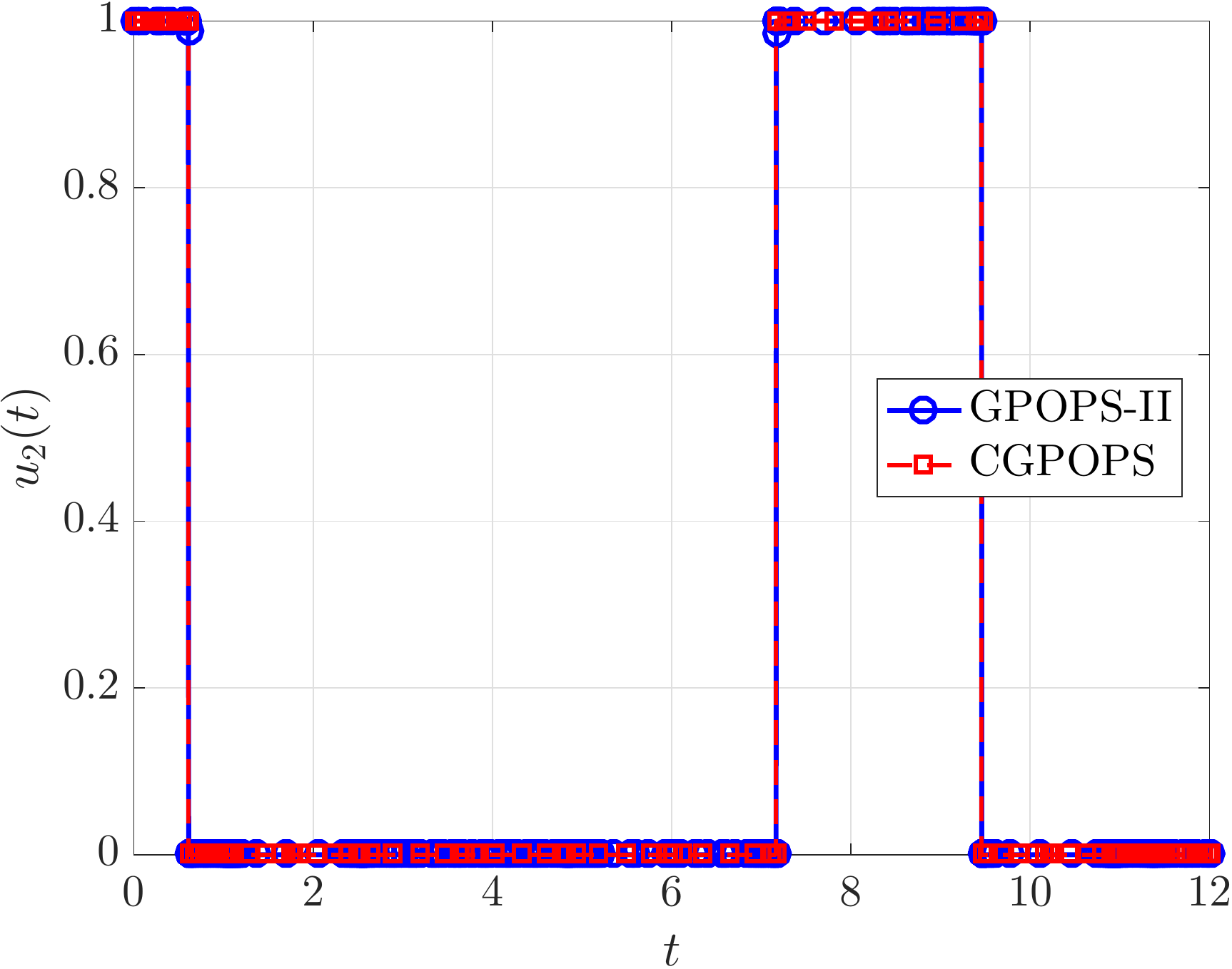,height=2in}}

\subfloat[$u_3(t)$ vs.~$t$.\label{fig:freeFlyingRobotControlU3}]{\epsfig{figure=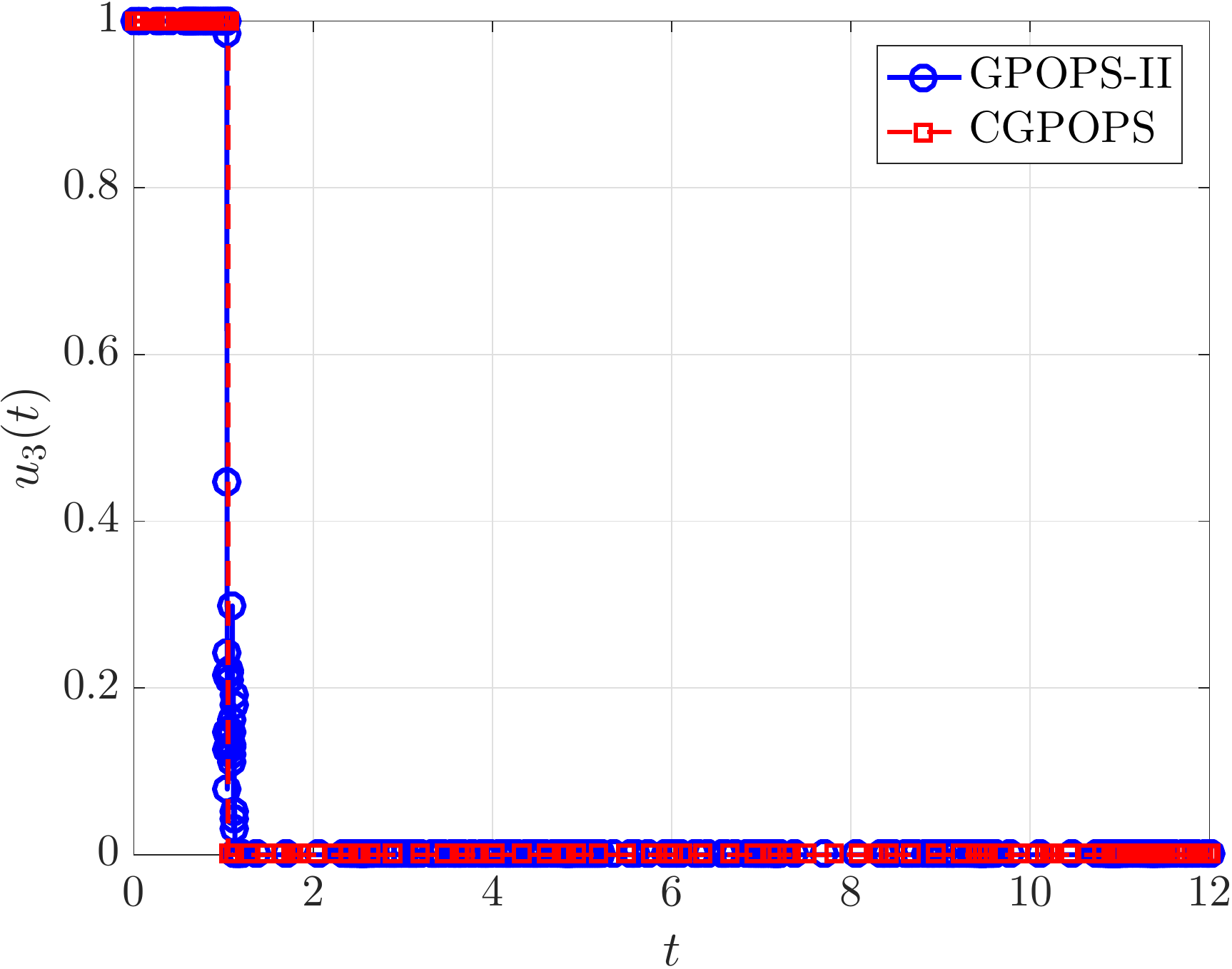,height=2in}}~~~\subfloat[$u_4(t)$ vs.~$t$. \label{fig:freeFlyingRobotControlU4}]{\epsfig{figure=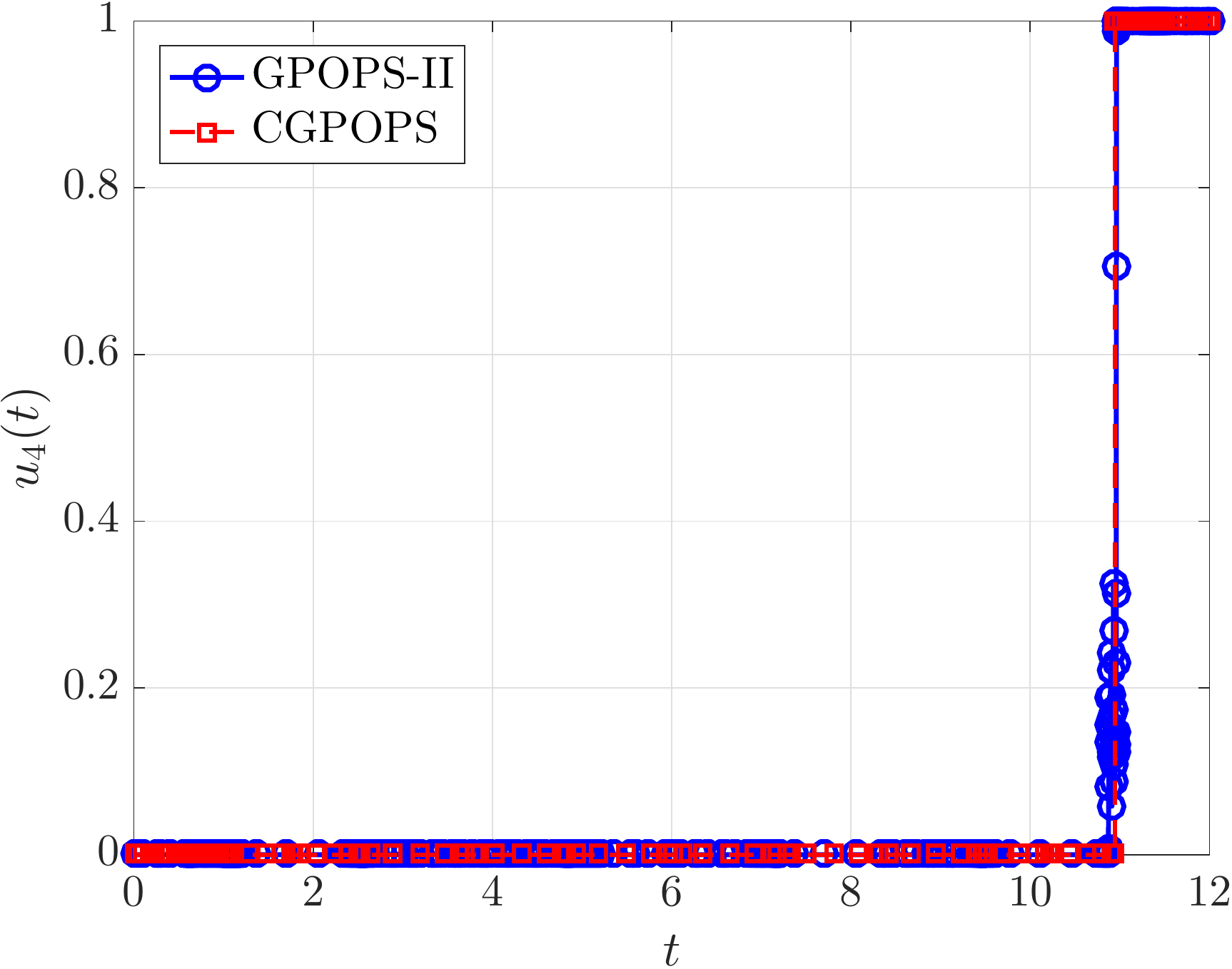,height=2in}}

\subfloat[$F_1(t)$ vs.~$t$.\label{fig:freeFlyingRobotControlF1}]{\epsfig{figure=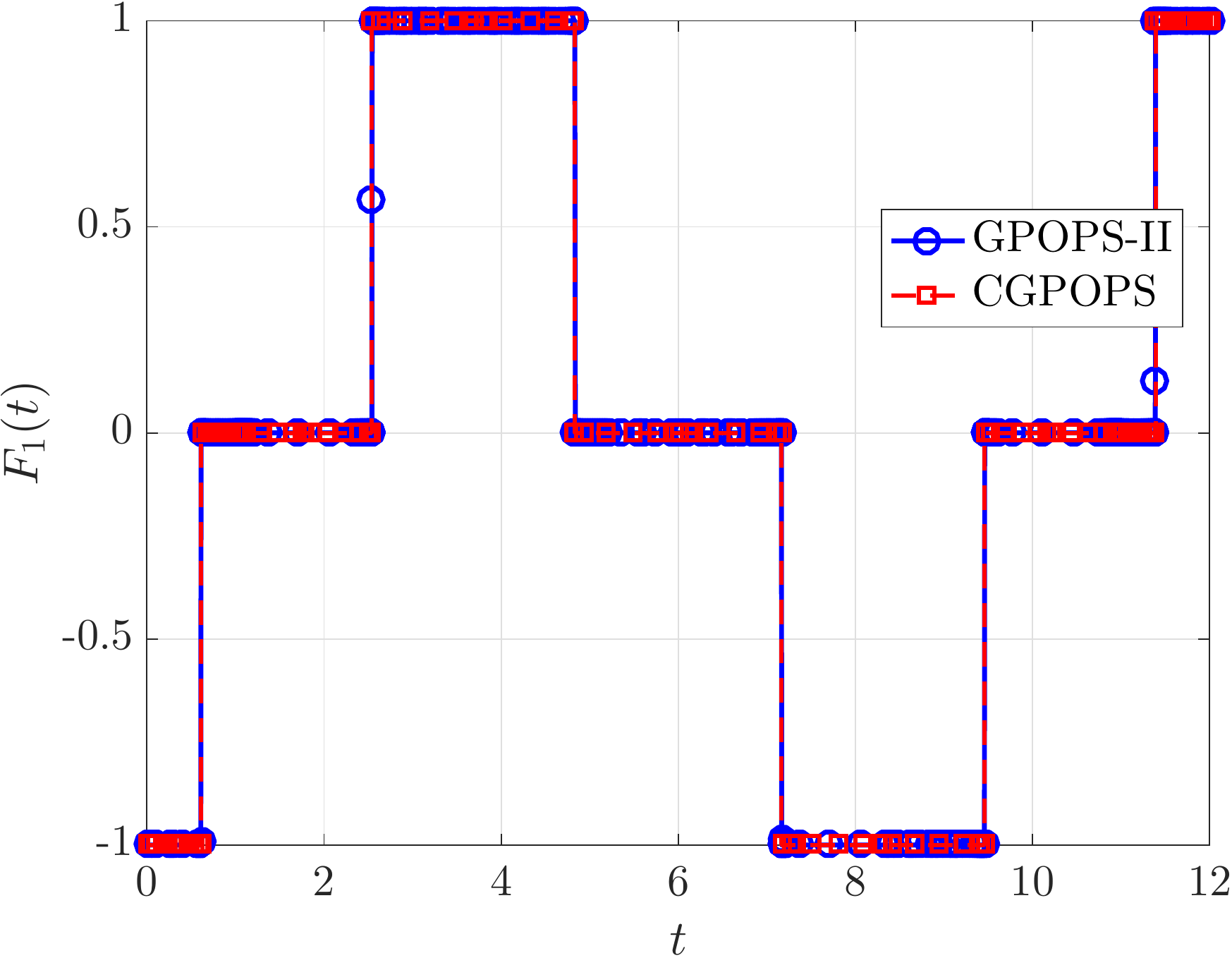,height=2in}}~~~\subfloat[$F_2(t)$ vs.~$t$. \label{fig:freeFlyingRobotControlF2}]{\epsfig{figure=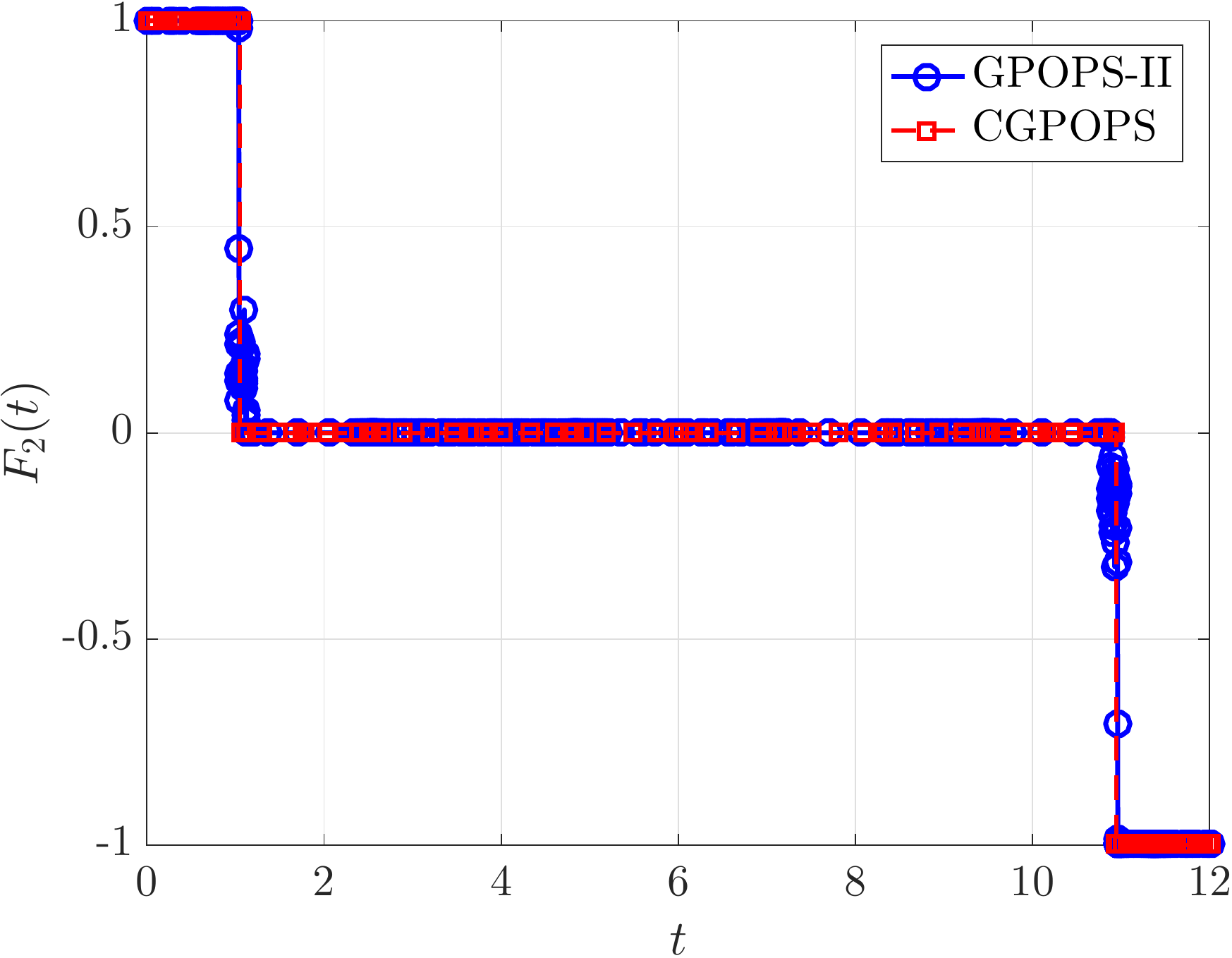,height=2in}}
\caption{$\CGPOPS$ and $\mathbb{GPOPS-II}$ Control Solutions to Example 4 using $hp$-BB(3,10) and $hp$-II(3,10), respectively. \label{fig:freeFlyingRobotControlSolution}}
\end{figure}

The free-flying robot optimal control problem was solved using $\CGPOPS$ with the mesh refinement methods $hp$-I(3,10), $hp$-II(3,10), $hp$-III(3,10), $hp$-IV(3,10), and $hp$-BB(3,10) on an initial mesh of ten evenly spaced mesh intervals with five LGR points per mesh interval.  The NLP solver and mesh refinement accuracy tolerances were set to $10^{-9}$ and $10^{-7}$, respectively.  The solution obtained using $\CGPOPS$ with the $hp$-BB(3,10) method is shown in Figs.~\ref{fig:freeFlyingRobotStateSolution} and \ref{fig:freeFlyingRobotControlSolution} alongside the solution obtained with $\mathbb{GPOPS-II}$ \cite{Patterson2014} with the $hp$-II(3,10) method.  It is seen that the $\CGPOPS$ and $\mathbb{GPOPS-II}$ solutions are in excellent agreement.  Moreover, the optimal cost obtained using $\CGPOPS$ and $\mathbb{GPOPS-II}$ are extremely close, with values of $7.9101471$ and $7.9101421$, respectively, agreeing to the fifth decimal place.  Additionally, the computational time required by $\CGPOPS$ and $\mathbb{GPOPS-II}$ to solve the optimal control problem was $0.6313$ seconds and $9.1826$ seconds, respectively.  In order to demonstrate how $\CGPOPS$ is capable of more effectively capturing the bang-bang control profile of the optimal solution, Fig.~\ref{fig:freeFlyingRobotControlSolution} shows the control solutions obtained using $\CGPOPS$ and $\mathbb{GPOPS-II}$ with their most effective mesh refinement methods for this example, $hp$-BB(3,10) and $hp$-II(3,10), respectively.  It is seen that $\CGPOPS$ exactly captures the switching points for all eight discontinuities, while the solution obtained using $\mathbb{GPOPS-II}$ has fluctuations in the control near the discontinuities for the third and fourth components, as may be observed from Figs.~\ref{fig:freeFlyingRobotControlU3},\ref{fig:freeFlyingRobotControlU4}, and \ref{fig:freeFlyingRobotControlF2}.  Additionally, Fig.~\ref{fig:freeFlyingRobotMeshRefinement} shows the mesh refinement history for $\CGPOPS$ with the $hp$-BB(3,10) method and $\mathbb{GPOPS-II}$ with the $hp$-II(3,10) method where $\CGPOPS$ only requires a single mesh refinement iteration to solve, while $\mathbb{GPOPS-II}$ takes nine mesh refinement iterations to solve.  Finally, Tables~\ref{tab:freeFlyingRobotMeshRefinement-hpPatterson} -- \ref{tab:freeFlyingRobotMeshRefinement-hpBangBang} show the estimated error on each mesh, where it is seen that the solution error decreases steadily with each mesh refinement iteration for all $hp$ methods employed.
 
\begin{table}[htp]
\centering
\caption{Performance of $\CGPOPS$ on Example 4 using $hp$-I(3,10). \label{tab:freeFlyingRobotMeshRefinement-hpPatterson}}
\footnotesize
\begin{tabular}{|c|c|c|c|c|} \hline
{\bf Mesh} & {\bf Estimated}&{\bf Number of} &{\bf Estimated }&{\bf Number of} \\
{\bf Iteration} & {\bf Error ($\CGPOPS$)}&{\bf Collocation}&{\bf Error ($\mathbb{GPOPS-II}$)}&{\bf Collocation} \\
{\bf Number} & {\bf $hp$-I(3,10)}&{\bf Points}&{\bf $hp$-I(3,10)}&{\bf Points} \\\hline\hline
$1$	& $5.7636 \times 10^{-4}$	& $50$ 	&	$5.7636 \times 10^{-4}$	& $50$		\\ \hline
$2$	& $2.3428 \times 10^{-4}$	& $82$ 	&	$1.2977 \times 10^{-4}$ 	& $82$  	\\ \hline
$3$	& $7.5065 \times 10^{-5}$	& $122$	&	$2.3256 \times 10^{-4}$	& $120$ 	\\ \hline
$4$	& $6.2091 \times 10^{-5}$		&$157$ 	&	$1.1175 \times 10^{-5}$	&$161$ 	\\ \hline
$5$	& $9.4236 \times 10^{-6}$	& $184$ 	& $6.2093 \times 10^{-5}$	& $188$  	\\ \hline
$6$	& $3.9835 \times 10^{-6}$	& $209$ 	& $4.8405 \times 10^{-6}$	& $212$  	\\ \hline
$7$	& $2.8105 \times 10^{-6}$		& $224$ 	& $2.8104 \times 10^{-6}$	& $234$  	\\ \hline
$8$	& $8.3276 \times 10^{-7}$	& $237$ 	& $1.5139 \times 10^{-6}$	& $253$  	\\ \hline
$9$	& $5.4493 \times 10^{-7}$	& $250$ 	& $6.9960 \times 10^{-7}$	& $261$  	\\ \hline
$10$	& $3.4339 \times 10^{-7}$	& $258$ 	& $7.5178 \times 10^{-7}$	& $268$  	\\ \hline
$11$	& $3.4145 \times 10^{-7}$		& $268$ 	& $2.7108 \times 10^{-7}$	& $281$  	\\ \hline
$12$	& $1.3458 \times 10^{-7}$		& $274$	& $5.5799 \times 10^{-7}$	& $287$  	\\ \hline
$13$	& $2.3812 \times 10^{-7}$		& $275$ 	& $2.3815 \times 10^{-7}$	& $295$  	\\ \hline
$14$	& $9.0332 \times 10^{-8}$	& $278$	& $9.0299 \times 10^{-8}$	& $297$  	\\ \hline
\end{tabular}
\normalsize
\end{table}

\begin{table}[htp]
\centering
\caption{Performance of $\CGPOPS$ on the Example 4 using $hp$-II(3,10). \label{tab:freeFlyingRobotMeshRefinement-hpDarby}}
\footnotesize
\begin{tabular}{|c|c|c|c|c|} \hline
{\bf Mesh} & {\bf Estimated}&{\bf Number of} &{\bf Estimated }&{\bf Number of} \\
{\bf Iteration} & {\bf Error ($\CGPOPS$)}&{\bf Collocation}&{\bf Error ($\mathbb{GPOPS-II}$)}&{\bf Collocation} \\
{\bf Number} & {\bf $hp$-II(3,10)}&{\bf Points}&{\bf $hp$-II(3,10)}&{\bf Points} \\\hline\hline
$1$	& $5.7636 \times 10^{-4}$	& $50$		&	$5.7636 \times 10^{-4}$	& $50$		\\ \hline
$2$	& $2.3718 \times 10^{-4}$		& $98$ 	& $1.1649 \times 10^{-4}$	& $98$ \\ \hline
$3$	& $6.4909 \times 10^{-5}$	& $162$	& $9.3164 \times 10^{-5}$	& $146$ \\ \hline
$4$	& $2.1470 \times 10^{-5}$		&$219$ 	&	$1.1244 \times 10^{-4}$	&$207$ \\ \hline
$5$	& $9.3539 \times 10^{-6}$	& $263$ 	& $3.2283 \times 10^{-6}$	& $267$ \\ \hline
$6$	& $1.0198 \times 10^{-6}$		& $297$ 	& $3.5320 \times 10^{-7}$	& $302$  \\ \hline
$7$	& $1.7028 \times 10^{-7}$		& $310$ 	& $2.3505 \times 10^{-7}$	& $320$  \\ \hline
$8$	& $9.8413 \times 10^{-8}$	&$315$ 	&	$1.3862 \times 10^{-7}$	&$322$ \\ \hline
$9$	& --										& --			& $1.0431 \times 10^{-7}$	& $325$ \\ \hline
$10$	& --										& -- 			& $9.5122 \times 10^{-8}$	& $328$  \\ \hline
\end{tabular}
\normalsize
\end{table}
     
\begin{table}[htp]
\centering
\caption{Performance of $\CGPOPS$ on Example 4 using $hp$-III(3,10). \label{tab:freeFlyingRobotMeshRefinement-hpLiu}}
\footnotesize
\begin{tabular}{|c|c|c|c|c|} \hline
{\bf Mesh} & {\bf Estimated}&{\bf Number of} &{\bf Estimated }&{\bf Number of} \\
{\bf Iteration} & {\bf Error ($\CGPOPS$)}&{\bf Collocation}&{\bf Error ($\mathbb{GPOPS-II}$)}&{\bf Collocation} \\
{\bf Number} & {\bf $hp$-III(3,10)}&{\bf Points}&{\bf $hp$-III(3,10)}&{\bf Points} \\\hline\hline
$1$	& $5.7636 \times 10^{-4}$	& $50$ 	&	$5.7636 \times 10^{-4}$	& $50$		\\ \hline
$2$	& $1.8489 \times 10^{-4}$	& $68$ 	&	$1.8489 \times 10^{-4}$ 						& $68$  \\ \hline
$3$	& $5.8497 \times 10^{-5}$	& $185$	&	$5.8497 \times 10^{-5}$	& $185$ \\ \hline
$4$	& $4.3708 \times 10^{-6}$	&$275$ 	&	$4.3709 \times 10^{-6}$	&$264$ \\ \hline
$5$	& $8.2894 \times 10^{-7}$	& $349$ 	& $2.3747 \times 10^{-6}$	& $324$  \\ \hline
$6$	& $4.5337 \times 10^{-7}$	& $395$ 	& $2.4780 \times 10^{-7}$	& $389$ \\ \hline
$7$	& $8.1069 \times 10^{-8}$	& $460$ 	& $1.5231 \times 10^{-7}$	& $410$ \\ \hline
$8$	& --										& -- 			& $1.0142 \times 10^{-7}$	& $436$ \\ \hline
$9$	& --										& -- 			& $2.1817 \times 10^{-7}$	& $437$ \\ \hline
$10$	& --										& -- 			& $8.0985 \times 10^{-8}$	& $458$ \\ \hline
\end{tabular}
\normalsize
\end{table}
      
\begin{table}[htp]
\centering
\caption{Performance of $\CGPOPS$ on Example 4 using $hp$-IV(3,10). \label{tab:freeFlyingRobotMeshRefinement-hpLegendre}}
\footnotesize
\begin{tabular}{|c|c|c|c|c|} \hline
{\bf Mesh} & {\bf Estimated}&{\bf Number of} &{\bf Estimated }&{\bf Number of} \\
{\bf Iteration} & {\bf Error ($\CGPOPS$)}&{\bf Collocation}&{\bf Error ($\mathbb{GPOPS-II}$)}&{\bf Collocation} \\
{\bf Number} & {\bf $hp$-IV(3,10)}&{\bf Points}&{\bf $hp$-IV(3,10)}&{\bf Points} \\\hline\hline
$1$	& $5.7636 \times 10^{-4}$	& $50$ 	&	$5.7636 \times 10^{-4}$	& $50$		\\ \hline
$2$	& $7.2614 \times 10^{-5}$ 	& $100$ &	$7.2614 \times 10^{-5}$ 						& $100$  \\ \hline
$3$	& $5.8350 \times 10^{-5}$	& $163$&	$5.8350 \times 10^{-5}$	& $163$ \\ \hline
$4$	& $7.0276 \times 10^{-6}$	&$212$ &	$3.9712 \times 10^{-6}$	&$203$ \\ \hline
$5$	& $2.9097 \times 10^{-6}$	& $259$ & $1.9372 \times 10^{-6}$	& $249$  \\ \hline
$6$	& $5.0338 \times 10^{-7}$	& $317$ & $7.0224 \times 10^{-6}$	& $301$  \\ \hline
$7$	& $2.1987 \times 10^{-7}$		& $362$ & $1.1880 \times 10^{-6}$	& $328$  \\ \hline
$8$	& $9.8979 \times 10^{-8}$	& $376$ & $7.4092 \times 10^{-7}$	& $347$  \\ \hline
$9$	& --										& -- 		& $1.9947 \times 10^{-7}$	& $360$  \\ \hline
$10$	& --										& -- 		& $9.1526 \times 10^{-8}$	& $373$  \\ \hline
\end{tabular}
\normalsize
\end{table}

\begin{table}[htp]
\centering
\caption{Performance of $\CGPOPS$ on Example 4 using $hp$-BB(3,10). \label{tab:freeFlyingRobotMeshRefinement-hpBangBang}}
\footnotesize
\begin{tabular}{|c|c|c|c|c|} \hline
{\bf Mesh} & {\bf Estimated}&{\bf Number of} &{\bf Estimated }&{\bf Number of} \\
{\bf Iteration} & {\bf Error ($\CGPOPS$)}&{\bf Collocation}&{\bf Error ($\mathbb{GPOPS-II}$)}&{\bf Collocation} \\
{\bf Number} & {\bf $hp$-BB(3,10)}&{\bf Points}&{\bf $hp$-II(3,10)}&{\bf Points} \\\hline\hline
$1$	& $5.7636 \times 10^{-4}$	& $50$		&	$5.7636 \times 10^{-4}$	& $50$		\\ \hline
$2$	& $6.2675 \times 10^{-9}$		& $108$ 	& $1.1649 \times 10^{-4}$	& $98$ \\ \hline
$3$	& --										& --			& $9.3164 \times 10^{-5}$	& $146$ \\ \hline
$4$	& --										& --			&	$1.1244 \times 10^{-4}$	&$207$ \\ \hline
$5$	& --										& --			& $3.2283 \times 10^{-6}$	& $267$ \\ \hline
$6$	& --										& --			& $3.5320 \times 10^{-7}$	& $302$  \\ \hline
$7$	& --										& --			& $2.3505 \times 10^{-7}$	& $320$  \\ \hline
$8$	& --										& --			&	$1.3862 \times 10^{-7}$	&$322$ \\ \hline
$9$	& --										& --			& $1.0431 \times 10^{-7}$	& $325$ \\ \hline
$10$	& --										& -- 			& $9.5122 \times 10^{-8}$	& $328$  \\ \hline
\end{tabular}
\normalsize
\end{table}

\clearpage

\subsection{Example 5: Multiple-Stage Launch Vehicle Ascent Problem\label{subsect:launch}}


The problem considered in this section is the ascent of a multiple-stage launch vehicle.  The objective is to maneuver the launch vehicle from the ground to the target orbit while maximizing the remaining fuel in the upper stage.  It is noted that this example is found verbatim in Refs.~\cite{Benson1}, \cite{Rao8}, and \cite{Betts3}.  The problem is modeled using four phases where the objective is to maximize the mass at the end of the fourth phase, that is maximize
\begin{equation}\label{eq:launch-cost}
  \cal{J} = m(t_f^{(4)})~,
\end{equation}
subject to the dynamic constraints 
\begin{equation}\label{eq:launch-dynamics}
\begin{array}{rcl}
  \dot{\textbf{r}}^{(p)} &=& \mathbf{v}^{(p)}~, \\
  \dot{\textbf{v}}^{(p)} &=& -\displaystyle\frac{\mu}{\|\textbf{r}^{(p)}\|^3}\mathbf{r}^{(p)} +
  \displaystyle\frac{T^{(p)}}{m^{(p)}}\mathbf{u}^{(p)} + \displaystyle\frac{\mathbf{D}^{(p)}}{m ^{(p)}}~, \\
  \dot{m}^{(p)} & = & -\displaystyle\frac{T ^{(p)}}{g_0I_{sp}}~, \\
\end{array}\quad (p=1,2,3,4)~,
\end{equation}
the initial conditions
\begin{equation}\label{eq:launch-ICs}
\begin{array}{rcl}
\mathbf{r}(t_0) &=& \mathbf{r}_0 = (5605.2,0,3043.4)\times 10^{3}\textrm{ m}~, \\
\mathbf{v}(t_0) &=& \mathbf{v}_0 = (0,0.4076,0)\times 10^{3}\textrm{ m/s}~, \\
m(t_0) &=& m_0 = 301454 \textrm { kg}~.
\end{array}
\end{equation}
the interior point constraints
\begin{equation}\label{eq:launch-IntPnt}
\begin{array}{rcl}
\mathbf{r}^{(p)}(t_f^{(p)})-\mathbf{r}^{(p+1)}(t_0^{(p+1}) &=& \mathbf{0}~, \\
\mathbf{v}^{(p)}(t_f^{(p)})-\mathbf{v}^{(p+1)}(t_0^{(p+1)}) &=& \mathbf{0}~, \quad (p=1,2,3)~, \\
m^{(p)}(t_f^{(p)})-m_{\textrm{dry}}^{(p)}-m^{(p+1)}(t_0^{(p+1)}) &=& 0~, \\
\end{array}
\end{equation}
the terminal constraints (corresponding to a geosynchronous transfer orbit)
\begin{equation}\label{eq:launch-FCs}
\begin{array}{lclclcl}
 a(t_f^{(4)}) & = & a_f  = 24361.14 \textrm{ km}~, & & 
 e(t_f^{(4)})  & = & e_f  =  0.7308~, \\
 i(t_f^{(4)})  & = & i_f  = 28.5 \deg~, & & 
 \theta(t_f^{(4)}) & = & \theta_f = 269.8\deg~, \\
 \phi(t_f^{(4)}) & = & \phi_f =  130.5\deg~,
\end{array}
\end{equation}
and the path constraints
\begin{equation}\label{eq:launch-xpath}
\begin{array}{lcl}
  |\mathbf{r}^{(p)}|_2 & \geq & R_e~, \\
  \|\mathbf{u}^{(p)}\|_2^2 & = & 1~, 
\end{array}\quad (p=1,\ldots,4)~.
\end{equation}
In each phase $\mathbf{r}=(x,y,z)$ is the position relative to the center of the Earth expressed in earth-centered inertial (ECI) coordinates, $\mathbf{v} = (v_x,v_y,v_z)$ is the inertial velocity expressed in ECI coordinates, $\mu$ is the gravitational parameter, $T$ is the vacuum thrust, $m$ is the mass, $g_0$ is the acceleration due to gravity at sea level, $I_{sp}$ is the specific impulse of the engine, $\mathbf{u} = (u_x,u_y,u_z)$ is the  thrust direction expressed in ECI coordinates, and $\mathbf{D}=(D_x,D_y,D_z)$ is the drag force expressed in ECI coordinates.  The drag force is defined as
\begin{equation}\label{eq:launch-drag-def}
  \m{D} = -{\textstyle \frac{1}{2}}C_D S\rho \|\mathbf{v}_{\textrm{rel}}\|\mathbf{v}_{\textrm{rel}}~,
\end{equation}
where $C_D$ is the drag coefficient, $S$ is the vehicle reference area, $\rho=\rho_0\exp(-h/H)$ is the atmospheric density, $\rho_0$ is the sea level density, $h=r-R_e$ is the altitude, $r=\|\m{r}\|_2=\sqrt{x^2+y^2+z^2}$ is the geocentric radius, $R_e$ is the equatorial radius of the Earth, $H$ is the density scale height, $\mathbf{v}_{\textrm{rel}}=\mathbf{v}-\boldsymbol{\omega} \times \mathbf{r}$ is the velocity as viewed by an observer fixed to the Earth expressed in ECI coordinates, and $\boldsymbol\omega=(0,0,\Omega)$ is the angular velocity of the Earth as viewed by an observer in the inertial reference frame expressed in ECI coordinates.  Furthermore, $m_{\textrm{dry}}$ is the dry mass of phases 1, 2, and 3 and is defined $m_{\textrm{dry}}=m_{\textrm{tot}}-m_{\textrm{prop}}$, where $m_{\textrm{tot}}$ and $m_{\textrm{prop}}$ are, respectively, the total mass and dry mass of phases 1, 2, and 3.  Finally, the quantities $a$, $e$, $i$, $\theta$, and $\phi$ are, respectively, the semi-major axis, eccentricity, inclination, longitude of ascending node, and argument of periapsis, respectively.  The vehicle data for this problem and the numerical values for the physical constants can be found in Tables~\ref{tab: launch vehicle properties} and \ref{tab:dynamics properties}, respectively.

\begin{table}[htp]
\centering
\caption{Vehicle Properties for Multiple-Stage Launch Vehicle Ascent Problem. \label{tab: launch vehicle properties}}
\footnotesize
\begin{tabular}{|c|c|c|c|}
\hline
Quantity & Solid Boosters & Stage 1 & Stage 2 \\
 \hline \hline
 $m_{\textrm{tot}}$ (kg) & 19290 & 104380 & 19300  \\
 \hline
 $m_{\textrm{prop}}$ (kg) & 17010 & 95550 & 16820 \\
 \hline
 $T$  (N) & 628500 & 1083100 & 110094 \\
 \hline
 $I_{sp}$ (s) & 283.3 & 301.7 & 467.2 \\
 \hline
 Number of Engines & 9 & 1 & 1 \\
 \hline
 Burn Time (s) & 75.2 & 261 & 700 \\
 \hline
\end{tabular}
\normalsize
\end{table}

\begin{table}[htp]
\centering
\caption{Constants Used in the Launch Vehicle Ascent Optimal Control Problem. \label{tab:dynamics properties}}
\footnotesize
\begin{tabular}{|c|c|}
\hline
Constant & Value \\
\hline \hline
Payload Mass & $4164$ kg \\
\hline
$S$ & $4\pi$ m${}^2$ \\
\hline
$C_D$ & $0.5$ \\
\hline
$\rho_0$ & $1.225$ kg/m${}^3$ \\
\hline
$H$ & $7200$ m\\
\hline
$t_1$ & $75.2$ s\\
\hline
 $t_2$ & $150.4$ s \\
\hline
 $t_3$ & $261$ s \\
\hline
 $R_e$ & $6378145$ m \\
\hline
 $\Omega$ & $7.29211585\times10^{-5}$ rad/s \\
\hline
 $\mu$ & $3.986012\times10^{14}$ m${}^3$/s${}^2$ \\
\hline
$g_0$ & $9.80665$ m/s${}^2$ \\ \hline
\end{tabular}
\normalsize
\end{table}

\begin{figure}[htp]
\centering
\subfloat[$h(t)$ vs.~$t$.\label{fig:launchAltitude}]{\epsfig{figure=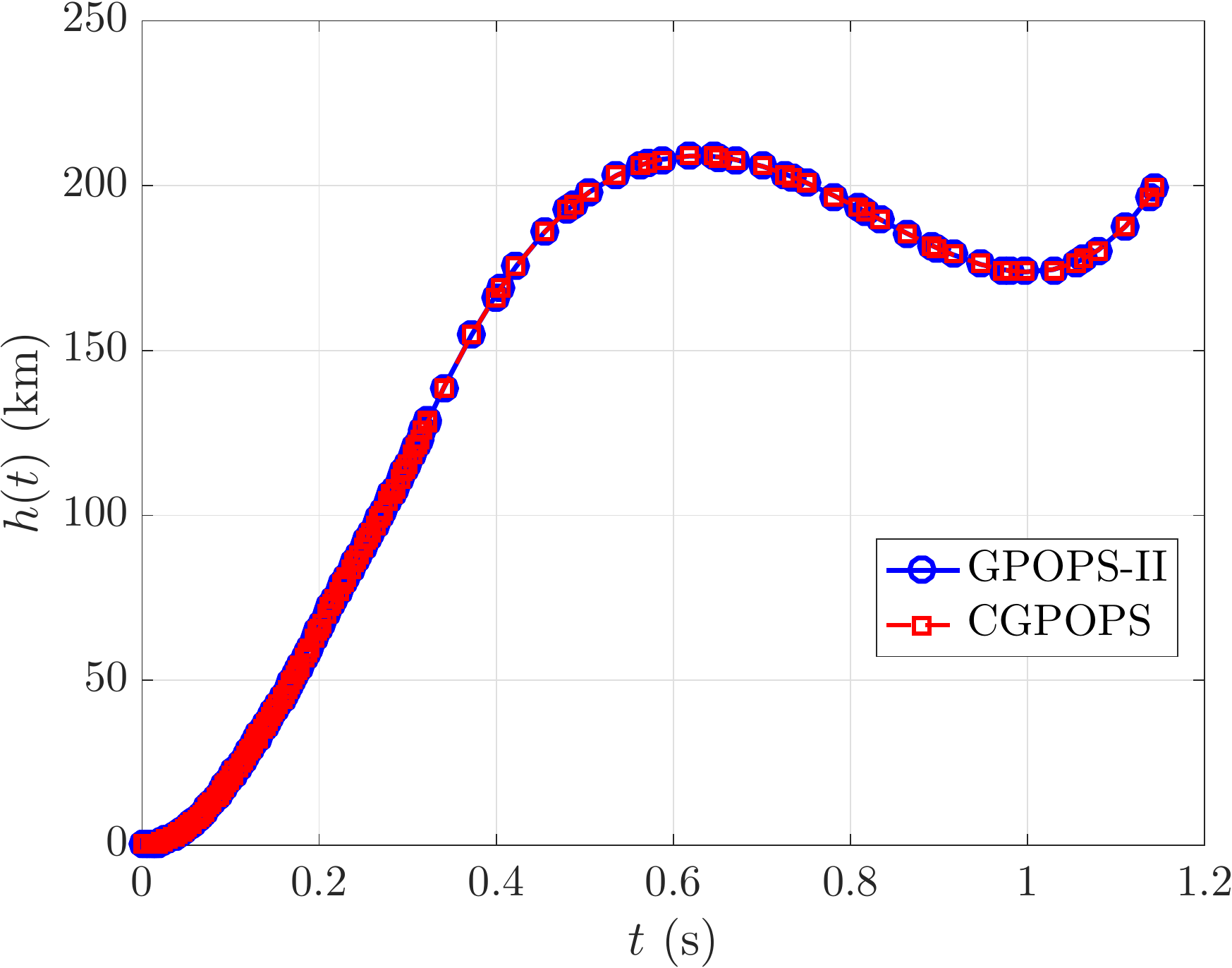,height=2in}}~~~\subfloat[$v(t)$ vs.~$t$. \label{fig:launchSpeed}]{\epsfig{figure=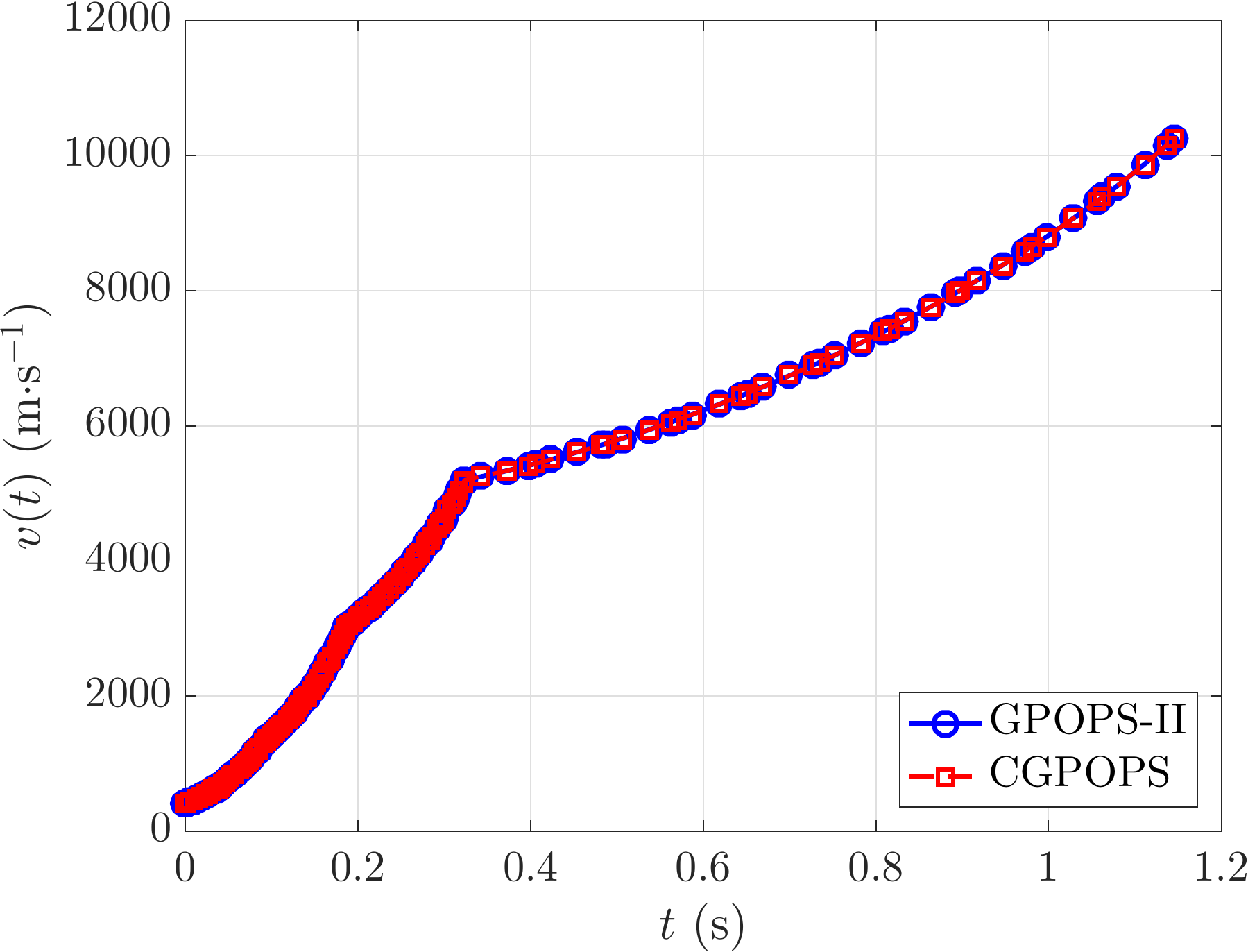,height=2in}}

\subfloat[$m(t)$ vs.~$t$.\label{fig:launchMass}]{\epsfig{figure=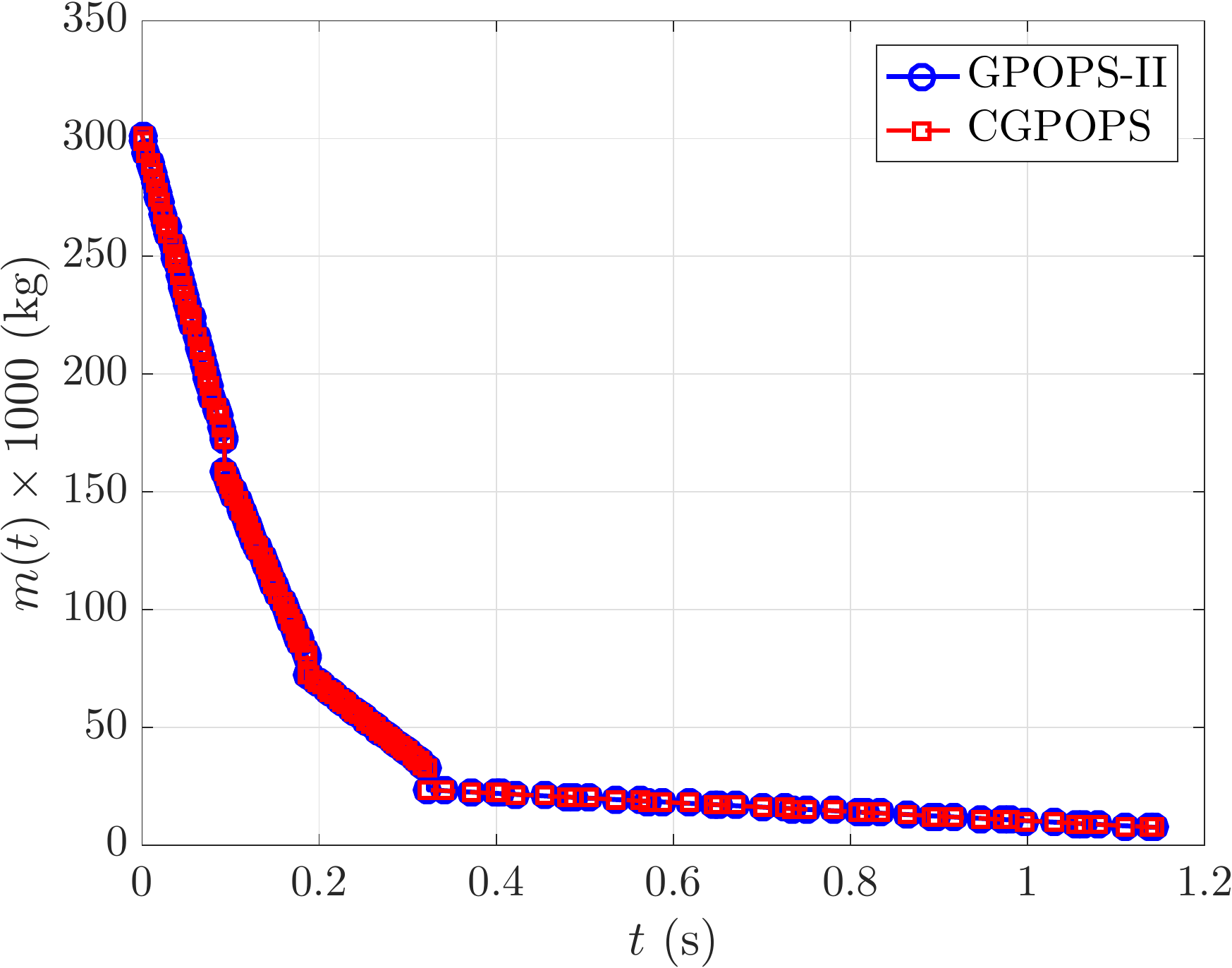,height=2in}}~~~\subfloat[$u_1(t)$ vs.~$t$.\label{fig:launchControl1}]{\epsfig{figure=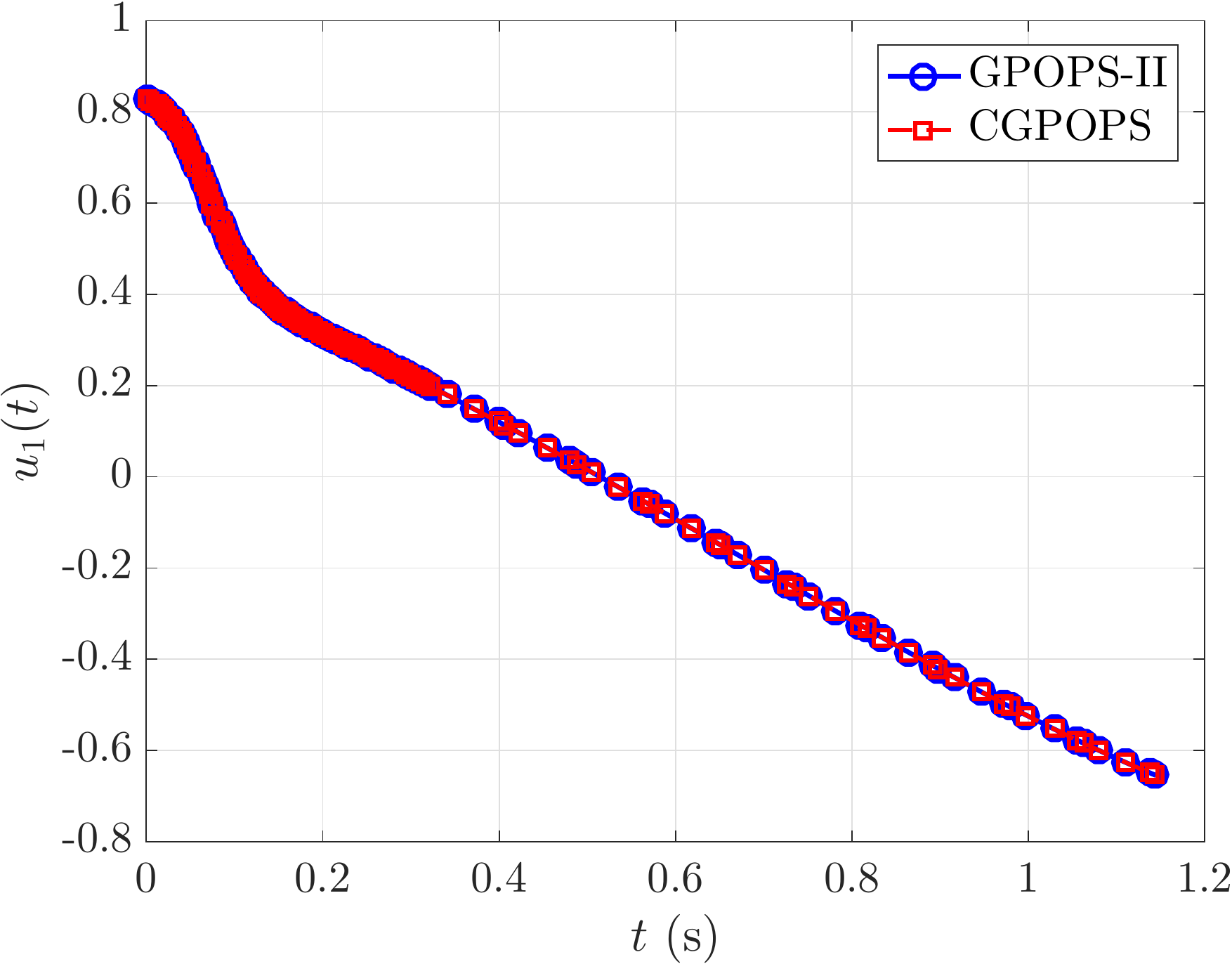,height=2in}}

\subfloat[$u_2(t)$ vs.~$t$.\label{fig:launchControl2}]{\epsfig{figure=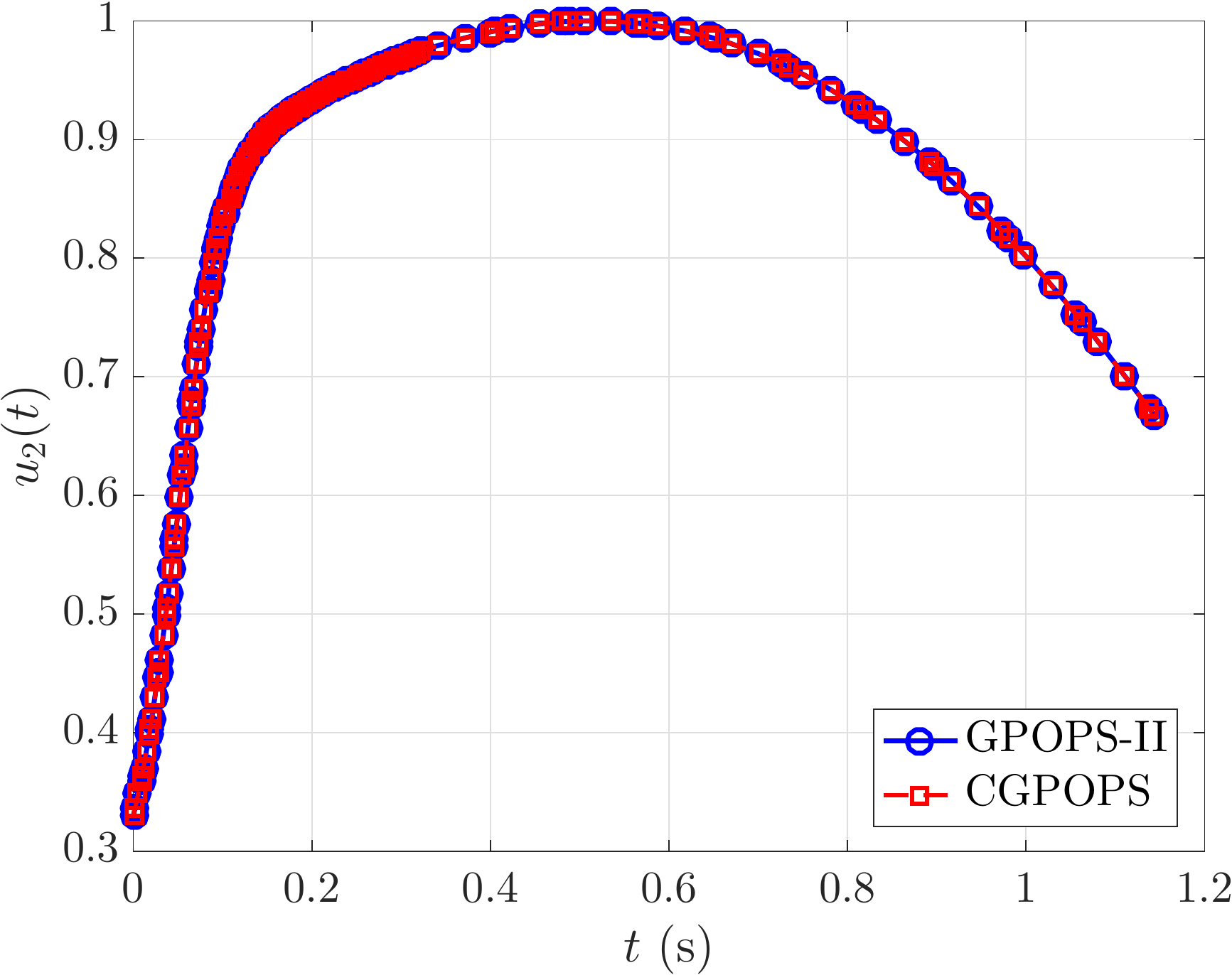,height=2in}}~~~\subfloat[$u_3(t)$ vs.~$t$.\label{fig:launchControl3}]{\epsfig{figure=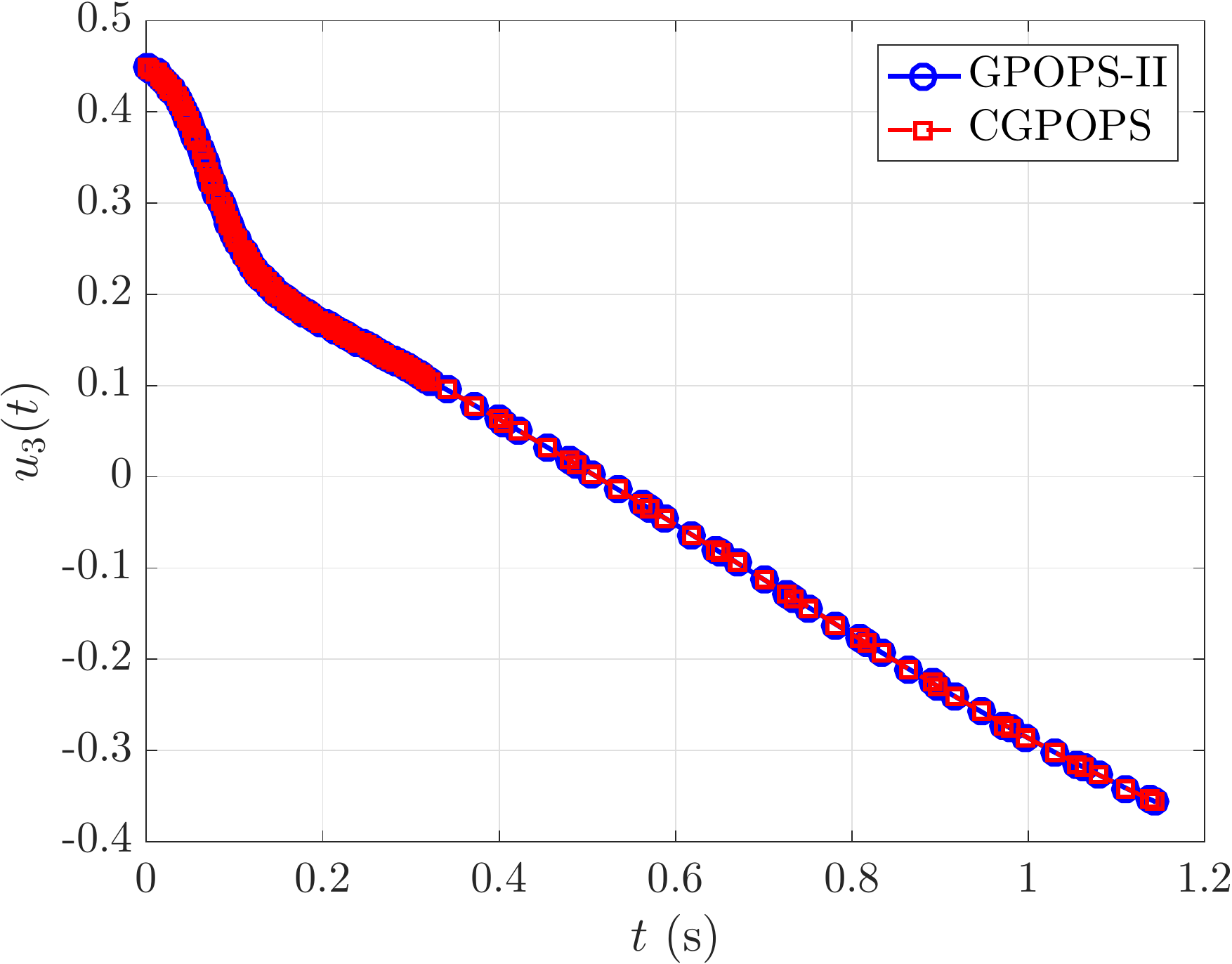,height=2in}}

\caption{Solution of Example 5 Using $\CGPOPS$. \label{fig:launchSolution}}
\end{figure}

The multiple-stage launch vehicle ascent optimal control problem was solved using $\CGPOPS$ with an initial mesh in each phase consisting of ten uniformly spaced mesh intervals with four LGR points per mesh interval.  The NLP solver and mesh refinement accuracy tolerances were set to $10^{-7}$ and $10^{-6}$, respectively.  The initial guess of the solution was constructed such that the initial guess of the position and the velocity in phases 1 and 2 was constant at $(\m{r}(0),\m{v}(0))$ as given in Eq.~\eqref{eq:launch-ICs} while in phases 3 and 4 the initial guess of the position and velocity was constant at $(\tilde{\m{r}},\tilde{\m{v}})$, where $(\tilde{\m{r}},\tilde{\m{v}})$ are obtained via a transformation from orbital elements to ECI coordinates using the five known orbital elements of Eq.~\eqref{eq:launch-FCs} and a true anomaly of zero.  Furthermore, in all phases the initial guess of the mass was a straight line between the initial and final mass, $m(t_0^{(p)})$ and $m(t_f^{(p)})$ ($p\in[1,\ldots,4]$).  Moreover, in all phases the guess of the control was constant at $\m{u}=(0,1,0)$.  The $\CGPOPS$ solution is shown in Fig.~\ref{fig:launchSolution}.  In this example the mesh refinement accuracy tolerance of $10^{-6}$ is satisfied on the initial mesh using both $\CGPOPS$ and $\mathbb{GPOPS-II}$, so no mesh refinement is necessary.  The solution obtained using $\CGPOPS$ matches closely with the solution obtained using the software $\mathbb{GPOPS-II}$ \cite{Patterson2014}, where it is noted that the optimal objective values obtained using $\CGPOPS$ and $\mathbb{GPOPS-II}$ are $7547.9729$ and $7547.9739$, respectively.  Finally, the computational time required by $\CGPOPS$ and $\mathbb{GPOPS-II}$ to solve the optimal control problem was $2.9466$ seconds and $18.9401$ seconds, respectively.

\clearpage



\section{Capabilities of $\CGPOPS$\label{sect:capabilities}}


The various capabilities of $\CGPOPS$ were highlighted by the five examples provided in Section~\ref{sect:examples}.  First, the capabilities of the $hp$ mesh refinement methods were demonstrated on the hyper-sensitive example where the the number of mesh intervals and collocation points were increased in regions where the dynamics changed rapidly.  Additionally, the ability of the $hp$ methods to satisfy a specified accuracy tolerance using fewer collocation points than might be possible using a fixed low-order collocation method was indicated by the reusable launch vehicle entry problem.  Second, the flexibility of the software to achieve better performance by modifying the default settings of the mesh initialization and/or refinement process is demonstrated by the space station attitude control and free-flying robot examples.  Third, all five examples demonstrate the increased computational efficiency of implementing the optimal control framework developed in Sections~\ref{sect:multi-phase} to \ref{sect:RPM} in C++ as opposed to MATLAB.  In particular, the space station attitude control example shows the computational benefits of using an exact NLP Lagrangian Hessian sparsity pattern obtained by the derivative dependency identification capabilities of the hyper-dual (or bicomplex-step) derivative supplier as described in Section~\ref{subsect:dependencies}.  Furthermore, the free-flying robot example demonstrates how the $hp$-BB mesh refinement method implemented in $\CGPOPS$ is significantly more effective for solving bang-bang optimal control problems compared to previously developed $hp$ methods.  The different examples also demonstrate the wide variety of problems that can be formulated in $\CGPOPS$.  Such problems range from one-phase problems with a Lagrange cost (for example, the hyper-sensitive, space station attitude control, and free-flying robot examples), Mayer cost (for example, the reusable launch vehicle entry and launch vehicle ascent examples), and problems that include bang-bang controls (for example, the free-flying robot example).  Moreover, it was shown that $\CGPOPS$ has the ability to solve challenging multiple-phase problems (for example, the launch vehicle ascent example).  The fact that $\CGPOPS$ is capable of solving the challenging benchmark optimal control problems shown in this paper shows the general utility of the software on problems that may arise in different application areas.



\section{Limitations of $\CGPOPS$\label{sect:limitations}}


As with any software, $\CGPOPS$ has limitations.  First, it is assumed that all functions used to formulate an optimal control problem of interest have continuous first and second derivatives.  It is noted, however, that for some applications the functions may have discontinuous derivatives while the functions themselves are continuous.  In cases where the derivatives are discontinuous $\CGPOPS$ may have difficulty obtaining a solution because the NLP solver operates under the assumption that all first and second derivatives are continuous.  Second, because $\CGPOPS$ is a direct collocation method, the ability to obtain a solution depends upon the NLP solver that is used.  In particular, while the NLP solver IPOPT \cite{Biegler2} used with $\CGPOPS$ may be effective for some examples, other NLP solvers (for example, SNOPT \cite{Gill1} or KNITRO \cite{Byrd1}) may be more effective than IPOPT for certain problems.  Moreover, for problems with high-index path constraints, the constraint qualification conditions may not be satisfied when the mesh becomes extremely fine.  In such cases, unique NLP Lagrange multipliers may not exist or, in some cases, these Lagrange multipliers may be unbounded.  Furthermore, it may be difficult to obtain a solution to a poorly scaled problem.  Finally, as is true for any optimal control software, optimal control problems whose solutions lie on a singular arc can create problems due to the inability to determine the optimal control along the singular arc.  Moreover, the problems associated with a singular optimal control problem are exacerbated with mesh refinement.  Thus, when solving a singular optimal control problem, it may be necessary to modify the original problem by including the higher-order optimality conditions that define the control on the singular arc.



\section{Conclusions\label{sect:conclusions}}


A general-purpose C++ software program called $\CGPOPS$  has been described for solving multiple-phase optimal control problems using adaptive direct orthogonal collocation methods.  In particular, the software employs a Legendre-Gauss-Radau quadrature orthogonal collocation where the continuous control problem is transcribed to a large sparse nonlinear programming problem.  The software implements five previously developed adaptive mesh refinement methods that allow for flexibility in the number and placement of the collocation and mesh points in order to achieve a specified accuracy.  In addition, the software is designed to compute all derivatives required by the NLP solver using one of four derivative estimation methods for the optimal control functions.  The key components of the software have been described in detail and the utility of the software is demonstrated on five benchmark optimal control problems.  The software described in this paper provides researchers a transitional platform upon which to solve a wide variety of complex constrained optimal control problems for real-time applications. 


\section*{Acknowledgments}

The authors gratefully acknowledge support for this research from the U.S.~Office of Naval Research under grant N00014-15-1-2048 and from the U.S.~National Science Foundation under grants DMS-1522629, DMS-1924762, and CMMI-1563225.  

\renewcommand{\baselinestretch}{1.0}
\normalsize\normalfont
\bibliographystyle{aiaa}

\end{document}